\documentclass[leqno,11pt]{amsart}
\usepackage{latexsym,amsfonts,amsmath,epic,eepic,psfig,epsf}

\newtheorem{theorem}{Theorem}[section]

\newtheorem{lemma}[theorem]{Lemma}
\newtheorem{proposition}[theorem]{Proposition}
\theoremstyle{definition}
\newtheorem{definition}{Definition}[section]
\newtheorem{example}{Example}
\theoremstyle{remark}
\newtheorem{remark}{Remark}[section]
\numberwithin{equation}{section}
\newcommand{\labelm}[1]{\label{#1}}
\newcommand{\C}{{\mathbb C}} \newcommand{\R}{{\mathbb R}}
\newcommand{\Z}{{\mathbb Z}} \newcommand{\N}{{\mathbb N}}
\newcommand{\F}{{\mathbb F}} \newcommand{\CA}{\mathcal{A}}
\newcommand{\CB}{\mathcal{B}} \newcommand{\ba}{\mathbf{a}}
 
\newcommand{\vol}{\operatorname{vol}} \newcommand{\Sm}{\mathcal{S}}
\newcommand{\Smm}{\mathcal{P}} \newcommand{\SP}{\mathrm{SP}}
\newcommand{\RSP}{\mathrm{RSP}}
\newcommand{\Coeff}{\operatorname{Coeff}}
\newcommand{\Tres}{\mathrm{Tres}}
\newcommand{\vsing}{C^\Delta_\mathrm{sing}}
\newcommand{\vreg}{C^\Delta_\mathrm{reg}}
\newcommand{\cnull}{\c^{\mathrm{null}}}
\newcommand{\bdelta}{\mathcal{B}(\Delta)}
\newcommand{\bfi}{\mathcal{B}(\Phi)} \renewcommand{\c}{\mathfrak{c}}
 \newcommand{\jk}[3]{J\langle #1,
  #2\rangle_{#3}} \newcommand{\JK}[3]{J\left\langle #1,
    #2\right\rangle_{#3}}
\newcommand{\res}{\operatornamewithlimits{Res}}
 \newcommand{\vc}{V_\C}
\newcommand{\bfj}{\mathbf{j}}
\newcommand{\volga}{\mathrm{vol}_{\Gamma^*}}
\newcommand{\hatr}{\widehat{R}_\mathrm{hp}}
\newcommand{\bh}{\mathbf{h}}\newcommand{\by}{\mathbf{y}}\newcommand{\bx}{\mathbf{x}}
\newcommand{\even}{\mathrm{even}}\newcommand{\odd}{\mathrm{odd}}
\newcommand{\EP}{\mathrm{EP}} \newcommand{\bfk}{\mathbf{k}}
\title[Residue formulae for vector partitions]{Residue formulae for
  vector partitions and Euler-MacLaurin sums}
\author{Andr{\'a}s Szenes and Mich{\`e}le Vergne} \date{\today}
\begin{document}
\maketitle

\setcounter{section}{-1}
\section{Introduction}
\labelm{sec:intro}

Let $V$ be an $n$-dimensional real vector space endowed with a
rank-$n$ lattice $\Gamma$.  The dual lattice $\Gamma^*$ is
naturally a subset of the dual vector space $V^*$. Let
$\Phi=[\beta_1, \beta_2,\ldots,\beta_N]$ be a sequence of not
necessarily distinct elements of $\Gamma^*$, which span $V^*$ and
lie entirely in an open halfspace of $V^*$.  In what follows, the
order of elements in the sequence will not matter.

The closed cone $C(\Phi)$ generated by the elements of $\Phi$ is
an acute convex cone, divided into open conic chambers by the
$(n-1)$-dimensional cones generated by linearly independent
$(n-1)$-tuples of elements of $\Phi$.  Denote by $\Z \Phi$ the
sublattice of $\Gamma^*$ generated by $\Phi$. Pick a vector $a\in
V^*$ in the cone $C(\Phi)$, and denote by $\Pi_\Phi(a)\subset
\R_+^N$ the convex polytope consisting of all solutions
$\bx=(x_1,x_2, \ldots, x_N)$ of the equation $\sum_{k=1}^N x_k
\beta_k=a$ in {\em nonnegative real} numbers $x_k$. This is a
closed convex polytope  called the {\em partition polytope}
associated to $\Phi$ and $a$. Conversely, any  closed convex
polytope can be realized as a partition polytope.

If $\lambda\in \Gamma^*$, then the vertices of the partition polytope
$\Pi_\Phi(\lambda)$ have rational coordinates.  We denote by
$\iota_\Phi(\lambda)$ the number of points with integral coordinates
in $\Pi_\Phi(\lambda)$. Thus $\iota_\Phi(\lambda)$ is the number of
solutions of the equation $\sum_{k=1}^N x_k \beta_k=\lambda$ in {\em
  nonnegative integers} $x_k$. The function $\lambda\mapsto
\iota_\Phi(\lambda)$ is called the {\em vector partition function}
associated to $\Phi$. Obviously, $\iota_\Phi(\lambda)$ vanishes if
$\lambda$ does not belong to $C(\Phi)\cap \Z\Phi$.

Let $\EP(\R^N)$ be  the ring of complex functions on $\R^N$ generated
by exponentials and polynomials. Thus any  $f\in  \EP(\R^N)$ is of the
form
$$f(\bx)=\sum_{j=1}^m e^{\langle \by_j,\bx\rangle} P_j(\bx),$$
where $\by_1,\ldots, \by_m\in \C^N$, and the functions $P_1,\ldots,
P_m$ are polynomials with complex coefficients. If the elements $\{\by_j\,,1\leq
j\leq m\}$ are such that there exists an integer $M$ with
$M\by_j\in 2\pi i\Z^N$, then the function $f$ is said to be {\em
periodic-polynomial}. The restriction of such a function to any
coset $\bx+ M\Z^N$ of $\Z^N/M\Z^N$ is plainly polynomial.

A generalization of $\iota_\Phi(\lambda)$ is the sum of the values
of a function  $f\in \EP(\R^N)$ over the integral points of
$\Pi_\Phi(\lambda)$:
$$\Sm[f,\Phi](\lambda)=\sum_{\xi\in \Pi_{\Phi}(\lambda)\cap \Z^N} f(\xi).$$
Indeed, if $f=1$, the function $\Sm[f,\Phi]$ is just the function
$\iota_\Phi$. Such a sum $\Sm[f,\Phi]$ will be called an
Euler-MacLaurin sum.

In this paper, we will search for ``explicit'' formulae for the
function $\lambda\mapsto \Sm[f,\Phi](\lambda)$ on $\Gamma^*$.  Let
us recall some qualitative results about this function. We start
with the following result of Ehrhart: for a rational polytope
$\Pi$ in $\R^r$, consider the function $k\mapsto\#(k\Pi\cap
\Z^r)$, where $\#S$ stands for the cardinality of the set $S$.
Ehrhart proved that this function is given by a
periodic-polynomial formula for all $k\geq 0$. More precisely (see
\cite{eh3} and references therein), if $M$ is an integer such that
all the vertices of the polytope $M\Pi$ are in $\Z^r$, then there
exist polynomial functions $P_j, 0\leq j\leq M-1$, such that
$\#(k\Pi\cap \Z^r)=\sum_{j=0}^{M-1}e^{\frac{2i\pi j}{M}
  k}P_j(k) $.  If $f$ is a polynomial, then $\Sm[f,\Phi](\lambda)$
consists of summing up the values of a polynomial over the
integral points of the rational polytope $\Pi_{\Phi}(\lambda)$. If
$f$ is an exponential $\bx\mapsto e^{\langle \by,\bx\rangle}$,
then $\Sm[e^\by,\Phi](\lambda)$ is the sum $\sum_{\xi\in
  \Pi_{\Phi}(\lambda)\cap \Z^N} e^{\langle \by,\xi\rangle}$;  such sums
were evaluated ``explicitly'' by M. Brion \cite{b} and by A.I.
Barvinok \cite{ba} for generic exponentials.

Assume first that $\Phi$ consists of $n=\dim V$ linearly
independent vectors of $\Gamma^*$. Denote by $\rho$ the linear
isomorphism from $\R^n$ to $V^*$ defined by $\rho(\bx)=\sum_{i=1}^n
x_i\beta_i$.  The set $\Pi_\Phi(\lambda)$ is nonempty if and only
if $\lambda\in C(\Phi)\cap \Z\Phi$. In this case, the set
$\Pi_\Phi(\lambda)$ coincides with $\rho^{-1}(\lambda)$, and our
function $\lambda\mapsto \Sm[f,\Phi](\lambda)$ on $\Gamma^*$ is
just the function $\lambda\mapsto f(\rho^{-1}(\lambda))$
restricted to $C(\Phi)\cap \Z\Phi$. In general, the map $\rho:
\R^N\rightarrow V^*$ defined by $\rho(\bx)=\sum_{i=1}^N x_i\beta_i$
is a surjection, and the following qualitative statement holds:
\begin{theorem}
For each conic chamber $\c$ of the cone  $C(\Phi)$, there exists
an exponential-polynomial function $\Smm[\c,f,\Phi]$ on $V^*$ such
that for each $\lambda\in \overline{\c}\cap \Gamma^*$, we have
$$\Sm[f,\Phi](\lambda)= \Smm[\c,f,\Phi](\lambda).$$
\end{theorem}

This theorem follows, for example, from \cite{bv2}, and there are
many antecedents of this result in particular cases.  The
periodic-polynomial behavior of $\iota_{\Phi}(\lambda)$ on
closures of conic chambers of the cone $C(\Phi)$ is proved in
\cite{s}. If $f$ is a polynomial function, then the sum $\sum_{\xi
\in k\Pi}f(\xi)$ is a polynomial function of $k$ for $k\geq 0$ if
the vertices of $\Pi$ have integral coordinates \cite{eh3,b,cs2}.
Let $\Pi_1, \Pi_2,\ldots,\Pi_N$ be rational polytopes in $\R^r$. For
a sequence $[k_1,\ldots,k_N]$ of nonnegative integers, denote by
$k_1\Pi_1+k_2\Pi_2+\cdots+k_N\Pi_N$ the weighted Minkowski sum of
the polytopes $\Pi_i$. Then, as proved in \cite{mm2}, there exists
an periodic-polynomial function $\Smm$ on $\R^N$ such that
$$
\#((k_1\Pi_1+k_2\Pi_2+\cdots+k_N\Pi_N)\cap
\Z^r)=\Smm(k_1,k_2,\ldots,k_N).$$ We explain in Section
\ref{sec:mink} how to pass from the setting of Minkowski sums to
the setting of partition polytopes.

Most of the investigations of the function $\Sm[f,\Phi]$
(\cite{kp,cs2,bv2}), starting with the Euler-MacLaurin formula
evaluating the sum $\sum_{A}^B f(k)$ of the values of a function
$f$ at all integral points of an interval $[A,B]$, were dedicated
to the fascinating relation of $\Sm[f,\Phi](\lambda)$ with the
integral of $f$ on the polytopes $\Pi_{\Phi}(a)$, when $a$ varies
near $\lambda$. This relation uses Todd differential operators,
which leads to a Riemann-Roch calculus for $\Sm[f,\Phi]$
initiated by Khovanskii-Pukhlikov \cite{kp}. In fact, there is a
dictionary between rational polytopes and line bundles on toric
varieties, which inspired these results.

Introduce the convex polytope
$$\Box(\Phi)=\sum_{i=1}^N [0,1]\beta_i.$$
We obtain a residue formula for $\Sm[f,\Phi]$ which implies  the
following qualitative result.
\begin{theorem}
For each conic chamber $\c$ of the cone  $C(\Phi)$, there exists
an exponential-polynomial function $\Smm[\c,f,\Phi]$ on $V^*$ such
that, for each $\lambda\in (\c-\Box(\Phi))\cap \Gamma^*$, we have
$$\Sm[f,\Phi](\lambda)= \Smm[\c,f,\Phi](\lambda).$$
\end{theorem}

 We assumed that $\Phi$
linearly generates $V^*$, hence the set $\c-\Box(\Phi)$ contains
$\overline\c$. In particular, the function $\iota_\Phi(\lambda)$
is periodic-polynomial on the neighborhood $\c-\Box(\Phi)$ of the
closure of the conic chamber $\c$ (this neighborhood is usually
much larger than $\overline \c$, see the pictures in the
Appendix). We give specific residue formulae on each of these
sectors $\c-\Box(\Phi)$.  Our main theorems are Theorem
\ref{thm:main} and its various corollaries: the residue formulae
of Theorem \ref{thm:vpartition} for $\iota_{\Phi}(\lambda)$ and
the residue formulae of Theorem \ref{thm:final} for
$\Sm[f,\Phi](\lambda)$. If $f$ is a generic exponential
$\bx\mapsto e^{\langle\by,\bx \rangle}$, then the residue formula
of Theorem \ref{exp} implies that  Formula 3.4.1 of Brion-Vergne
\cite{bv2} holds on the neighborhood $\c-\Box(\Phi)$ of
$\overline{\c}$.

 The residue formula  makes the exponential-polynomial
behavior of $\Sm[f,\Phi](\lambda)$ in each of these sectors clear.
More specifically, in  section \ref{sec:residue}, we construct an
exponential-polynomial function $E[f,\Phi]$ on the entire vector
space $V^*$ with values in a finite dimensional vector space $S$,
the space of {\em simple elements}, and linear functionals
$J_{\c}:S\rightarrow\C$ depending on the conic chamber $\c$, such
that $\Sm[f,\Phi](\lambda)=\langle J_{\c},
E[f,\Phi](\lambda)\rangle$ for $\lambda$  in a specified
neighborhood of the chamber $\c$ depending on $f$ and containing
$\c-\Box(\Phi)$.
 Moreover, from
the comparison with the Jeffrey-Kirwan expression for the volume
of $\Pi_{\Phi}(a)$, which is given by a  very similar residue
formula on each conic chamber (cf. \cite{bav}), one immediately
obtains the Riemann-Roch formula of \cite{kp,cs2,bv2} for
$\Sm[f,\Phi]$.

Conversely, applying Todd operators to the Jeffrey-Kirwan residue
expression, we could deduce our main theorem from \cite{cs2} or
\cite{bv2}. However, our present result is an explicit formula
which holds on a region larger than $\overline{\c}$, and the path
followed in the present article to obtain this result is direct.
Furthermore, our result has the advantage that it provides
independent and very similar residue formulae for volumes and for
Ehrhart polynomials of polytopes. These computations are quite
efficient: we give a few illustrative examples in the Appendix.
We refer to \cite{bav} for examples
of calculations of volumes by residue methods and examples of
application of change of variables in residue for expressions of
Ehrhart polynomials.

Our method is based on a detailed study of the generating function
$$\frac{1}{\prod_{k=1}^N(1-e^{\beta_k})}$$
for the partition function or, more generally, of periodic
meromorphic functions with poles on an affine arrangement of
hyperplanes.  As a main tool, we will use a separation theorem due to
the first author \cite{sz2}. We review these results in Section 1.

As stated before, the equation $\Sm[f,\Phi](\lambda)=
\Smm[\c,f,\Phi](\lambda)$ holds for $\lambda$ belonging to a
specified ``neighborhood'' of $\c$, which, in general, is strictly
larger than $\overline{\c}$. This neighborhood depends on $f$ and
$\Phi$. As a result the polynomials $\Smm[\c,f,\Phi](\lambda)$ for
two neighboring chambers will coincide along a thick strip near
their common boundary.  We illustrate our residue formula and
this effect with an example here.

\begin{example}
We set $V^*=\R^2$ with standard basis vectors $e_1,e_2$ and
corresponding coordinates $a_1,a_2$. Let $$\Phi_h=[e_1,e_1,\ldots,
e_1,e_2,e_2,\ldots,e_2,e_1+e_2,e_1+e_2,\ldots,e_1+e_2],$$ where
each vector $e_1$, $e_2$, $e_1+e_2$ is repeated $h$-time. There
are two chambers contained in $C(\Phi_h)$: $\c_1=\{a_1>a_2>0\}$ and
$\c_2=\{a_2>a_1>0\}$.

Our  residue formula in this case reduces to the following iterated
residues:
  $$\iota_{\Phi_h}(a_1,a_2)=
\res_{z_2=0}\left(\res_{z_1=0}\left(\frac{e^{a_1 z_1+a_2
z_2}\;dz_1\,dz_2}
{(1-e^{-z_1})^h(1-e^{-z_2})^h(1-e^{-(z_1+z_2)})^h}\right)\right)$$
for any $(a_1,a_2)\in S_{1,h}=\c_1-\Box(\Phi_h)$, while

$$\iota_{\Phi_h}(a_1,a_2)=
\res_{z_1=0}\left(\res_{z_2=0}\left(\frac{e^{a_1 z_1+a_2
z_2}\;dz_1\,dz_2}{(1-e^{-z_1})^h(1-e^{-z_2})^h(1-e^{-(z_1+z_2)})^h}\right)\right)$$
for any $(a_1,a_2)\in S_{2,h}=\c_2-\Box(\Phi_h)$.

Pictures of the chambers and of the sets $\Box(\Phi_h)$,
$S_{1,h}$, $S_{2,h}$, $S_{1,h}\cap S_{2,h}$ are given on figures 4
through 8 in the Appendix.

Let us give the explicit result for $h=3$. We denote by
$\iota[\c,\Phi_3]$ the polynomial function of $(a_1,a_2)$, which
coincides with the vector partition function $\iota_{\Phi_3}$ on the
chamber $\c$.

The function
$\iota[\c_1,\Phi_3]$ is equal to
\[
\frac1{14}\binom{a_2+5}5(7{a_1}^{2} - 7 {a_1}{a_2} + 2{a_2}^{2} +
21{a_1} - 9{a_2} + 14 ),
\]
so it vanishes along the lines $a_2=-1,-2,-3,-4,-5$. By symmetry,
the function
\[
\iota[\c_2,\Phi_3]=\frac1{14}\binom{a_2+5}5(2{a_{1}}^{2} - 7
{a_{1}}{a_{2}} + 7{a_{2}}^{2} - 9{a_{1}} + 21{a_{2}} + 14 )
\]
vanishes along the lines $a_1=-1,-2,-3,-4,-5$.  These vanishing
properties may be deduced from the Ehrhart reciprocity Theorem.
Our results show that the  functions $\iota[\c_1,\Phi_3]$ and $\iota[\c_2,\Phi_3]$ coincide on
the integral points in $S_{1,h}\cap S_{2,h}$. Indeed, we have
$$\iota[\c_1,\Phi_3]-\iota[\c_2,\Phi_3]=\frac1{14}\binom{a_1-a_2+2}5
 (2{a_{1}}^{2} + 3{a_{1}}{a_{2}} + 2{a_{2}}^{2} + 21
{a_{1}} + 21{a_{2}} + 59),
$$
thus the two polynomial functions $\iota[\c_1,\Phi_3]$ and
$\iota[\c_2,\Phi_3]$ coincide along the lines
$a_1-a_2=-2,-1,0,1,2$.
\end{example}

{\bf Acknowledgments:} We would like to thank Michel Brion,
Velleda Baldoni-Silva for helpful discussions and Charles Cochet
for his drawings and his careful reading of the manuscript.

\section{Partial fraction decompositions}
\labelm{sec:parfrac}

\subsection{Complex hyperplane arrangements.}
\labelm{sec:cxhpa}

Let $E$ be a $n$-dimensional complex vector space. If $\alpha\in
E^*$ is a nonzero linear form on $E$, then we denote by
$H_{\alpha}$ the hyperplane $\{z\in E|\, \langle \alpha,z\rangle
=0\}$.

An {\em arrangement} $\CA$ of hyperplanes in $E$ is a finite collection
of hyperplanes.  Thus one may associate an arrangement
$\CA(\Delta)$ to any finite subset $\Delta\subset E^*$ of nonzero
linear forms; this arrangement consists of the set of hyperplanes
$H_{\alpha}$,  where $\alpha$ varies in $\Delta$. Conversely,
given an arrangement $\CA=\{H_1,\ldots, H_{N}\}$ of hyperplanes,
we choose for each hyperplane $H_i\in \CA$ a linear form
$\alpha_i\in E^*$ such that $H_i=H_{\alpha_i}$.  Note that such a
linear form $\alpha_i$ is defined only up to proportionality.

We will call a set $\{H_i\}_{i=1}^m$ of $m$ hyperplanes in $E$ {\em
  independent} if $\dim \cap H_i = n-m$. This is equivalent to saying
that the corresponding linear forms are linearly independent. We
will say that an hyperplane $L_0$ is {\em dependent} on an
arrangement $\{L_i\}_{i=1}^R$, if the linear form $\alpha_0$
defining $L_0$ can be expressed as a linear combination of the
forms $\alpha_i$ $(1\leq i\leq R)$ defining $L_i$. An arrangement
$\CA=\{H_1,\ldots, H_N\}$ is called {\em essential} if $\cap_{i}
H_i=\{0\}$. Writing $\CA=\CA(\Delta)$, this condition means that
the set of vectors $\Delta$ generates $E^*$.

Let $\CA=\{H_1,\ldots, H_{N}\}$ be an arrangement of hyperplanes
and $\Delta=\{\alpha_1,\ldots,\alpha_N\}$ be a set of linear forms
such that $\CA=\CA(\Delta)$. Let us denote by $R_{\CA}$ the ring
of rational functions on $E$ with poles along $\cup_{i=1}^NH_i$.
Then each element $F\in R_{\CA}$ can be written as
$F=P/\prod_{i=1}^R\beta_i$, where $P$ is a polynomial and
$[\beta_1,\ldots, \beta_R]$ is a sequence of elements of $\Delta$.
The algebra $R_{\CA}$ is $\Z$-graded by the degree.  Denote by
$\CB(\Delta)$ the set of $n$-element subsets of $\Delta$ which are
bases of $E^*$. Given $\sigma\in \CB(\Delta)$, we can form the
following elements of $R_{\CA}$:
\begin{equation}
  \label{eq:fsigma}
f_{\sigma}(z):=\frac{1}{\prod_{\alpha\in
    \sigma}\alpha(z)}.
\end{equation}

Clearly, the vector space spanned by the functions $f_\sigma$ for
$\sigma\in \CB(\Delta)$ depends only on $\CA$.
\begin{definition}
The subspace $S_{\CA}$ of $R_{\CA}$ spanned by the functions
$f_{\sigma},~ \sigma\in\CB(\Delta),$ is called the space of
{\em simple elements} of $R_{\CA}$: $$ S_{\CA}=\sum_{\sigma\in
 \CB(\Delta)}\C f_{\sigma}. $$
\end{definition}

The vector space $S_{\CA}$ is contained in the homogeneous component
of degree $-n$ of $R_{\CA}$.  If $\CA$ is not an essential
arrangement, then the set $\CB(\Delta)$ is empty and $S_\CA=\{0\}$.

We let vectors $v\in E$ act on $R_{\CA}$ by  differentiation:
$$
\partial_vf(z):= \frac{d}{d\epsilon}f(z+\epsilon v)|_{\epsilon=0}. $$

Then the following holds (\cite{bv4}, Proposition 7).
\begin{theorem}\labelm{bv2}
There is a direct sum decomposition
  $$
  R_\CA= \left(\sum_{v\in E}\partial_v R_{\CA}\right)\oplus S_{\CA}. $$
\end{theorem}

 As a corollary of this decomposition, we can define the
projection map
$$ \Tres_{\CA}: R_{\CA}\to S_{\CA}, $$
called the {\em total residue}.  The following assertion is
obvious.

\begin{lemma}\label{obvious}
 Assume that $\CA$ is a subset of $\CB$. Then
$$R_\CA\subset R_\CB, \quad S_\CA\subset S_\CB.$$
Furthermore, if $f\in R_\CA$, then $\Tres_\CB(f)$ belongs to
$S_\CA$ and $$\Tres_\CB(f)=\Tres_\CA(f). $$
\end{lemma}

We denote by $R_\mathrm{hp}$ the space of rational functions on
$E$ with poles along hyperplanes. In other words, $R_\mathrm{hp}$
is the union of the spaces $R_{\CA}$ as $\CA$ varies over all
arrangements of hyperplanes in $E$. The preceding lemma shows that
the assignment $ \Tres f=\Tres_{\CA}f$, for $f\in R_{\CA}$, is
well defined on $R_\mathrm{hp}$.  For $f\in R_\mathrm{hp}$, the
function $\Tres f$ is a linear combination of functions
$f_{\sigma}$, defined in~\eqref{eq:fsigma},  where the set
$\sigma$ is a basis of $E^*$ such that $\CA(\sigma)$ is contained
the set of poles of $f$. The map $\Tres$ vanishes on the space
$R_\mathrm{hp}(m)$ of homogeneous fractions of degree $m$ unless
$m+n=0$. In particular, if $\phi=f_{\sigma} P$, where $P$ is a
polynomial and $\sigma$ is a basis of $E^*$, then the total
residue of $\phi$ is $P(0)f_\sigma$.

The total residue also vanishes on all homogeneous elements of
degree $-n$ of the form ${P}/{\prod_{i=1}^R\beta_i}$, where $P$ is
a homogeneous polynomial of degree $R-n$ and vectors
$\{\beta_i\}_{i=1}^R$ do not generate $E^*$.

Denote by $\hatr$ the space of formal meromorphic functions on $E$
near zero, with poles along hyperplanes. In other words, any
element of $\hatr$ can be written as ${P}/{\prod_{i=1}^R\beta_i}$,
where $P$ is a formal power series and $[\beta_1,\ldots, \beta_R]$
is a sequence of elements of $E^*$. The total residue
extends to the space $\hatr$ by defining
$$\Tres\left(\frac{P}{\prod_{i=1}^R\beta_i}\right)=\Tres
\left(\frac{P_{[R-n]}}{\prod_{i=1}^R\beta_i}\right),$$ where  $P_{[R-n]}$ is
the homogeneous component of $P$ of degree $R-n$.

For example, if $a\in E^*$, then the element $e^a$ denotes the
power series $\sum_{k=0}^{\infty}{a^k}/{k!}$ and the total residue
of ${e^a}/{\prod_{i=1}^R \beta_i}$ is, by definition,  equal to
the total residue of ${a^{R-n}}/{((R-n)!\prod_{i=1}^R \beta_i)}$.
Again, this total residue vanishes if the linear forms
$\{\beta_i\}_{i=1}^R$ do not span $E^*$.

\begin{example}
Consider the function
$$g(z_1,z_2)=\frac{e^{z_1}}
{(1-e^{-z_1})(1-e^{-z_2})(1-e^{-(z_1-z_2)})^2}.$$ Thus we write
$$g=\frac{P}{z_1 z_2 (z_1-z_2)^2}\text{ with } P=e^{z_1}
\frac{z_1}{1-e^{-z_1}}\frac{z_2}{1-e^{-z_2}}
\left(\frac{z_1-z_2}{1-e^{-(z_1-z_2)}}\right)^2.$$

To compute the total residue of $g$, we need the term of degree $2$
in the expansion of $P$ at the origin. This is $P_{[2]}:=3 z_1^2
- \frac {13}{12}z_{1}z_{2} $. Then
$$\frac{P_{[2]}}{z_1
  z_2(z_1-z_2)^2}=\frac{23}{12}\frac{1}{(z_1-z_2)^2}+3\frac{1}{z_2(z_1-z_2)}.$$
The total residue of the first fraction is equal to $0$, and we obtain
the answer
$$ \Tres\, g=\frac{3}{z_2(z_1-z_2)}.$$
\end{example}

The following statements follow from the discussion above. We will use
them later.
\begin{lemma}
\labelm{poles}
Consider the meromorphic function $F$ on $E$ expressed as
$$F(z)=\frac{e^{\langle a,z\rangle}}{\prod_{i=1}^N(1-u_i e^{-\langle
    \beta_i,z\rangle})},$$
where $[\beta_1,\ldots, \beta_N]$ is a sequence of elements of
$E^*$ and the coefficients $u_i$, $i=1,\ldots,N$, are nonzero
complex numbers. Then
\begin{itemize}
  \item $\Tres F=0$ if those $\beta_j$ for which $u_j=1$ do not span $E^*$.
  \item If the set $\sigma=\{\beta_j|\, u_j=1\}$ forms a basis of
    $E^*$, then
  $$(\Tres F)(z)=\frac{1}{\prod_{\beta_j\in \sigma}
\langle\beta_j,z\rangle}\frac{1}{\prod_{\beta_k\notin
\sigma}(1-u_k)}.$$
\end{itemize}
\end{lemma}

\subsection{Rational hyperplanes arrangements}
\labelm{sec:rathpa}

Let $V$ be a  real vector space of dimension $n$. For $\alpha\in
V^*$, we  denote by $H_{\alpha}=\{v\in  V|\, \langle
\alpha,v\rangle =0\}$, this time,  the real hyperplane determined by $\alpha$.

Again, let $\Gamma$ be a rank-$n$ lattice in $V$ and denote by
$\Gamma^*\subset V^*$ the dual lattice. This means that if
$\alpha\in \Gamma^*$ and $\gamma\in \Gamma$, then $\langle
\alpha,\gamma\rangle \in \Z$. We denote by $\C[\Gamma^*]$ the ring
of functions on $V_\C$ generated by the exponential functions
$z\mapsto e^{\langle
  \xi,z\rangle }$, $\xi\in \Gamma^*$.

\begin{definition}
  An arrangement $\CA$ of real hyperplanes in $V$ is $\Gamma$-{\em
    rational} if $\CA=\CA(\Delta)$ for some finite subset $\Delta$ of
  $\Gamma^*$. We simply say that $\CA$ is {\em rational} if $\Gamma$
  has been fixed.
\end{definition} Given a rational arrangement $\CA=\{H_1,\ldots,
H_{N}\}$ of hyperplanes, for each hyperplane $H_i\in \CA$ we choose a
linear form $\alpha_i\in \Gamma^*$ such that $H_i=H_{\alpha_i}$. If
$H_{\alpha}=H_{\beta}$ with both $\alpha$ and $\beta$ in $\Gamma^*$,
then $\alpha$ and $\beta$ are proportional with a {\em rational}
coefficient of proportionality.

For any $u\in \C^*$, $\alpha\in \Gamma^*$, consider the meromorphic
function on $V_\C$ defined by
$$g[\alpha,u](z)=\frac{1}{1-ue^{\langle \alpha,z\rangle }}.$$
If $u=e^a$ with $a\in \C$, then the set of poles of the function
$g[\alpha,u]$ is the set $\{z\in V_\C|\, \langle \alpha, z\rangle
+a\in 2i\pi \Z\}$.

\begin{definition}
  We denote by $M^{\Gamma}$ the ring of meromorphic functions on
$V_{\C}$ generated by $\C[\Gamma^*]$ and by the functions
$g[\alpha,u]$, where $u$ varies in $\C^*$ and $\alpha$ in
$\Gamma^*$.
Given  a finite subset $\Delta$ of nonzero elements of
$\Gamma^*$, denote by $M^{\Gamma\Delta}$ the ring of meromorphic
functions on $V_\C$ generated by the ring $\C[\Gamma^*]$  and by
the meromorphic functions $g[\alpha,u]$, where $u$ varies in $\C^*$
and now $\alpha$ is restricted  to be a member of  the finite set
$\Delta$.
\end{definition}
 Thus, to be explicit, a function $F\in M^{\Gamma}$
can  be written, by reducing to a common denominator, as
$$F(z)=\frac{\sum_{\xi\in I} c_\xi e^{\langle \xi,z\rangle }}
{\prod_{k=1}^R(1-u_k e^{\langle \alpha_k,z\rangle })}$$ where $I$
is a finite subset of  $\Gamma^*$;  $u_k,c_\xi\in \C^*$, and the
elements $\alpha_k$ are in $\Gamma^*$. If in addition
$\alpha_k\in\Delta$, then this function is in $M^{\Gamma\Delta}$.

If we write $z=x+iy$ with $x,y\in V$, then the function
$y\mapsto F(x+iy)$ is periodic: $F(x+i(y+2\pi
\gamma))=F(x+iy)$ for any $\gamma \in \Gamma$. Thus functions $F\in
M^{\Gamma}$ induce functions on the complexified torus $V_\C/2i\pi
\Gamma$.

\begin{lemma}\labelm{bete}
  Let $\Delta$ and $\Delta'$ be two finite subsets of $\Gamma^*$ such
  that $\CA(\Delta)=\CA(\Delta')$. Then we have
  $M^{\Gamma\Delta}=M^{\Gamma\Delta'}$. Thus the ring
  $M^{\Gamma\Delta}$ depends only on the rational hyperplane
  arrangement $\CA(\Delta)$.
\end{lemma}
\begin{proof}
Let us note the following identities:
 $$\frac{1}{(1-e^a e^{k z})}=\frac{1}{\prod_{\zeta,
\zeta^k=1}(1-\zeta e^{a/k}e^z)},$$
$$\frac{1}{1-ue^z}=\frac{1+ue^z+u^2e^{2z}+\cdots+u^{(n-1)}e^{(n-1)z}}{1-u^ne^{nz}},$$
$$\frac{1}{1-ue^z}=\frac{u^{-1}e^{-z}}{u^{-1}e^{-z}-1},$$
where
$n,k\in\Z$, and $a,u,z\in \C$.

These identities show that $M^{\Gamma\Delta}$ does not change if
we multiply one of the elements of $\Delta$ by a non-zero integer.
This implies the Lemma since any two sets $\Delta,\Delta'\subset\Gamma^*$ such that
$\CA(\Delta)= \CA(\Delta')$ may be transformed into each other by
such an operation.
\end{proof}

Now we can give the following definition:
\begin{definition}
  Let $\CA$ be a $\Gamma$-rational hyperplane arrangement in a vector
  space $V$. Define $M^{\Gamma\CA}$ to be the ring $M^{\Gamma\Delta}$,
  where $\Delta\subset \Gamma^*$ is an arbitrary subset such that
  $\CA=\CA(\Delta)$.
\end{definition}

It is clear that, if $\CB$ is a subset of $\CA$, then
$M^{\Gamma\CB}$ is a subring of  $M^{\Gamma\CA}$.

\subsection{Behavior at $\infty$}
\labelm{sec:behave}

Consider a function $F\in M^{\Gamma}$. The function of the real
variable $y\mapsto F(x+iy)$ is $2\pi\Gamma$-periodic. In this
paragraph, we study the behavior of the function of the real variable
$x\mapsto F(x+iy)$ at $\infty$.

  Let $z_0\in V_\C$ be not a pole of $F$.
Then, for all $v\in V$, the function
 $s\mapsto F(z_0+sv)$ is well-defined when $s$ is a sufficiently large
 real number.

\begin{definition}
Let $F\in M^{\Gamma}$. Assume that $z_0\in V_\C$ is not a pole of
$F$. Define $\Box(z_0,F)$ to be the set of $\mu\in V^*$ such that
for every $v\in V$, the function $s\mapsto e^{s\langle
\mu,v\rangle }F(z_0+sv)$ remains bounded when $s$ is real and
tends to $+\infty$.
\end{definition}

{\bf Example}. Let $F(z)={1}/(1-e^z)$; pick $z_0\notin 2i\pi \Z$.
Then  $\Box(z_0,F)=[0,1]$. Indeed, the function
$\theta(s,v)={e^{s\mu v}}/(1-e^{z_0+s v})$ is bounded as $s$ tends
to $\infty$ if and only $0\leq \mu\leq 1$: when $v=0$, the
function $\theta(s,v)$ is the constant ${1}/(1-e^{z_0})$; if
$v>0$, we obtain the condition $\mu\leq 1$; if $v<0$, we obtain
the condition $\mu\geq 0$. Note that if $v\neq 0$, then, for
$\mu\in ]0,1[$, the function $s\mapsto \theta(s,v)$ tends to $0$
when $s$ tends to $\infty$.
\begin{definition}
  For two subsets $A$ and $B$ of a real vector space, we denote by
  $A+B$ their Minkowski sum:
\[ A+B = \{a+b|\, a\in A, \, b\in B\}.\]
\end{definition}
Note that the sum of convex sets is convex.

\begin{proposition}\labelm{desc}
  Let $F\in M^\Gamma$ be written in the form $$F(z)=\frac{\sum_{\xi\in
      I} c_\xi e^{\langle \xi,z\rangle }} {\prod_{i=1}^R(1-u_i
    e^{\langle \alpha_i,z\rangle })},$$
  where $I$  is a finite subset of $\Gamma^*$,
  $\alpha_i$ are in $\Gamma^*$, and all the constants $c_\xi$ and $u_i$ are
  nonzero complex numbers.  Assume that $z_0\in V_\C$ is such that
  $\prod_{i=1}^R(1-u_i e^{\langle \alpha_i,z_0\rangle })\neq 0.$ Then
  \begin{equation}
    \label{eq:box}
\Box(z_0,F)=\{\mu\in V^*\,| \,\mu+\xi\in \sum_{i=1}^R
[0,1]\alpha_i\, {\rm for\, all\,} \xi\in I \}.
  \end{equation}
\end{proposition}
\begin{proof}
  The set described on the right hand side of \eqref{eq:box} is easily
  seen to be contained in $\Box(z_0,F)$. Indeed, let $\mu\in V^*$ be such
  that $\mu+\xi$ belongs to the set $\sum_{i=1}^R [0,1]\alpha_i$ for
  all $\xi\in I.$ We write $F(z)=\sum_{\xi\in I}c_\xi F_{\xi}(z)$ with
$$F_\xi(z)=\frac{ e^{\langle \xi,z\rangle }} {\prod_{i=1}^R(1-u_i
e^{\langle \alpha_i,z\rangle })}.$$ Let us show that for each
$\xi\in I$, the function $s\mapsto e^{s\langle \mu,v\rangle }F_\xi(z_0+sv)$
remains bounded when $s$ tends to $\infty$. We have
$\mu+\xi=\sum_{i=1}^R t_i \alpha_i$ with $0\leq t_i\leq 1$ and we may
write $e^{s\langle \mu,v\rangle }F_\xi(z_0+sv)$ as
$$e^{\langle \xi,z_0\rangle
}\frac{e^{s\langle \mu+\xi,v\rangle }}{\prod_{i=1}^R(1-u_i
e^{\langle \alpha_i,z_0\rangle }e^{s\langle
\alpha_i,v\rangle})}=e^{\langle \xi,z_0\rangle}
 \prod_{i=1}^R\frac{e^{st_i\langle \alpha_i,v\rangle }} {(1-u_i
e^{\langle \alpha_i,z_0\rangle }e^{s\langle \alpha_i,v\rangle
})}.$$ As each of the factors on the right hand side remains
bounded when $s$ tends to $\infty$, we have shown that $\mu\in
\Box(z_0,F)$.

We now prove the converse. Let $\mu$ be such that the function
$$s\mapsto e^{s\langle \mu,v\rangle}F(z_0+sv)$$
is bounded as $s\rightarrow\infty$ for any $v\in V$. Assume that,
nevertheless, there exists $\nu$ in the set $I$ such that
$\mu+\nu$ is not in the convex polytope $\Pi:=\sum_{i=1}^R
[0,1]\alpha_i$. The vectors $\alpha_{\bfk}=\sum_{i\in
\bfk}\alpha_i$, where $\bfk$ is a subset of $\{1,\ldots,R\}$, are
all in the polytope $\Pi$.  Thus there exists $w\in V$ and $a\in
\R$ such that $\langle \sum_{i\in \bfk} \alpha_i,w\rangle <a$ for
all subsets $\bfk$ of $\{1,2,\ldots,R\}$, while $\langle
\mu+\nu,w\rangle > a$. The set of such vectors $w$ is an open set
in $V$.

We write $$e^{s\langle \mu,v\rangle}F(z_0+sv)=\frac{P(s,v)}{D(s,v)},$$
with $P(s,v)=\sum_{\xi} c_\xi e^{\langle \xi,z_0\rangle }e^{s\langle
  \mu+\xi,v\rangle }$ and $ D(s,v)=\prod_{i=1}^R(1-u_i e^{\langle
  \alpha_i,z_0+sv\rangle })$.  Then $D(s,v)=\prod_{i=1}^R(1-u_i
e^{\langle \alpha_i,z_0+sv\rangle })=\sum_{\bfk} c_{\bfk} e^{\langle
  \alpha_{\bfk},z_0\rangle }e^{s\langle \alpha_{\bfk},v\rangle }$,
  for some constants $c_{\bfk}$.
   Note that the function
$s\mapsto D(s,v)$ does not vanish identically, as
$D(0,v)=\prod_{i=1}^R(1-u_i e^{\langle \alpha_i,z_0\rangle }).$
Thus
 for any $w$ such that $\langle \sum_{i\in \bfk}
\alpha_i,w\rangle <a$, the denominator $D(s,w)$ can be rewritten
as a finite sum of exponentials $\sum_k h_k e^{b_k s}$ with
distinct exponents $b_k$ and nonzero coefficients $h_k$. We
clearly have $\max_k(b_k)<a$.

Consider now
$$P(s,w)=\sum_{\xi} c_\xi e^{\langle \xi,z_0\rangle }e^{s\langle
  \mu+\xi,w\rangle }.$$
Since the set $\{w\in V\,|\,\langle \mu+\nu,w\rangle > a, \langle
\sum_{i\in \bfk} \alpha_i,w\rangle <a\}$ is open, we can choose an
element $w_0$ in it such that the numbers $\langle
\mu+\xi,w_0\rangle $ are {\em distinct} for all $\xi\in I$. Then
the numerator $P(s,w_0)$ may be rewritten as a sum of exponentials
$\sum_j c_j e^{a_js}$ with nonzero constants $c_j$, and distinct
exponents $a_j$ such that $\max_j(a_j)>a$.  Thus the function
$F(s,w_0)$ is equal to the quotient $\sum_jc_je^{a_js}/\sum_k h_k
e^{b_k s}$, which is equivalent to $c e^{s(\max_ja_j -\max_k
b_k)}$ as $s\rightarrow+\infty$ ($c\neq0$).  The exponent is
positive, hence the function $s\mapsto F(s,w_0)$ tends to $\infty$
when $s$ tends to $+\infty$.  This contradicts our assumption on
$\mu$, and the proof of Proposition \ref{desc} is complete.
\end{proof}

Let $F\in M^{\Gamma}$. As a consequence of Proposition \ref{desc},
 the set $\Box(z_0,F)$ is independent of the choice of $z_0$.
 \begin{definition}
 Let $F\in M^{\Gamma}$ and $z_0$ be an arbitrary element of
 $V_\C$ which is not a pole of $F$. We denote by $\Box(F)$ the set
 $\Box(z_0,F)$.
 \end{definition}

 The set $\Box(F)$ is easy to determine, using any presentation of $F$
 as a fraction.

 \begin{example}
Let
$$F(z)=\frac{1}{1-e^z} =\frac{1+e^z}{1-e^{2z}}.$$
Using the first expression, we obtain $\Box(F)=[0,1]$. Using the
second expression, we obtain $\Box(F)=[0,2]\cap [-1,1]$.
 \end{example}

\begin{lemma}
Let $F\in M^{\Gamma}$ and $\mu\in V^*$. Assume that $\mu$ is in the
interior of $\Box(F)$ and that  $z_0$ is not a pole of $F$.  Then  for
all {\em nonzero } $v\in V$,  the function $s\mapsto e^{s\langle
\mu,v\rangle }F(z_0+sv)$ tends to zero  when $s$ is real and tends
to $+\infty$.
\end{lemma}

\begin{proof}
  Consider $F\in M^{\Gamma}$ written as in Proposition \ref{desc} and
  let us return to the first part of the proof of this Proposition. If
  the interior of $\Box(F)$ is nonempty, then the linear forms
  $\alpha_i$ necessarily generate $V^*$.  Furthermore, if $\mu$ is in
  the interior of $\Box(F)$, then, for each $\xi\in I$, we can write
  $\mu+\xi=\sum_{i=1}^Rt_i\alpha_i$ with $0<t_i<1$. Each factor
  ${e^{st_i\langle \alpha_i,v\rangle }}/{(1-u_i e^{\langle
      \alpha_i,z_0\rangle }e^{s\langle \alpha_i,v\rangle })}$ remains
  bounded when $s$ tends to $\infty$. Since $v$ is not equal to $0$,
  there exists at least one linear form $\alpha_j$ with $\langle
  \alpha_j,v\rangle \neq 0$. The corresponding factor tends to $0$
  when $s$ tends to $\infty$, and we obtain the Lemma.
\end{proof}

\begin{definition}
  \begin{itemize}
  \item For $\mu\in V^*$, we denote by $M^{\Gamma}(\mu)$ the set of $F\in
  M^{\Gamma}$ such that $\mu\in \Box(F)$.
Similarly,  for $\mu\in V^*$ and a $\Gamma$-rational arrangement
$\CA$, let
$$M^{\Gamma\CA}(\mu)=\{F\in M^{\Gamma\CA}|\,\mu\in\Box(F)\}.$$
\item Let $F\in M^{\Gamma}$ and $\mu\in \Box(F)$. A decomposition
  $F=\sum F_i$ of $F$ into a sum of terms from $M^\Gamma$ will be called
  {\em $\mu$-admissible} if $\mu\in \Box(F_i)$ for every $i$.
  \end{itemize}
\end{definition}

We have the following obvious inclusions:
\begin{lemma}
  \labelm{twoparts}
Let $F,G\in M^{\Gamma}$. Then
$$\Box(F)\cap \Box(G)\subset \Box(F+G)\quad\text{and}\quad
\Box(F)+\Box(G)\subset\Box(FG).$$
\end{lemma}

\begin{remark} \labelm{samedenom}
  A consequence of Proposition \ref{desc} is that if $F=\sum_{i\in
    I}{P_i}/{D}$ is a sum of fractions from $M^{\Gamma}$ with the
  same denominator, then $F\in M^{\Gamma}(\mu)$ if and only if
  ${P_i}/{D}\in M^{\Gamma}(\mu)$ for each $i\in I$. However, for a
  decomposition $F=\sum_i {P_i/D_i}$ with different denominators, the
  inclusion $\cap_i \Box(P_i/D_i)\subset \Box(F)$ is strict in
  general.
\end{remark}

\begin{example}
Set
$$F_1=\frac{1}{(1-e^{z_1})(1-e^{z_2})},\,
F_2=\frac{1}{(1-e^{z_1+z_2})(1-e^{z_2})},\,
F_3=\frac{1}{(1-e^{z_1+z_2})(1-e^{-z_1})}.$$ Then we have
$F_1=F_2-F_3$. Figure 1  shows the three parallelograms
$\Box(F_1)$, $\Box(F_2)$ and $ \Box(F_3)$. Clearly, $\Box(F_2)\cap
\Box(F_3)$ is strictly smaller than $\Box(F_1)$ (cf. Figure 2).
\end{example}

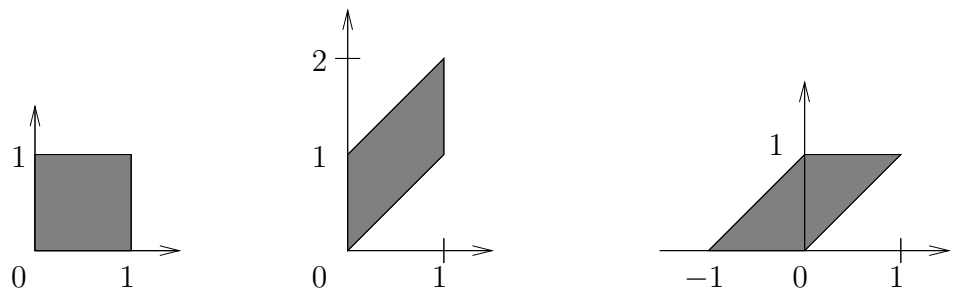
\begin{figure}[ht] \label{boxes}
\begin{center}
\setlength{\unitlength}{0.00083333in}
\begingroup\makeatletter\ifx\SetFigFont\undefined%
\gdef\SetFigFont#1#2#3#4#5{%
  \reset@font\fontsize{#1}{#2pt}%
  \fontfamily{#3}\fontseries{#4}\fontshape{#5}%
  \selectfont}%
\fi\endgroup%
{\renewcommand{\dashlinestretch}{30}
\begin{picture}(5862,1752)(0,-10)
\texture{44555555 55aaaaaa aa555555 55aaaaaa aa555555 55aaaaaa aa555555 55aaaaaa 
	aa555555 55aaaaaa aa555555 55aaaaaa aa555555 55aaaaaa aa555555 55aaaaaa 
	aa555555 55aaaaaa aa555555 55aaaaaa aa555555 55aaaaaa aa555555 55aaaaaa 
	aa555555 55aaaaaa aa555555 55aaaaaa aa555555 55aaaaaa aa555555 55aaaaaa }
\path(150,225)(1050,225)
\path(150,225)(1050,225)
\path(930.000,195.000)(1050.000,225.000)(930.000,255.000)
\path(150,225)(150,1125)
\path(150,225)(150,1125)
\path(180.000,1005.000)(150.000,1125.000)(120.000,1005.000)
\path(2100,225)(3000,225)
\path(2100,225)(3000,225)
\path(2880.000,195.000)(3000.000,225.000)(2880.000,255.000)
\path(2100,225)(2100,1725)
\path(2100,225)(2100,1725)
\path(2130.000,1605.000)(2100.000,1725.000)(2070.000,1605.000)
\path(4050,225)(5850,225)
\path(4050,225)(5850,225)
\path(5730.000,195.000)(5850.000,225.000)(5730.000,255.000)
\path(4950,225)(4950,1275)
\path(4950,225)(4950,1275)
\path(4980.000,1155.000)(4950.000,1275.000)(4920.000,1155.000)
\shade\path(150,825)(750,825)(750,225)
	(150,225)(150,825)
\path(150,825)(750,825)(750,225)
	(150,225)(150,825)
\shade\path(2100,825)(2700,1425)(2700,825)
	(2100,225)(2100,825)
\path(2100,825)(2700,1425)(2700,825)
	(2100,225)(2100,825)
\shade\path(4350,225)(4950,825)(5550,825)
	(4950,225)(4350,225)
\path(4350,225)(4950,825)(5550,825)
	(4950,225)(4350,225)
\path(2025,1425)(2175,1425)
\path(2025,1425)(2175,1425)
\path(2700,300)(2700,150)
\path(2700,300)(2700,150)
\path(5550,300)(5550,150)
\path(5550,300)(5550,150)
\path(4950,225)(4950,825)
\path(4950,225)(4950,825)
\put(2625,0){\makebox(0,0)[lb]{\smash{{{\SetFigFont{12}{14.4}{\rmdefault}{\mddefault}{\updefault}$1$}}}}}
\put(675,0){\makebox(0,0)[lb]{\smash{{{\SetFigFont{12}{14.4}{\rmdefault}{\mddefault}{\updefault}$1$}}}}}
\put(4200,0){\makebox(0,0)[lb]{\smash{{{\SetFigFont{12}{14.4}{\rmdefault}{\mddefault}{\updefault}$-1$}}}}}
\put(4875,0){\makebox(0,0)[lb]{\smash{{{\SetFigFont{12}{14.4}{\rmdefault}{\mddefault}{\updefault}$0$}}}}}
\put(5475,0){\makebox(0,0)[lb]{\smash{{{\SetFigFont{12}{14.4}{\rmdefault}{\mddefault}{\updefault}$1$}}}}}
\put(1875,750){\makebox(0,0)[lb]{\smash{{{\SetFigFont{12}{14.4}{\rmdefault}{\mddefault}{\updefault}$1$}}}}}
\put(1875,1350){\makebox(0,0)[lb]{\smash{{{\SetFigFont{12}{14.4}{\rmdefault}{\mddefault}{\updefault}$2$}}}}}
\put(1875,0){\makebox(0,0)[lb]{\smash{{{\SetFigFont{12}{14.4}{\rmdefault}{\mddefault}{\updefault}$0$}}}}}
\put(0,750){\makebox(0,0)[lb]{\smash{{{\SetFigFont{12}{14.4}{\rmdefault}{\mddefault}{\updefault}$1$}}}}}
\put(0,0){\makebox(0,0)[lb]{\smash{{{\SetFigFont{12}{14.4}{\rmdefault}{\mddefault}{\updefault}$0$}}}}}
\put(4725,825){\makebox(0,0)[lb]{\smash{{{\SetFigFont{12}{14.4}{\rmdefault}{\mddefault}{\updefault}$1$}}}}}
\end{picture}
}

\caption{The sets $\Box(F_1)$, $\Box(F_2)$  and $\Box(F_3)$.}
\end{center}
\end{figure}

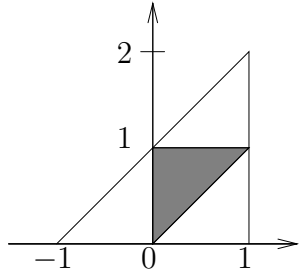
\begin{figure}[ht] \label{compare}
\begin{center}

\setlength{\unitlength}{0.00083333in}
\begingroup\makeatletter\ifx\SetFigFont\undefined%
\gdef\SetFigFont#1#2#3#4#5{%
  \reset@font\fontsize{#1}{#2pt}%
  \fontfamily{#3}\fontseries{#4}\fontshape{#5}%
  \selectfont}%
\fi\endgroup%
{\renewcommand{\dashlinestretch}{30}
\begin{picture}(1824,1677)(0,-10)
\texture{44555555 55aaaaaa aa555555 55aaaaaa aa555555 55aaaaaa aa555555 55aaaaaa 
	aa555555 55aaaaaa aa555555 55aaaaaa aa555555 55aaaaaa aa555555 55aaaaaa 
	aa555555 55aaaaaa aa555555 55aaaaaa aa555555 55aaaaaa aa555555 55aaaaaa 
	aa555555 55aaaaaa aa555555 55aaaaaa aa555555 55aaaaaa aa555555 55aaaaaa }
\path(12,150)(1812,150)
\path(12,150)(1812,150)
\path(1692.000,120.000)(1812.000,150.000)(1692.000,180.000)
\path(912,150)(912,1650)
\path(912,150)(912,1650)
\path(942.000,1530.000)(912.000,1650.000)(882.000,1530.000)
\path(1512,150)(1512,1350)(312,150)
\shade\path(912,750)(1512,750)(912,150)(912,750)
\path(912,750)(1512,750)(912,150)(912,750)
\path(837,1350)(987,1350)
\put(687,750){\makebox(0,0)[lb]{\smash{{{\SetFigFont{12}{14.4}{\rmdefault}{\mddefault}{\updefault}$1$}}}}}
\put(687,1275){\makebox(0,0)[lb]{\smash{{{\SetFigFont{12}{14.4}{\rmdefault}{\mddefault}{\updefault}$2$}}}}}
\put(162,0){\makebox(0,0)[lb]{\smash{{{\SetFigFont{12}{14.4}{\rmdefault}{\mddefault}{\updefault}$-1$}}}}}
\put(837,0){\makebox(0,0)[lb]{\smash{{{\SetFigFont{12}{14.4}{\rmdefault}{\mddefault}{\updefault}$0$}}}}}
\put(1437,0){\makebox(0,0)[lb]{\smash{{{\SetFigFont{12}{14.4}{\rmdefault}{\mddefault}{\updefault}$1$}}}}}
\end{picture}
}
\caption{The sets $\Box(F_1)$, $\Box(F_2)$, $\Box(F_3)$ and
$\Box(F_2)\cap \Box(F_3)$.}
\end{center}
\end{figure}

 The following lemma will allow us to
obtain admissible decompositions of certain specific elements of
$M^{\Gamma}(\mu)$.
\begin{lemma}\labelm{crucial}
 Let  $\alpha_1,\ldots,\alpha_r$ be  nonzero linear forms,
 and let $\alpha_0=-(\alpha_1+\alpha_2+\cdots +\alpha_r)$.
 Let $u_1,\ldots,u_r$ be nonzero complex numbers, and  let
$$F=\frac{1}{\prod_{i=1}^r (1-u_i e^{\alpha_i})}.$$
Set $\mu=\sum_{i=1}^rt_i\alpha_i\in\Box(F)$ with $0\leq t_1\leq t_2\leq
\cdots \leq t_r\leq 1$.

Assume that either the linear form $\alpha_0$ is not identically
zero, or if $\alpha_0=0$ then the product $u_1\cdots u_r\neq 1$.
This ensures that the function $(1-u_1\cdots u_r e^{-\langle
\alpha_0,z\rangle})^{-1}$ is well defined as a meromorphic
function on $V_\C$. Then we have
$$F=\sum_{i=1}^r F_i$$
where $$F_i=(-1)^{i+1}\frac{1}{(1-u_1u_2\cdots u_re^{-\alpha_0})}
\prod_{j=1}^{i-1}\frac{1}{(1-u_j^{-1}e^{-\alpha_j})}\prod_{j=i+1}^r
\frac{1}{(1-u_je^{\alpha_j})},$$ and $\mu\in \Box(F_i)$ for each
$1\leq i\leq r$.
\end{lemma}

\begin{proof}
The equality $F=\sum_{i=1}^r F_i$ is easily verified by
multiplying by  $(1-u_1u_2\cdots u_re^{-\alpha_0})$. The resulting formula in another form is
\begin{equation}\label{obvious2}
 \sum_{i=1}^r\frac{\prod_{j=1}^{i-1}u_j e^{\alpha_j}}{\prod_{j\neq
 i}(1-u_je^{\alpha_j})}=
\frac{1-\prod_{i=1}^r u_ie^{\alpha_i}}
{\prod_{i=1}^r(1-u_ie^{\alpha_i})}.
 \end{equation}

It remains to check
that $\mu\in \Box(F_i)$ for each $i$. We have
$$\mu=\sum_{j=1}^{i-1}t_j\alpha_j+ t_i\alpha_i+\sum_{j=i+1}^r t_j
\alpha_j
=-t_i\alpha_0+\sum_{j=1}^{i-1}(t_j-t_i)\alpha_j+\sum_{j=i+1}^r
(t_j-t_i) \alpha_j$$ where the coefficients of
$-\alpha_0,-\alpha_1,\ldots,-\alpha_{i-1},\alpha_{i+1},\ldots,\alpha_r$
are between $0$ and $1$. Considering the form of the functions $F_i$,
this is exactly the criterion of being in $\Box(F_i)$, and the proof
is complete.
\end{proof}

\begin{remark}
We would like to stress here that the $\mu$-admissible
decomposition of $F$ given in Lemma \ref{crucial} depends on the
position of the element $\mu$ in $\Box(F)$ in an essential manner.
  \end{remark}

  \begin{example}
Let $$F({z_1},z_2)=\frac{1}{(1-e^{z_1})(1-e^{z_2})}$$
and $\mu=(\mu_1,\mu_2)$ in $\Box(F)$, i.e. $0\leq \mu_1\leq 1$ and
$0\leq \mu_2\leq 1$.
Then if $\mu_1\leq \mu_2$, we write
$$F=F_1-F_2$$
with
$$F_1({z_1},z_2)=\frac{1}{(1-e^{{z_1}+z_2})(1-e^{z_2})}\quad\text{and}\quad
F_2({z_1},z_2)=\frac{1}{(1-e^{{z_1}+z_2})(1-e^{-{z_1}})},$$
so that $\mu\in \Box(F_1)\cap \Box(F_2)$.

If $\mu_1\geq \mu_2$,  then the  roles of ${z_1}$ and $z_2$ are reversed, and
$$F=F'_1-F'_2$$
with
$$F'_1({z_1},z_2)=\frac{1}{(1-e^{{z_1}+z_2})(1-e^{{z_1}})}\quad\text{and}\quad
F'_2({z_1},z_2)=\frac{1}{(1-e^{{z_1}+z_2})(1-e^{-z_2})}.$$
 Again, we have $\mu\in \Box(F'_1)\cap \Box(F'_2)$.
  \end{example}

  \begin{example}
Let $$F(z)=\frac{1}{(1-u^{-1} e^{-z})(1-ve^{z})}$$ with $u\neq
v$.
Let $\mu\in [-1,1]$. Then if $0\leq \mu\leq 1$, we write
$$F(z)=F_1(z)-F_2(z)$$
with
$$F_1(z)=\frac{1}{(1-u^{-1}v)}\frac{1}{(1-ve^z)},\quad
F_2(z)=\frac{1}{(1-u^{-1}v)}\frac{1}{(1-ue^{z})}$$
and $\mu\in \Box(F_1)\cap \Box(F_2)$.
If $-1\leq \mu\leq 0$, then we exchange the roles of $z$ and $-z$ and
write
$$F(z)=F'_1(z)-F'_2(z)$$
with $$F'_1(z)=\frac{1}{(1-u^{-1}v)}\frac{1}{(1-u^{-1}e^{-z})},\quad
F'_2(z)=\frac{1}{(1-u^{-1}v)}\frac{1}{(1-v^{-1}e^{-z})},$$
where again $\mu\in \Box(F'_1)\cap \Box(F'_2)$.
  \end{example}

\subsection{Separating variables}
\labelm{sec:sepvar}

Let $\Gamma\subset V$ be a lattice of full rank and $\CA$ be a
$\Gamma$-rational arrangement of hyperplanes in $V$. Recall that
given $\mu\in V^*$,  we
defined $M^{\Gamma\CA}(\mu)$  to be the subspace of functions
 $F$ in $M^{\Gamma\CA}$ such that $\mu\in \Box(F)$.

\begin{lemma}{\rm (}The exchange Lemma{\rm )}
Let $\CA=\{H_1,\ldots,H_m\}$ be a rational arrangement of
hyperplanes and let $H_0$ be a rational hyperplane, which is
dependent on $\CA$.  Denote by  $\CA_i$ the arrangement $\{H_0,
H_1,\ldots,\hat H_i,\ldots, H_m\}$, where we have replaced the
hyperplane $H_i$  by the hyperplane $H_0$. Then, for any
$\mu\in V^*$, we have
$$M^{\Gamma \CA}(\mu)
 \subset \sum_{i=1}^mM^{\Gamma \CA_i}(\mu).$$
\end{lemma}
\begin{proof}
The dependence of $H_0$ on $\CA$ means that there are linear forms
$\alpha_0, \ldots,\alpha_m\in \Gamma^*$ with
$H_0=H_{\alpha_0},H_1=H_{\alpha_1},\ldots,H_m=H_{\alpha_m}$,  such
that $\alpha_0$ may be expressed as a linear combination of the rest
of the $\alpha$s. By using multiples of these linear forms to describe our
hyperplanes and reordering the hyperplanes in $\CA$ if necessary, we
may assume that the relation takes the form
$\alpha_0+\alpha_1+\cdots+\alpha_r=0$, where $r\leq m$.

Let $F\in M^{\Gamma\CA}(\mu)$. By Lemma \ref{bete}, we may write
$$F=\sum_{\xi\in I} c_\xi \frac{e^{\xi}}{D},\quad
D=\prod_{j=1}^R(1-u_j e^{\beta_j}),$$ where $[\beta_1,\ldots,
\beta_R]$ is a sequence of not necessarily distinct elements of
the set $\{\alpha_1,\ldots,\alpha_m\}$ and $c_\xi$ is a nonzero
complex number for $\xi\in I$.  According to Remark
\ref{samedenom}, each of the terms ${e^{\xi}}/{D}$ is in
$M^{\Gamma\CA}(\mu)$, so we may assume that $F$ is of the form
${e^{\xi}}/{D}$ to begin with.

We argue by induction on the length $R$ of the sequence
$[\beta_1,\ldots,\beta_R]$.  If the set $\{\beta_j|\,1\leq j\leq
R\}$ of elements occurring in the sequence is strictly smaller
than the set $\{\alpha_1,\ldots,\alpha_m\}$, then one of the
linear forms $\alpha_i$ does not appear in the sequence
$[\beta_1,\ldots,\beta_R]$, and thus $F$ is already in
$\sum_{i=1}^mM^{\Gamma\CA_i}(\mu)$. Otherwise, reordering the
sequence, we may assume that
$\beta_1=\alpha_1,\beta_2=\alpha_2,\ldots,\beta_r=\alpha_r$.  We
write $D'=\prod_{i=1}^r(1-u_i e^{\alpha_i})$,
$D''=\prod_{j=r+1}^R(1-u_j e^{\beta_j})$, so that $D=D'D''$. As
$\mu\in \Box(F)$, we write $\mu=\mu' +\mu''$ with $\mu'\in
\sum_{i=1}^r[0,1]\alpha_i$ and $\mu''\in -\xi+
\sum_{i=r+1}^R[0,1]\beta_i$. Now
$$\frac{e^{\xi}}{D}=\frac{1}{D'}\frac{e^{\xi}}{D''}\,\,\, \text{ with }
\frac{1}{D'}\in M^{\Gamma}(\mu')\text{  and } \frac{e^{\xi}}{D''}\in
 M^{\Gamma}(\mu'').$$

We may suppose that, after reordering  the first $r$ elements of
the sequence if necessary, we have $\mu'=\sum_{k=1}^r t_k
\alpha_k$  with $0\leq t_1\leq t_2\leq \cdots\leq t_r$.

Using Lemma \ref{crucial}, we write ${1}/{D'}=\sum_{k=1}^r F'_k$,
with $\mu'\in \Box(F'_k)$. Thus we obtain a $\mu$-admissible
decomposition ${e^{\xi}}/{D}=\sum_{k=1}^rF'_k{e^{\xi}}/{D''}$.
More explicitly, writing $u_0=(u_1\cdots u_r)^{-1}$, we obtain the
$\mu$-admissible decomposition
 $$\frac{e^{\xi}}{D}=\sum_{k=1}^r\frac{1}{(1-u_0^{-1}e^{-\alpha_0})}G_k$$
 with
 $$G_k=
 \frac{e^{\xi}}{\prod_{j=1}^{k-1}(1-u_j^{-1}e^{-\alpha_j})\prod_{j=k+1}^r
   (1-u_je^{\alpha_j})\prod_{j=r+1}^R(1-u_j e^{\beta_j})}.$$
 For each
 $1\leq k\leq r$, we have $\mu=t_k\alpha_0+\mu'_k$ with $0\leq t_k\leq
 1$ and $\mu'_k\in \Box(G_k)$. We can apply our induction hypothesis
 to $G_k\in M^{\Gamma\CA}(\mu'_k)$ since the length of the denominator
 of $G_k$ is $R-1$. We then obtain an admissible decomposition of
 $G_k$ as $\sum_{i=1}^m G_k^i$ with $\mu'_k\in \Box(G_k^i)$ and
 $G_k^i\in M^{\Gamma\CA_i}$. According to Lemma \ref{twoparts}, the
 function $G_k^i/(1-u_0^{-1}e^{-\alpha_0})$ is in
 $M^{\Gamma\CA_i}(t_k\alpha_0+\mu'_k)=M^{\Gamma\CA_i}(\mu)$. Hence the proof
 is now complete.
\end{proof}

Clearly, if $\mathcal{B}\subset\CA$ then $M^{\Gamma\CB}(\mu)\subset
M^{\Gamma\CA}(\mu)$. The following crucial partial fraction
decomposition type result holds in the reverse direction:
\begin{theorem}[\cite{sz2}]
\labelm{thm:parfrac}
For each $\mu\in V^*$, we have the equality
$$M^{\Gamma\CA}(\mu)=\sum M^{\Gamma\ba}(\mu),\; \ba
\text{ independent subarrangement  of } \CA.$$
\end{theorem}
\begin{proof}
We use induction on the number $N$ of elements in $\CA$. If $\CA$
is linearly independent, we are done. If not, we assume that the
statement is known for arrangements with $N-1$ elements, and write
$\CA= \{H_{\alpha_1},H_{\alpha_2},\ldots, H_{\alpha_N}\}$. As
$\CA$ is not independent, there is a hyperplane, say
$H_{\alpha_N}$, which is  dependent on the rest of the system
$$\CA'=\{H_{\alpha_1},\ldots,H_{\alpha_{N-1}}\}.$$
 For $1\leq i\leq (N-1)$ we let
$$\CA'_i=\{H_{\alpha_1},H_{\alpha_2},\ldots,\hat H_{\alpha_i},\ldots,H_{\alpha_N}\}.$$
Note that  each $\CA'_i$ has $N-1$ elements. A  function $F$ of
$M^{\Gamma\CA}(\mu)$ may be written in the form  $F=P/D$ with
$P=\sum_{\xi\in I} c_\xi e^{\xi}$ and $D=D' D_N$, where
 $$D'=\prod_{j=1}^{R}(1-u_je^{\beta_j})\quad\text{and}\quad
D_N=\prod_{j=1}^{n_N}(1-v_je^{\alpha_N}).$$ In the factorization
of $D'$, the elements $\beta_j$ belong to the set
$\{\alpha_1,\ldots,\alpha_{N-1}\}$.

Each of the terms $F_\xi={e^{\xi}}/{D}$ of $F$ is in
$M^{\Gamma\CA}(\mu)$. We may split $\mu$ as $\mu=\mu'+\mu_N$, with
$\mu_N=t_N\alpha_N$, $0\leq t_N\leq n_N$ and
$$\frac{e^{\xi}}{D'}\in M^{\Gamma\CA'}(\mu').$$
Applying the exchange lemma to $H_{\alpha_N}$ and the system
$\CA'$, we obtain an admissible decomposition of ${e^{\xi}}/{D'}$
as a sum of elements $F'_i\in M^{\Gamma \CA'_i}(\mu')$.  Then
$F_\xi$ is a sum of terms of the form $F'_i/{D_N}$, each of which
is in $M^{\Gamma \CA'_i}(\mu)$. Since the system $\CA'_i$ is
composed of $N-1$ hyperplanes, we may conclude the proof of the
Theorem by our induction hypothesis.
\end{proof}

\begin{remark}
  A fixed total order $\prec$ on the arrangement $\CA$ of hyperplanes in an 
  $n$-dimensional vector space selects a subset $NBC(\CA,\prec)$ of the set of
  $n$-tuples of  independent  hyperplanes in $\CA$. This subset is
  called  the {\em no-broken-circuit basis} of $\CA$
 (cf. \cite{sz2} for details). The arguments used in the proof of
  the above Theorem may be used to show that, in fact,
$$M^{\Gamma \CA}(\mu)=\sum_{\ba\in
NBC(\CA,\prec)}M^{\Gamma\ba}(\mu).$$
Moreover, the sets $NBC(\CA,\prec)$ are minimal with respect to this property.
\end{remark}

Now we analyze the set $M^{\Gamma \ba}(\mu)$ when the arrangement
$\ba$ is independent. Thus let $\ba$ be a set of $m$ independent
hyperplanes. We choose $\alpha_k\in \Gamma^*$, $k=1,\ldots,m$, such
that $\ba=\{H_{\alpha_1},\ldots,H_{\alpha_m}\}$. Then
$\phi=[\alpha_1,\ldots,\alpha_m]$ is a sequence of linearly
independent linear forms.  Let ${\bf h}=[h_1, h_2,\ldots, h_m]$ be
a sequence of nonnegative integers and ${\bf u}= [u_1,u_2,\ldots
,u_m]$ be a sequence of nonzero complex numbers. We introduce the
function
$$g(\xi,\phi, {\bf h},{\bf u})=\frac{e^{\xi}}{\prod_{i=1}^m (1-u_i
e^{\alpha_i})^{h_i}},$$
where $\xi\in \Gamma^*$.

\begin{proposition}\labelm{prop:simplest}
For an independent arrangement $\ba=\CA(\phi)$, each function $F\in
M^{\Gamma \ba}(\mu)$ may be represented as a $\mu$-admissible linear
combination of the functions $g(\xi,\phi,{\bf h},{\bf u})$.
\end{proposition}
\begin{proof}
Clearly, it is sufficient to prove this statement for the case $|\phi|=1$.
 The general case will follow by taking the product of the
 linear combinations
 for each participating linear form.

 Set $\phi=\{\alpha\}$ and $\ba=\CA(\phi)$.
An element $F\in M^{\Gamma \ba}(\mu)$ is a  linear
combination of elements $F_\xi={e^{\xi}}/{D}\in M^{\Gamma
\ba}(\mu)$, where $D=\prod_{i=1}^R(1-u_i e^{\alpha})$. We need to
show that each function $F_\xi$ may be represented as a linear
combination of elements of $M^{\Gamma
  \ba}(\mu)$ of the form ${e^{\zeta}}/{(1-ve^{\alpha})^h}$.

We use induction  on $R$. If all the $u_i$ are equal, then
$F_\xi$ already has the required form. If not, up to reordering,
we can assume that $u_1\neq u_2$. We write $D=D_{12}D'$ with
$D_{12}=(1-u_1e^{\alpha})(1-u_2e^{\alpha})$.
Factor $F_\xi$ as  $F_\xi=G/{D_{12}}$ with
$G={e^{\xi}}/{\prod_{i=3}^R(1-u_i e^{\alpha})}$, and let
$\mu=\mu'+\mu''$, where $\mu'=t\alpha$ with $0\leq t\leq 2$ and
$\mu''\in \Box(G)$.

There are two cases: if $0\leq t\leq 1$, we write
$$\frac{1}{(1-u_1e^{\alpha})(1-u_2e^{\alpha})}=\frac{1}{(1-u_2
u_1^{-1})}\frac{1}{(1-u_1e^{\alpha})}+\frac{1}{(1-u_1
u_2^{-1})}\frac{1}{(1-u_2 e^{\alpha})};$$
if $1\leq t\leq 2$, we write
$$\frac{1}{(1-u_1e^{\alpha})(1-u_2e^{\alpha})}=\frac{-1}{(1-u_2
  u_1^{-1})}\frac{u_1^{-1}e^{-\alpha}}{(1-u_2e^{\alpha})}+\frac{-1}{(1-u_1
  u_2^{-1})}\frac{u_2^{-1}e^{-\alpha}}{(1-u_1 e^{\alpha})}.$$
In both
cases, we obtain a $\mu$-admissible decomposition of
${e^{\xi}}/{D}$ into a sum  $G_1+G_2$, where
$$G_1=c_1\frac{e^{\xi'}}{\prod_{i\neq 1}(1-u_i e^{\alpha})}
\text{ and } G_2=c_2\frac{e^{\xi'}}{\prod_{i\neq 2}(1-u_i
e^{\alpha})}.$$ This allows us to reduce the number of
factors in $F_\xi$  by one. Our statement now follows by the inductive
hypothesis.
\end{proof}

\subsection{Essential arrangements and nonspecial elements}
\labelm{sec:essential}

Now we formulate a version of Theorem \ref{thm:parfrac} in a form
which incorporates Proposition \ref{prop:simplest} and excludes
some degenerate cases.

Let again $\Gamma$ be a lattice of full rank in the
$n$-dimensional vector space $V$, and let $\CA=\{H_1,\ldots,
H_N\}$ be an essential $\Gamma$-rational arrangement of
hyperplanes in $V$. Fix a set $\Delta$ of representative linear
forms for $\CA$; thus we have $\CA=\CA(\Delta)$.  Define $\mu\in
V^*$ to be $\Gamma$-{\em special} with respect to $\CA$ if $\mu =
\lambda+\sum_{i=1}^N t_i \alpha_i$, where $\lambda\in\Gamma^*$,
$t_i\in\R$, $\alpha_i\in \Delta$ and at most $n-1$ of the
coefficients $t_i$ are nonzero. This property depends both on
$\Gamma$ and on $\CA$. The set of nonspecial elements is a
$\Gamma^*$-invariant union of open polyhedral chambers in $V^*$.

Note that if $F\in M^{\Gamma\CA}$, then the boundary of $\Box(F)$
is contained in the set of special elements. Thus if $\mu$ is
nonspecial and $F\in M^{\Gamma\CA}(\mu)$, then $\mu$ is in the
interior of $\Box(F)$. We arrive at the following
proposition.
\begin{proposition}\labelm{prop:simplestbasis}
  Let $\CA$ be an essential arrangement of rational hyperplanes. Let
  $F\in M^{\Gamma\CA}$. Let $\mu\in \Box(F)$ be a nonspecial element.
  Then there exists a set $\CB$ consisting of independent
  $n$-tuples of hyperplanes and a $\mu$-admissible decomposition
  $F=\sum_{\ba \in \CB} F_{\ba}$, with $F_{\ba}\in M^{\Gamma\ba}(\mu)$.
  Furthermore, choosing a basis $\sigma=\{\alpha_1,\ldots,\alpha_n\}$
  such that $\ba=\CA(\sigma)$, the function $F_\ba$ is a linear
  combination of functions $g(\xi,\sigma, {\bf h},{\bf u})$ with
  $\mu+\xi=\sum_{i=1}^n t_i\alpha_i$ and $0< t_i < h_i$.
\end{proposition}

This proposition allows us to write $F$ as a linear combination of
those functions  $g(\xi,\sigma, {\bf h},{\bf u})$, for which $\mu$ belongs
to the interior of $\Box(g(\xi,\sigma,{\bf h}, {\bf u}))$.

\section{Expansion and inversion formula}
\labelm{sec:expansion}

\subsection{Expansion of functions}
\labelm{sec:expfus}

Let $V$ be a real vector space of dimension $n$ endowed with a lattice
$\Gamma$. A choice of a nonzero vector $v$ in $V$ induces a choice of
an open  half space $V^*_+\subset V^*$ of linear forms which take positive
values on $v$. We fix such a half space and  consider a finite subset
$\Delta$ of elements of $V^*_+\cap\Gamma^*$. We
assume that $\Delta$ linearly spans the vector space $V^*$ and thus
generates a closed acute $n$-dimensional cone
$C(\Delta)$:
$$C(\Delta) =\sum_{\alpha\in \Delta}\R_{\geq0}\alpha.$$

Recall that we denoted by $\bdelta$ the set of those subsets of
$\Delta$ which are bases of $V^*$. Following \cite{agz}, we will
call a vector in $V^*$ {\em singular} with respect to $\Delta$ if
it is in a cone $C(\nu)$ generated by a subset $\nu\subset\Delta$
of cardinality strictly less than $n$. The set of singular vectors
will be denoted by $\vsing$ and the vectors in the complement
$\vreg=V^*\setminus\vsing$ will be called {\em regular}. The
connected components of $\vreg$ are conic chambers called {\em big
chambers}.  This term is chosen to differentiate them from the
smaller chambers cut out by special elements defined in section
\ref{sec:essential}. We might call big chambers simply {\em chambers}, whenever this does not  cause confusion.

A big chamber is an open cone. Note that there might be
regular elements which are {\em special} in the sense of paragraph
\ref{sec:essential}: if such $\mu$ is written $\mu =\sum_{i=1}^N
t_i \alpha_i$, $t_i\in\R$, $\alpha_i\in \Delta$ and at most $n-1$
of the coefficients $t_i$ are nonzero, then at least one of the
coefficients $t_i$ is strictly negative.

\begin{figure}[ht]
\begin{center}

\setlength{\unitlength}{0.00062500in}
\begingroup\makeatletter\ifx\SetFigFont\undefined%
\gdef\SetFigFont#1#2#3#4#5{%
  \reset@font\fontsize{#1}{#2pt}%
  \fontfamily{#3}\fontseries{#4}\fontshape{#5}%
  \selectfont}%
\fi\endgroup%
{\renewcommand{\dashlinestretch}{30}
\begin{picture}(3963,3360)(0,-10)
\path(2025,3075)(1950,2100)
\path(1350,2550)(3075,1125)
\path(675,2025)(2550,2100)
\path(2025,3075)(1950,2100)
\path(675,2025)(2025,3075)(3075,1125)(675,2025)
\path(675,2025)(675,75)
\path(675,75)(3075,1125)
\dashline{60.000}(1425,1725)(2025,3075)
\path(1350,2550)(2550,2100)
\dashline{60.000}(675,75)(1425,1725)
\put(1800,1800){\makebox(0,0)[lb]{\smash{{{\SetFigFont{9}{10.8}{\rmdefault}{\mddefault}{\updefault}$\frac{1}{3}e_1$}}}}}
\put(3150,1125){\makebox(0,0)[lb]{\smash{{{\SetFigFont{9}{10.8}{\rmdefault}{\mddefault}{\updefault}$e_3$}}}}}
\put(1725,3150){\makebox(0,0)[lb]{\smash{{{\SetFigFont{9}{10.8}{\rmdefault}{\mddefault}{\updefault}$e_2-e_3$}}}}}
\put(2625,2100){\makebox(0,0)[lb]{\smash{{{\SetFigFont{9}{10.8}{\rmdefault}{\mddefault}{\updefault}$\frac{1}{2}e_2$}}}}}
\put(375,2625){\makebox(0,0)[lb]{\smash{{{\SetFigFont{9}{10.8}{\rmdefault}{\mddefault}{\updefault}$\frac{1}{2}(e_1-e_3)$}}}}}
\put(0,2100){\makebox(0,0)[lb]{\smash{{{\SetFigFont{9}{10.8}{\rmdefault}{\mddefault}{\updefault}$e_1-e_2$}}}}}
\put(450,0){\makebox(0,0)[lb]{\smash{{{\SetFigFont{9}{10.8}{\rmdefault}{\mddefault}{\updefault}$0$}}}}}
\end{picture}
}
\caption{The Chambers of the system $A_3^+$}
\end{center}
\end{figure}
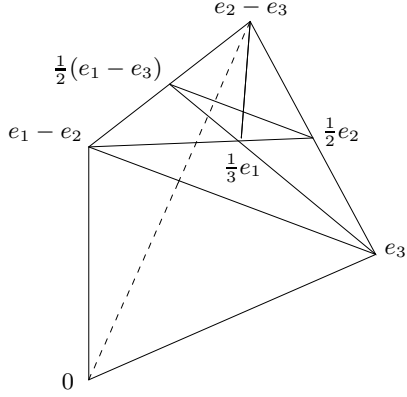

\vspace{10mm}

If $\c$ is a big chamber and $\sigma\in\bdelta$, then either
$\c\subset C(\sigma)$ or $\c\cap C(\sigma)=\emptyset$.  One of the big
chambers is the complement of the closed cone $C(\Delta)$; we denote
it by $\cnull$. Note that this convention is slightly different from
the convention adopted in \cite{bav}, where $\cnull$ 
was not considered a chamber.

If $\c$ is a big chamber contained in $C(\Delta)$,
then the closure of $\c$ may be represented as
$$ \overline{\c} = \bigcap C(\sigma),\quad \c\subset
C(\sigma),\,\sigma\in\bdelta.$$ In particular $\overline{\c}$ is a
closed convex polyhedral cone.

Denote by $\C[[\Gamma^*]]$ the set of complex, formal, possibly
infinite linear combinations of the exponentials $e^{\lambda}$,
where $\lambda\in\Gamma^*$.  If $\Theta=\sum_{\lambda\in
\Gamma^*}m_\lambda e^{\lambda}$ is an element of $\C[[\Gamma^*]]$,
then the {\em support} of $\Theta$ is the set of
$\lambda\in\Gamma^*$ such that $m_\lambda\neq 0$. The coefficient
$m_\lambda$ of $e^{\lambda}$ in $\Theta$ will be denoted by
$\Coeff(\Theta,\lambda)$.

Let $\C_\Delta[[\Gamma^*]]$ be the subspace of
$\C[[\Gamma^*]]$ spanned by the elements $\Theta$ with supports
contained in  sets of the form  $I+C(\Delta)$,  where $I$ is a finite subset
of $\Gamma^*$. This subspace forms a ring which contains the ring
$\C[\Gamma^*]$ of finite linear combinations of elements $e^{\xi}$,
$\xi\in \Gamma^*$.

Consider the arrangement of hyperplanes $\CA=\CA(\Delta)$ and recall
the definition of the algebra $M^{\Gamma\CA}$ from \S1.   Every
function $F\in M^{\Gamma\CA}$ can be written in the form
$$F=\frac{\sum_{\xi\in I} c_\xi e^{\xi}} {\prod_{k=1}^R(1-u_k e^{
\beta_k})},$$
where $I$ is a finite subset of $\Gamma^*$,
$u_k,c_\xi\in \C^*$, and the exponents $\beta_k$ are in $\Delta$.

For $\alpha\in \Delta$ and $u\in \C^*$, define the expansion
$$r^+\left(\frac{1}{1-u e^{\alpha}}\right)=\sum_{k=0}^{\infty} u^k
e^{k\alpha},$$
where the right hand side is interpreted as a formal series.
This expansion map extends to an injective ring homomorphism
 $$r^+:M^{\Gamma\CA}\to \C_\Delta[[\Gamma^*]]$$
given by $$r^+(F)=\left(\sum_{\xi\in I} c_\xi e^{\xi}\right)\prod_{k=1}^R
r^+\left(\frac{1}{1-u_k e^{\beta_k}}\right).$$ We call $r^+(F)$ the expansion
of $F$.

The  aim  of this section is to give a residue formula for the
coefficient $\Coeff(r^+(F), \lambda)$ for $F\in M^{\Gamma\CA}$
and $\lambda\in \Gamma^*$.

\subsection{The residue transform}
\labelm{sec:residue}

We start with a general definition of exponential-polynomial functions.
\begin{definition}
\begin{itemize}
  \item For a $\Z$-module $W$ and a field $\F$, define the space of
  polynomial functions $P(W,\F)$ to be the subring of
  $\F$-valued functions on $W$ generated by the additive
  $\F$-valued characters of $W$.
  \item For a $\Z$-module $W$ and a field $\F$, define the space of
  exponential-polynomial functions $\EP(W,\F)$ to be the subring of
  $\F$-valued functions on $W$ generated by the additive $\F$-valued
  and multiplicative $\F^*$-valued characters of $W$.
\end{itemize}
\end{definition}
Clearly, an exponential-polynomial function is a linear combination of
multiplicative characters  (exponentials) with polynomial coefficients.
Usually, we will set $\F=\C$, and in this case we will write $\EP(W)$ for
$\EP(W,\F)$. In our applications, $W$ will be either a vector
space or a lattice.

When $W$ is a lattice of full rank in a vector space $E$, a polynomial
function $f$ on $W$ extends in an unique way to a polynomial function
on $E$. Exponential-polynomial functions also extend to
exponential-polynomial functions on $E$, but the extension is not
unique. For example, if $W=\Z\subset \R=E$, then the function
$n\mapsto (-1)^n n$ is an exponential-polynomial function on $\Z$,
which can be extended on $\R$ as the exponential-polynomial function
$x\mapsto e^{i(2k+1)\pi x} x$ for any integer $k$.

When $W$ is a lattice and a function $f\in \EP(W)$ is such that the
multiplicative characters which appear in it take values in roots
of unity, then such a function is called {\em periodic-polynomial}
or sometimes, {\em quasipolynomial}.

We continue with the setup of a lattice $\Gamma\subset V$, an
arrangement $\CA$ and a set of linear forms
$\Delta\subset\Gamma^*$ representing $\CA$. In this section we
associate to any $F\in M^{\Gamma\CA}$ an exponential-polynomial
function on $\Gamma^*$ with values in the space of simple
fractions $S_{\CA}$.

According to Lemma \ref{poles}, the total residue of a function $F\in
M^{\Gamma\CA}$, written in the form
\begin{equation}
  \label{eq:fform}
F=\frac{\sum_{\xi\in I}
c_\xi e^{\xi}} {\prod_{k=1}^R(1-u_k e^{ \beta_k})},
\end{equation}
vanishes unless the set of linear forms $\{\beta_k|\,u_k=1\}$ spans
the vector space $V^*$. Let us define the total residue of $F$ at some
point $p\in\vc$ as the total residue of the function $z\mapsto
F(z-p)$. Then we observe that the total residue of $F$ given in the
above form vanishes at $p\in\vc$ unless
$$
\text{the set of forms } \{\beta_k|\,
e^{\langle\beta_k,p\rangle}u_k=1\} \text{ linearly spans }V^*.$$
The linear forms $\beta_k$ are all in $\Gamma^*$, hence the set
$\SP(F,\Gamma)$ of those points $p\in\vc$ which satisfy this condition
is invariant under translations by elements of the lattice $2\pi i\Gamma$.
Consider two points in $\vc$ equivalent if they are related by such a
translation, and choose a set $\RSP(F,\Gamma)$ containing exactly one
point from each equivalence class of points in $\SP(F,\Gamma)$. It is
clear from the definitions that the set $\RSP(F,\Gamma)$ is finite; we
will call it a {\em reduced set of poles} of $F$.

This definition of the set $\RSP(F,\Gamma)$ is somewhat
informal: it depends on the presentation of $F$. The only properties
that we will need from it are that
\begin{itemize}
\item the set $\RSP(F,\Gamma)$ is finite;
\item if $p,q\in \RSP(F,\Gamma)$ and $p-q\in2\pi i\Gamma$, then $p=q$;
\item if the total residue of $F(z)G(z)$, where $G(z)$ is an entire
  function, does not vanish at some $q\in\vc$, then $q\in2\pi
  i\Gamma+\RSP(F,\Gamma)$.
\end{itemize}

Now we define a function $s[F,\Gamma]:\Gamma^*\rightarrow S_\CA$
with values in the space of simple fractions associated to
$\Delta$, whose value at $\lambda\in \Gamma^*$ is the sum of all
the total residues of the function $z\mapsto e^{\langle
\lambda,z\rangle}F(-z)$ taken at inequivalent points in $\vc$.
More precisely,
\begin{definition}
\labelm{def:s} For $F\in M^{\Gamma\CA}$ and $\lambda\in \Gamma^*$,
we introduce
\begin{equation}
  \label{eq:defs}
s[F,\Gamma](\lambda)=\sum_{p\in  \RSP(F,\Gamma)} \Tres(e^{\langle
\lambda,z-p\rangle}F(p-z)),
\end{equation}
 where the set
$\RSP(F,\Gamma)$ is  a reduced set of poles of $F$.
\end{definition}
Clearly, the definition does not depend on the choice of
representatives $\RSP(F,\Gamma)$.
\begin{lemma}
The function $\lambda\mapsto s[F,\Gamma](\lambda)$ is an
exponential-polynomial function on $\Gamma^*$ with values in the
space of simple fractions $S_\CA$.
\end{lemma}
 \begin{proof}
Let $F\in M^{\Gamma\CA}$ be given in the form \eqref{eq:fform}, and
pick an element $p\in\RSP(F,\Gamma)$.  For
$\lambda\in V^*$  consider the total residue
$$\Tres\left(e^{\langle\lambda,z\rangle}F(p-z)\right).$$

The function $e^{\langle\lambda,z\rangle}F(p-z)$ is in the space
$\hatr$ introduced at the end of section \ref{sec:cxhpa}.  As the
total residue depends only on the component of degree $-n$ of this
function, for the purpose of the calculation of its total residue, we
can replace the exponential $e^{\langle \lambda,z\rangle}$ by its
expansion truncated up to order $R-n$. Thus we have
$$\Tres\left(e^{\langle\lambda,z\rangle}F(p-z)\right)=
\Tres\left(\sum_{j=1}^{R-n}\frac{\langle\lambda,z\rangle^j}{j!}F(p-z)\right).$$
The right hand side here clearly depends polynomially on
$\lambda$, thus each term
$$
\Tres\left(e^{\langle\lambda,z-p\rangle}F(p-z)\right)=
e^{-\langle\lambda,p\rangle}
\Tres\left(e^{\langle\lambda,z\rangle}F(p-z)\right),$$ appearing
in the definition of $s[F,\Gamma]$ is an exponential-polynomial
function of $\lambda$. As the set $\RSP(F,\Gamma)$ is finite, this completes
the proof.
\end{proof}

Let us look at a few special cases.

\noindent {\bf 1}.  Let $F$ be of the form
$$F=\frac{\sum_{\xi\in I} c_\xi e^{\xi}} {\prod_{k=1}^R(1- e^{
    \beta_k})},$$
i.e.  let all constants $u_k$ be  equal to $1$. For a basis
$\sigma=\{\alpha_1,\ldots,\alpha_n\}$ of $V^*$, formed by elements
of the sequence  $[\beta_1,\ldots,\beta_R]$, the lattice $\Z\sigma$
is contained in $\Gamma^*$ and is usually different from
$\Gamma^*$. Consider $(\Z\sigma)^*\subset V$, the dual lattice to
$\Z\sigma$:
$$
(\Z\sigma)^*=\{s\in V|\,\langle s,\alpha_k\rangle\in \Z\text{ for
  }1\leq k\leq n\}.$$
If $p\in 2\pi i (\Z\sigma)^*$, then the set of linear forms
$\{\beta_k|\, e^{\langle\beta_k,p\rangle}=1\}$ linearly spans
$V^*$, since it contains $\sigma$. Then the set $\RSP(F,\Gamma)$
is a union of representatives of the finite groups $2\pi
i(\Z\sigma)^*/2\pi i\Gamma$ in $V_\C$, as $\sigma$ varies over
bases of $V^*$ formed by the $\beta_k$s.  As a result, for
$\lambda\in\Gamma^*$ and $p\in \RSP(F,\Gamma)$, the exponential
$e^{\langle\lambda,p\rangle}$ is a root of unity. This implies
that the function $s[F,\Gamma]$ on the lattice $\Gamma^*$ is {\em
periodic-polynomial}. More precisely, if $n_F$ is an integer such
that $n_F\Gamma\subset \Z\sigma$ for all bases $\sigma$ of $V^*$
formed by  $\beta_k$s, then the function $s[F,\Gamma]$ is {\em
polynomial} on all cosets of the form $\lambda+n_F\Gamma^*$.

\noindent {\bf 2}. There is an interesting special case of this
setup, when $s[F,\Gamma]$ is plainly polynomial: the unimodular case.

We will call a subset $\Delta\subset\Gamma^*$ {\em unimodular}, if
 every basis $\sigma\in\bdelta$ is a $\Z$-basis of $\Gamma^*$, i.e.
  the parallelepiped
$\sum_{\alpha\in\sigma}[0,1]\alpha$ contains no elements of the
lattice $\Gamma^*$ in its interior.

In this case, the integer $n_F$ mentioned above may be taken to be
equal to 1. We collect what we have found in the following
\begin{lemma}
Let
$$F=\frac{\sum_{\xi\in I} c_\xi e^{\xi}} {\prod_{k=1}^R(1-
  e^{\beta_k})}.$$
Then the function $\lambda\mapsto
s[F,\Gamma](\lambda)$ is a periodic-polynomial function on $\Gamma^*$.
If, furthermore, the elements of the sequence $[\beta_1,\beta_2,\ldots
,\beta_R]$ belong to a unimodular subset of $\Gamma^*$, then the
function $\lambda\mapsto s[F,\Gamma](\lambda)$ is a polynomial.
\end{lemma}

\noindent{\bf 3}.  Assume,  at the other extreme,  that the
constants $u_k$ are generic.

For $p\in \RSP(F,\Gamma)$, denote by  $\bfj(p)$ the subset of  the
set of indices $\{1,2,\ldots,R\}$ given by
$$
\bfj(p) = \{j\in\{1,2,\ldots,R\}|\,
u_je^{\langle\beta_j,p\rangle}=1\}.$$ If the constants $u_k$ are
generic, then the set $\{\beta_j\}, j\in \bfj(p)$, if nonempty,
consists of exactly $n$ linearly independent elements of $\Delta$.

We have
\begin{multline*}
e^{\langle\lambda,z-p\rangle}  F(p-z)=e^{\langle\lambda,z-p\rangle
} \left(\sum_{\xi\in I} c_\xi
    e^{\langle\xi,p-z\rangle}\right)\frac{1}{\prod_{j\in
      \bfj(p)}\langle\beta_j,z\rangle} \\
  \times\prod_{j\in\bfj(p)}
\frac{\langle\beta_j,z\rangle}{1-e^{-\langle\beta_j,z\rangle}}\prod_{k\notin
    \bfj(p)} \frac{1}{1-u_k e^{\langle\beta_k,p-z\rangle}}.
\end{multline*}
By Lemma \ref{poles}, the total residue of this function is the
simple fraction $\prod_{j\in {\bf j}(p)}\beta_j^{-1}$, multiplied
by the constant, which is obtained by setting $z$ to zero in the
rest of the expression. As a result we obtain the following
explicit formula:
$$
s[F,\Gamma](\lambda)= \sum_{p\in
\RSP(F,\Gamma)}e^{-\langle\lambda,p\rangle} \left(\sum_{\xi\in I}
c_\xi e^{\langle\xi,p\rangle}\right)\prod_{k\notin \bfj(p)}
\frac{1}{1-u_k e^{\langle\beta_k,p\rangle}}\times\frac{1}
{\prod_{j\in \bfj(p)}\beta_j},$$ which expresses the function
$\lambda\mapsto s[F,\Gamma](\lambda)$ as a linear combination of
exponentials.

\subsection{The residue formula}
\labelm{sec:resform}

We start with recalling the notion of residue introduced by
Jeffrey and Kirwan \cite{jk}. Let  again $\Delta$ be  a set of
vectors in an open  halfspace of an $n$-dimensional real vector
space $V^*$ and let $\CA=\CA(\Delta)$. We assume that $\Delta$
generate $V^*$. Fix a volume form $\vol$ on $V^*$. Given a big
chamber $\c$ of $\vreg$, one can construct a functional $f\mapsto
\jk\c f\vol$ on the space $S_\CA$ of simple fractions as follows.
For a simple fraction
$$f_{\sigma}=\frac{1}{\prod_{\alpha\in \sigma}\alpha}, \quad \sigma\in\bdelta,$$
 set
$$ \jk{\c}{f_\sigma}{\vol} =
\begin{cases}
{\vol(\sigma)}^{-1},&\text{if }\c\subset C(\sigma),\\
0, &\text{if }\c\cap C(\sigma)=\emptyset.
\end{cases}
$$
Here we denoted by $\vol(\sigma)$ the volume of the parallelepiped
$\sum_{\alpha\in\sigma} [0,1]\alpha$ with respect to  our chosen volume
form.

Now we  formulate our main result. Let $\Gamma$ be a rank-$n$
lattice in $V$, and let $\Gamma^*\subset V^*$ be its dual lattice.
As before, we assume that $\Delta\subset\Gamma^*$. Denote by
$\vol_{\Gamma^*}$ the measure on $V^*$ assigning volume 1 to a
minimal parallelepiped spanned by elements of $\Gamma^*$.

\begin{theorem}\labelm{thm:main}
Let $F\in M^{\Gamma\CA}$ and let $\Box(F)^{0}$ be the interior of
$\Box(F)$. Then for $\lambda\in \Gamma^*$ and any  big chamber
$\c$ such that $ (\lambda+\Box(F)^{0})\cap \c\neq\emptyset$, one has

\begin{equation}
  \label{eq:main}
  \Coeff(r^+(F),\lambda)=\jk{\c}{s[F,\Gamma](\lambda)}{\vol_{\Gamma^*}}.
\end{equation}
\end{theorem}
Before starting the proof, we analyze the 1-dimensional case. Let
$V=\R e$ and  $V^*=\R e^*$ with lattices $\Gamma=\Z e$ and
$\Gamma^*=\Z e^*$; let $\Delta=\{e^*\}$.  There are two chambers
in this case: $\c^+=\R_{>0} e^*$ and $\c^-=\R_{<0}e^*$.  We simply
write $F(z)$ for a function $F(ze)$ on $V_\C$. Then
$\jk{\c^+}{\Tres F}{\volga}=\res_{z=0}F(z)\,dz$, while
$\jk{\c^-}{\Tres F}{\volga}=0$.

Introduce the notation
$$c(k,R)=\binom{k+(R-1)}{R-1}=\frac{1}{(R-1)!}(k+1)(k+2)\cdots
(k+(R-1)).$$
We have the following simple generating function for $c(k,R)$:
\begin{lemma}
$$\res_{z=0}\frac{e^{kz}}{(1-e^{-z})^R}dz=c(k,R).$$
\end{lemma}
\begin{proof}
Using the  change of variables $y=e^z$ in the calculation of the
residue, we obtain
\begin{multline*}
\res_{z=0}\frac{e^{kz}}{(1-e^{-z})^R}dz=
\res_{y=1}\frac{y^k}{(1-y^{-1})^R}\frac{dy}{y}\\
=\res_{y=1}\frac{y^{R+k-1}}{(y-1)^R}dy=
\res_{x=0}\frac{(1+x)^{R+k-1}}{x^R}dx=c(k,R).
\end{multline*}
\end{proof}

Now consider the function
$$F(z)=\frac{e^{\xi z}}{(1-ue^z)^R},$$ where $\xi$ is an integer.
The following explicit formula holds for the expansion of $F$:
$$r^+(F)=e^{\xi z}\sum_{k=0}^{\infty}c(k,R) u^k e^{kz}.$$
Hence we have
\begin{equation}
  \label{eq:coeff}
\Coeff(F,\lambda) =
\begin{cases}
  0,& \text{ if }\lambda-\xi\in\Z_{<0};\\
u^{\lambda-\xi}c(\lambda-\xi,R),&\text{ if }
\lambda-\xi\in\Z_{\geq0}.
\end{cases}
  \end{equation}

Note that the relation
$$\Coeff(F,\lambda)=u^{\lambda-\xi}c(\lambda-\xi,R)$$
holds whenever
$\lambda-\xi\geq -(R-1)$, since both sides of this equality vanish for
$\lambda-\xi=-1,-2,\ldots,-(R-1)$.

Let us analyze our proposed formula \eqref{eq:main} in this example. We
first write out the element $s[F,\Gamma](\lambda)$ explicitly. The
function $F(z)$ has just one pole $p$ modulo $2\pi i\Gamma$; it is  given by
the equation $e^p=u^{-1}$. Thus  we have
 $$s[F,\Gamma](\lambda)=u^{\lambda-\xi} \Tres
 \frac{e^{(\lambda-\xi)z}}{(1-e^{-z})^R},$$
which leads to
$$s[F,\Gamma](\lambda)=u^{\lambda-\xi}c(\lambda-\xi,R)\frac{1}{z}.$$

Now assume that  $\lambda\in \Z$ and  we picked a chamber $\c$ such that
$(\lambda+\Box(F)^0)\cap \c$ is not empty.

First we consider the case $\c=\c^+$. Here
$$\jk{\c^+}{s[F,\Gamma](\lambda)}{\volga}=u^{\lambda-\xi}c(\lambda-\xi,R).$$
Since $(\lambda+\Box(F)^0)\cap \c^+$ is nonempty, there exists $\mu\in
\Box(F)^0$ such that $\lambda+\mu>0$. As $\mu+\xi=t$ with $0<t<R$,
this implies that $\lambda-\xi\geq -(R-1)$. This is consistent with
our computation of $\Coeff(F,\lambda)$ above.

Assume now that $(\lambda+\Box(F)^0)\cap \c^-$ is not empty. Now we
have
$$\jk{\c^-}{s[F,\Gamma](\lambda)}{\volga}=0.$$

Since $(\lambda+\Box(F)^0)\cap \c^-$ is nonempty, there exists
$\mu\in \Box(F)^0$ such that $\lambda+\mu<0$. As $\mu+\xi=t$ with
$0<t<R$, this implies that $\lambda-\xi<0 $. Again, this is
consistent with \eqref{eq:coeff}.

\medskip

We now return to the Proof of the Theorem.

\begin{proof}
If $(\lambda+\Box(F)^{0})\cap \c$ is nonempty, then we can choose
a nonspecial element $\mu\in \Box(F)^0$ such that $\lambda+\mu\in
\c$.

By Proposition \ref{prop:simplestbasis}, there is a
$\mu$-admissible decomposition of  $F$ as a sum of functions
$g(\xi,\sigma, {\bf h},{\bf u})$ with  $\sigma\in\bdelta$.
Furthermore, the element $\mu$ still belongs to
$\Box(g(\xi,\sigma, {\bf h},{\bf u}))^0$. It is thus sufficient to
prove the theorem in the case $F=g(\xi,\sigma,{\bf h},{\bf u})$.

Let $\sigma=\{\alpha_1,\alpha_2,\ldots,\alpha_n\}$ and $\xi\in
\Gamma^*$. Then  we have
$$F(z)=\frac{e^{\langle\xi,z\rangle}}{(1-u_1e^{\langle\alpha_1,z\rangle})^{h_1}\cdots
 (1-u_ne^{\langle\alpha_n,z\rangle})^{h_n}}.$$

Let $C(\sigma)$ be the cone generated by $\sigma$ and $\Z\sigma$
the sublattice of $\Gamma^*$ generated by $\sigma$. Then
$C(\sigma)\cap \Z\sigma$ is the set of elements $\lambda\in V^*$
of the form $\lambda=\sum_{i=1}^n k_i\alpha_i$, where $k_i$ are
nonnegative integers. Then it easily follows from the result
\eqref{eq:coeff} in the 1-dimensional case that we have
\begin{equation}
  \label{eq:coeffn}
\Coeff(F,\lambda) =
\begin{cases}
  0,&\text{if }k_i<0\text{ for some }i, \,1\leq i\leq n,\\
u_1^{k_1}c(k_1,h_1)\cdots u_n^{k_n}c(k_n,h_n),&\text{if } k_i\geq
1-h_i,\,i=1,\ldots,n,
\end{cases}
  \end{equation}
 where $\lambda-\xi=\sum_{i=1}^n k_i\alpha_i$ and $k_i\in\Z$,
 $i=1,\ldots,n$.

We now compute $s[F,\Gamma](\lambda)$. The set of poles
$\SP(F,\Gamma)$ of the function $F$ is given by
\begin{equation}
  \label{eq:spdef}
\SP(F,\Gamma)=\{p\in V_\C|\, e^{\langle\alpha_k,p\rangle}=u_k^{-1},
\text{ for } 1\leq k\leq n\}.
\end{equation}
Choose an element $p_0$ in this set and again denote by $(\Z\sigma)^*$ the
dual lattice to $\Z\sigma$.
 Then for any $s\in(\Z\sigma)^*$, the point $p_0+2i\pi s$ of $V_\C$ still satisfies
$e^{\langle\alpha_k,p_0+2i\pi s\rangle}=u_k^{-1}$.
Thus $\SP(F,\Gamma)= p_0+2i\pi(\Z\sigma)^*$.

We have
$$
s[F,\Gamma](\lambda)= \sum_{p\in
\RSP(F,\Gamma)}e^{\langle\xi-\lambda,p\rangle}\;
\Tres\frac{e^{\langle\lambda-\xi,z\rangle}}{\prod_{k=1}^n
(1-e^{-\langle\alpha_k,z\rangle})^{h_k}}.
$$

Since $\RSP(F,\Gamma)$ is a set of representatives of the set
$\SP(F,\Gamma)$ modulo the lattice $2\pi i\Gamma$, using
\eqref{eq:spdef}, we can write
$$\sum_{p\in \RSP(F,\Gamma)} e^{\langle\xi-\lambda,p\rangle}=
e^{\langle\xi-\lambda,p_0\rangle}\sum_{m\in (\Z\sigma)^*/ \Gamma}e^{2i\pi
  \langle \xi-\lambda,m\rangle}.$$

This sum is nonzero if and only if $\xi-\lambda\in \Z\sigma$.  If
$\lambda-\xi=\sum_{i=1}^n k_i\alpha_i$ with $k_i\in \Z $, then
$$e^{\langle\xi-\lambda,p_0\rangle}= e^{-\sum_{i=1}^n k_i\langle\alpha_i,p_0\rangle}=
u_1^{k_1}\cdots u_n^{k_n}.$$
We thus obtain:
\begin{itemize}
\item
$s[F,\Gamma](\lambda)$ is equal to $0$ if $\lambda-\xi$ is not in
$\Z\sigma$;
\item  If $\lambda-\xi=\sum_{i=1}^n k_i\alpha_i$ with $k_i\in \Z $, then
 $$s[F,\Gamma](\lambda)=|\Gamma^*/\Z \sigma|  \prod_{i=1}^n u_i^{k_i}
 c(k_i,h_i)\times\frac{1}{\prod_{i=1}^n\alpha_i}.$$
\end{itemize}

Now we proceed to computing the Jeffrey-Kirwan residues.
Let $\lambda\in \Gamma^*$ and let $\c$ be a chamber such that
$(\lambda+ \Box(F)^0)\cap \c\neq \emptyset$. This means that there is $\mu\in
V^*$ such that $\lambda+\mu\in\c$ and $\mu+\xi=\sum_{i=1}^n t_i\alpha_i$ with
$0<t_i<h_i$.

There are two cases: either $\c\subset C(\sigma)$ or $\c\cap
C(\sigma)=\emptyset$. If $\c\subset C(\sigma)$, then we can
conclude that $\lambda-\xi=\sum_{i=1}^n x_i\alpha_i$, where $x_i$
are rational numbers and $x_i> -h_i$, $i=1,\ldots,n$. On the other
hand, we have
\begin{itemize}
\item  $\jk\c{s[F,\Gamma](\lambda)}{\volga}=0$ if  $\lambda-\xi\notin\Z \sigma.$
\item  If $\lambda-\xi=\sum_{i=1}^n k_i\alpha_i$ with $k_i\in\Z$,
then $\jk\c{s[F,\Gamma](\lambda)}\volga$ factors:
$$|\Gamma^*/\Z \sigma|
\prod_{i=1}^nu_i^{k_i}
c(k_i,h_i)J\left\langle\c,\frac{1}{\prod_{i=1}^n\alpha_i}\right\rangle_{\volga}
=\frac{|\Gamma^*/\Z \sigma|}{\volga(\sigma)} \prod_{i=1}^nu_i^{k_i}
c(k_i,h_i).$$
\end{itemize}
It easy to see from the definitions that $\volga(\sigma)=|\Gamma^*/\Z
\sigma|$, hence
$$\jk\c{s[F,\Gamma](\lambda)}\volga=\prod_{i=1}^nu_i^{k_i}
c(k_i,h_i).$$
This is consistent with the expression \eqref{eq:coeffn} for $\Coeff(F,\lambda)$, as
$\lambda-\xi=\sum_{i=1}^n k_i\alpha_i$, where $k_i$ are integers and
$k_i\geq 1-h_i$.

In the case $\c\cap C(\sigma)=\emptyset$, we can conclude that
$\lambda-\xi=\sum_{i=1}^n x_i\alpha_i$, where at least one of the
numbers $x_i$ is negative. This is again consistent with
\eqref{eq:coeffn}, since by definition
$\jk\c{s[F,\Gamma](\lambda)}\volga=0$.

Thus we covered all cases and the theorem is proved.
\end{proof}

\section{ Ehrhart polynomials}
\labelm{sec:ehrhart}

\subsection{Partition polytopes and the vector partition function}
\labelm{sec:partpol}

Let $V$ be a real vector space of dimension $n$ endowed with a
lattice $\Gamma$, and let $\Phi$ be a sequence of {\em not
necessarily distinct} elements $[\beta_1,\ldots,\beta_N]$ of the
dual lattice $\Gamma^*\subset V^*$. We assume that $\Phi$
generates $V^*$. Denote by $\rho$ the surjective linear map from
$\R^N$ to the vector space $V^*$ defined by $\rho(w_k):=\beta_k$,
$1\leq k\leq N$, where $\{w_k\}_{k=1}^N$ is the standard basis of
$\R^N$.

The map $\rho$ may be written as
$$\rho(x_1,x_2,\ldots,x_N)=\sum_{i=1}^N x_i\beta_i.$$

 We
 denote by $C_N^+$ the closed convex cone in $\R^N$ generated by
 $w_1,\ldots,w_N$, and we set $C(\Phi):=\rho(C_N^+)$, the cone
 generated by $(\beta_1,\ldots,\beta_N).$ We assume here that
 $\rho^{-1}(0)\cap C_N^+=\{0\}$. Then $0$ is not in the convex hull of
 the vectors $\beta_k$ and $C(\Phi)$ is an acute cone.

 \begin{definition}
   For $a\in V^*$, we define the {\em partition polytope}
   $\Pi_{\Phi}(a)$ by
$$\Pi_{\Phi}(a):=\rho^{-1}(a)\cap    C_N^+~.$$
 \end{definition}
 The set $\Pi_{\Phi}(a)$ is the convex polytope consisting of all
 solutions $(x_1,x_2, \ldots, x_N)$ of the equation $$\sum_{k=1}^N x_k
 \beta_k=a$$
 in nonnegative real numbers $x_k$. In particular, the
 polytope $\Pi_{\Phi}(a)$ is empty if $a$ is not in the cone
 $C(\Phi)$.

 When $V=\R e$ is one-dimensional, and $\Delta=[b_1
 e^*,\ldots,b_Ne^*]$ where $b_k$ are positive integers,
  the polytope $\Pi_{\Phi}(a)$ is the
 $(N-1)$-dimensional simplex consisting of the intersection of the
 hyperplane $\sum_{i=1}^N b_i x_i=a$ with the positive quadrant.

 Denote by $\Z \Phi$ the lattice in $V^*$ generated by $\Phi$;
 naturally $\Z \Phi\subset\Gamma^*$. Then the map $\rho$ sends the
 standard lattice $\Z^N\subset\R^N$ to the lattice $\Z \Phi$.

For a general $\lambda$ in the lattice $\Gamma^*$, the vertices of
the polytope $\Pi_{\Phi}(\lambda)$ are only rational rather than
integral.

\begin{example}
  Set $V=\R e$ with $\Gamma=\Z e$, $\beta_1=2 e^*$ and $\beta_2= 3
  e^*$. Let $\lambda$ be a nonnegative integer. Then the polytope
  $\Pi_{\Phi}(\lambda e^*)$ consists of the set $\{(x_1,x_2)|\, x_1\geq
  0, x_2\geq 0, 2 x_1+3x_2=\lambda\}$. The vertices of $\Pi(\lambda)$
  are $(\frac{\lambda}{2},0)$ and $(0,\frac{\lambda}{3})$, so they are
  integral if and only if $\lambda$ is multiple of $6$.
\end{example}

Let $\Delta=\{\alpha_1,\ldots,\alpha_R\}$ be a set of linear forms
from $\Gamma^*$ such that
\begin{itemize}
\item each element of $\Delta$ is a {\em positive} multiple of an element of the sequence
$\Phi$;
\item for every $\beta_i$ in $\Phi$, there is a unique
  $\alpha_j\in\Delta$ which is a multiple of $\beta_i$.
\end{itemize}
Note that it is possible that $R<N$.

 For ${\bf j}$ a subset of $\{1,2,\ldots,N\}$, we denote by $C_{\bf
j}^+$  the closed convex cone in $\R^N$ generated by the set
$\{w_j|\, j\in {\bf j}\}$,  and by $C(\Phi_{\bf j})$ the closed
convex cone in $V^*$ generated by the set  $\{\beta_j|\,j\in {\bf
j}\}$.

For $\lambda\in \Gamma^*$, denote by $\iota_\Phi(\lambda)$ the number
of points with integral coordinates in $\Pi_\Phi(\lambda)$.  Thus
$\iota_\Phi(\lambda)$ is the number of solutions of the equation
$\sum_{k=1}^N x_k \beta_k=\lambda$ in {\em nonnegative integers}
$x_k$.  The function $\lambda\mapsto \iota_\Phi(\lambda)$ is called
the {\sl vector partition function} associated to $\Phi$. The number
$\iota_\Phi(\lambda)$ is zero if $\lambda$ does not belong to
$C(\Phi)\cap \Z\Phi$.

Recall the definition of the space of meromorphic functions
$M^{\Gamma\CA}$ defined in Section \ref{sec:rathpa} and the expansion
map $r^+$ defined in Section \ref{sec:expfus}.  Let
$$F_{\Phi}=\frac{1}{\prod_{i=1}^N(1-e^{\beta_i})}.$$ This function
is in the ring $M^{\Gamma\CA}$ and, almost by definition, the
expansion $r^+(F_{\Phi})$ is the generating function for
$\iota_\Phi$:
$$ r^+(F_{\Phi})=\sum_{\lambda\in \Gamma^*} \iota_\Phi(\lambda)
e^{\lambda}.$$

We can thus apply Theorem \ref{thm:main} and obtain a residue formula for
$\iota_{\Phi}(\lambda)$. We give this formula below in a slightly more precise form.

Similarly to the notation introduced earlier, we denote by $\bfi$
the set of linearly independent $n$-tuples of elements of the
sequence $\Phi$. For each such basis $\sigma\in\bfi$ of $V^*$, we
denote by  $C(\sigma)$ the cone generated by the elements of
$\sigma$ and by $G(\sigma,\Gamma)$ the lattice $2i\pi
(\Z\sigma)^*$ so that
$$G(\sigma,\Gamma)=\{p\in  V_\C|\,
e^{\langle\beta,p\rangle}=1,\text{ for all } \beta\in \sigma\}.$$
Clearly, $2\pi i\Gamma\subset G(\sigma,\Gamma)$, and we may choose
a finite, reduced set of elements $RG(\sigma,\Gamma)\subset
G(\sigma,\Gamma)$, which is in one-to one correspondence with the
finite factor group $G(\sigma,\Gamma)/2\pi i\Gamma$.

Given a chamber $\c$, we denote by $\CB(\Phi,\c)$ the set of
 $\sigma\in \CB(\Phi)$ such that $\c\subset C(\sigma)$ and
define
$$G(\Phi,\c,\Gamma)=\bigcup_{\sigma\in
  \CB(\Phi,\c)}G(\sigma,\Gamma)\;\text{ and }\;
RG(\Phi,\c,\Gamma)=\bigcup_{\sigma\in \CB(\Phi,\c)}RG(\sigma,\Gamma).$$

Introduce the convex polytope
$$\Box(\Phi)=\sum_{i=1}^N [0,1]\beta_i.$$

Now we are in position to formulate the appropriate version of our
Theorem~\ref{thm:main}.
\begin{theorem}
  \labelm{thm:vpartition} Denote by $\iota[\c,\Phi]$ the
  periodic-polynomial function on $\Gamma^*$ given by
$$ \sum_{p\in RG(\Phi,\c,\Gamma)}e^{-\langle\lambda,p\rangle}
\JK\c{\Tres\left(\frac{e^{\langle\lambda,z\rangle}}{\prod_{i=1}^N
(1-e^{\langle\beta_i,p\rangle}e^{-\langle\beta_i,z\rangle})}\right)}\volga.$$
Then, for any $\lambda\in (\c-\Box(\Phi))\cap \Gamma^*$, we have
\begin{equation}
  \label{eq:ii}
\iota_{\Phi}(\lambda)=\iota[\c,\Phi](\lambda).
\end{equation}
\end{theorem}
\begin{remark}
   We assumed that $\Phi$
linearly generates $V^*$, hence if $\c$ is a big chamber contained
in $C(\Phi)$, then the set $\c-\Box(\Phi)$ contains $\overline\c$.
This means that the formula  \eqref{eq:ii} is in particular true
for $\lambda\in \overline\c \cap \Gamma^*$. The set
$\cnull-\Box(\Phi)$ remains equal to $\cnull$ and does not touch
the boundary of $C(\Phi)$.
\end{remark}
\begin{proof}
The set $\c$ being open, the set $\c-\Box(\Phi)$ coincide with
$\c-\Box(\Phi)^0$. The sum appearing in the theorem is a
restricted version of the sum in~~\eqref{eq:defs} defining
$s[F_\Phi,\Gamma]$. Note that the set of poles
$\SP(F_\Phi,\Gamma)$ appearing in that definition specializes to
the set $\bigcup_{\sigma\in
  \CB(\Phi)}G(\sigma,\Gamma)$ in our case.

Thus in order to deduce the statement of the theorem from Theorem
\ref{thm:main}, we only need to check  that if $p\in V_\C$ is such
that $\jk\c{\Tres(e^{\langle\lambda,z\rangle}F_\Phi(p-z))}\volga$
does not vanish, then $p$ is necessarily in $G(\Phi,\c,\Gamma)$.
Indeed, by Lemma \ref{poles}, if
$\Tres(e^{\langle\lambda,z\rangle}F_\Phi(p-z))\neq0$, then the set
$\Delta(p)=\{\beta\in\Phi|\,\langle\beta,p\rangle\in2\pi i\Z\}$
  has to span $V^*$.   The function $z\mapsto
  e^{\langle\lambda,z\rangle}F_{\Phi}(p-z)$ is in the space
  $\widehat{R}_{\CA(\Delta(p))}$. Its total residue can be written as a sum
  of functions $\phi_{\sigma}\in S_{\CA(\Delta(p))}$ with
  $\sigma\in\mathcal{B}(\Delta(p))$.
Now if
$$\JK\c{\Tres(e^{\langle\lambda,z\rangle}F_\Phi(p-z))}\volga\neq0,$$
then there exists a basis $\sigma\in \mathcal{B}(\Delta(p))$ such
that $\c\subset C(\sigma)$. This implies that $p$ is in
$G(\Phi,\c,\Gamma)$.
\end{proof}

\begin{remark}
\label{khp}
 To compute the residue formula of Theorem \ref{thm:vpartition}
 for $\iota_{\Phi}(\lambda)$,
 a precise determination of the set
$RG(\Phi,\c,\Gamma)$ is not necessary. We can indeed sum over any
bigger set, the extra terms contributing  $0$ to the sum. For
example, we can sum over a set of representatives of the finite
group of $n_{\Phi}$-th roots of unity of the torus $V_\C/2i\pi
\Gamma$, where $n_\Phi$ is such that $n_\Phi \Gamma^*\subset
\Z\sigma$ for any $\sigma\in \CB(\Phi)$. In particular, if the
system $\Phi$ is unimodular, then our set of $p$'s reduces to a
single point $p=0$. In this case,  we obtain that the vector
partition function $\iota_{\Phi}(\lambda)$ is given by the
polynomial
$$\JK\c{\Tres\left(\frac{e^{\langle\lambda,z\rangle}}{\prod_{i=1}^N
(1-e^{-\langle\beta_i,z\rangle})}\right)}\volga.$$ on each sector
$(\c-\Box(\Phi))\cap \Gamma^*$.
\end{remark}

\bigskip
 Let us comment on the novel aspects of Theorem
\ref{thm:vpartition}. Given $\lambda\in  C(\Phi)\cap \Gamma^*$,
the function $k\mapsto \iota_{\Phi}(k\lambda)$, $k=0,1,2,\ldots$,
counts the number of integral points in the dilated polytope
$k\Pi_{\Phi}(\lambda)$ of the rational polytope
$\Pi_{\Phi}(\lambda)$.  Clearly, the ray $\{k\lambda\}$  remains
in the closure of a chamber of $C(\Phi)$, and $k\mapsto
\iota_{\Phi}(k\lambda)$ is a periodic polynomial function of $k$
called the Ehrhart periodic-polynomial \cite{eh1} of the rational
polytope $\Pi_{\Phi}(\lambda)$. When $V$ is one-dimensional, this
case corresponds to enumeration of lattice points in rational
simplices and is the cornerstone of Ehrhart's work (see
\cite{eh3}, and references there).
 The vector partition function in this case is called the
 restricted partition function.
 Our formula of Theorem \ref{thm:vpartition} for the restricted
 partition function
clearly  coincides with results summarized in Comtet (\cite{co},
Th{\'e}or{\`e}me B. page 122), since we use  the same method of generating
functions and partial fraction decompositions, in a multivariate
setting.

 Theorem \ref{thm:vpartition} gives an
explicit residue formula for the number of integral points in the
polytope $\Pi_{\Phi}(\lambda)$, when $\lambda$ now varies in the
cone $C(\Phi)$. If $\c$ is a big chamber contained in $C(\Phi)$,
this formula is periodic-polynomial on the ''neighborhood''
$\c-\Box(\Phi)$ of $\overline \c$. This is somewhat surprising, as
the combinatorial nature of the polytope $\Pi_{\Phi}(\lambda)$
changes, when crossing walls of the big chambers. Thus these
different periodic-polynomial functions for the vector partition
function on different sectors coincide for neighboring chambers in
a strip containing their common boundary. Precisely, for two
chambers $\c_1$ and $\c_2$, the periodic-polynomial functions
$\iota[\Phi,\c_1]$ and $\iota[\Phi,\c_2]$ are equal on the set
$\Gamma^*\cap(\c_1-\Box(\Phi) )\cap (\c_2-\Box(\Phi))$.  This
implies some divisibility properties of the function
$\iota[\Phi,\c_1]-\iota[\Phi,\c_2]$ on adjacent chambers. We give
some illustrative examples for these properties of
$\iota[\Phi,\c]$ in the Appendix.

The relation between the number of integral points and the volume
of the polytope $\Pi_{\Phi}(\lambda)$  has been the subject of
several investigations
(see e.g. \cite{cs1,cs2,bv2,ms,g,dr}), starting with the fascinating results of
Khovanskii-Pukhlikov \cite{kp}.

 Recall that in Baldoni-Vergne
(\cite{bav}, Theorem 9), we discussed the Jeffrey-Kirwan residue
formula for the volume of the polytope $\Pi_\Phi(a)$. Let $\c$ be
a big chamber contained in $C(\Phi)$. Denote by $v[\Phi,\c,\vol]$
the polynomial function

$$v[\Phi,\c,\vol](a)
=\JK\c{\Tres\left(\frac{e^{\langle a,z\rangle}}{\prod_{i=1}^N
\langle\beta_i,z\rangle}\right)}\vol.$$ The volume of the polytope
$\Pi_{\Phi}(a)$ is given by a locally polynomial formula in $a$.
Explicitly, for $a$ varying in the closure of the big chamber
$\c$:
\begin{equation}
  \label{eq:jk}\mathrm{volume}(\Pi_\Phi)(a)=v[\Phi,\c,\vol](a).
  \end{equation}

The residue formula of Theorem \ref{thm:vpartition} for
$\iota_{\Phi}(\lambda)$ on the closure of the  chamber $\c$ may be
immediately deduced from the results of Brion-Vergne \cite{bv2} or
Cappell-Shaneson \cite{cs2} by applying Todd operators to the
volume function given by the residue formula above (\ref{eq:jk}).
It is satisfying, however, to obtain  ``explicit'' and very
similar formulae for volumes of polytopes and for the number of
integral points in polytopes by residue methods, in a parallel
way. It is puzzling to see that the formula for the number of
points holds in a larger set than we would guess from its
continuous analogue, the volume function.

\subsection{Minkowski sum of rational convex  polytopes and families of
partition polytopes  }
\labelm{sec:mink}

In this section, we briefly describe how to realize any rational
convex polytope as a partition polytope $\Pi_{\Phi}(a)$.

Recall some standard conventions.  Faces of a polytope $\Pi$ of
dimension $r$ may have any codimension from $0$ to $r$. A face of
codimension $1$ is called a {\em facet}. A face of dimension $0$
is a {\em vertex}, a face of dimension $1$ is an {\em edge}. The
polytope $\Pi$ is said to be simple if each vertex of $\Pi$ is the
source of exactly $r$ edges. Given  a rational polytope $\Pi$ in a
vector space endowed with a lattice $\Theta$ of full rank, a face
$f$ of $\Pi$ is called {\em reticular} if the affine space spanned
by $f$ contains a point of $\Theta$. In particular, a vertex is
reticular if and only if it belongs to  $\Theta$. A rational
polytope $\Pi$ is {\em integral} if all its vertices belong to the
lattice $\Theta$.

Let $\Phi$ be again a sequence of $N$ linear forms
$[\beta_1,\ldots,\beta_N]$ generating $V^*$ and lying on the same
side of a  hyperplane. For $a$ in the interior of $C(\Phi)$, the
polytope $\Pi_{\Phi}(a)$ has dimension $N-n$.

We keep our earlier notations. For a basis $\sigma\in\CB(\Phi)$ of
$V^*$, we denote by $v_{\sigma}$ the map from $V^*$ to $\R^N$
defined by $v_{\sigma}(\beta_j)=w_j$ for all $\beta_j\in \sigma$.
Clearly, $\rho \circ v_{\sigma}$ is the identity on $V^*$. If
$\beta_k$ is not in $\sigma$, the vector $w_k-v_\sigma(\beta_k)$
is in the subspace $\rho^{-1}(0)$.

Recall
\begin{proposition}[\cite{bv2}]
\label{pro:vertex}
 Let $\c$ be a big chamber contained in $C(\Phi)$.
\begin{itemize}
\item For any $a\in\c$, the convex polytope $\Pi_{\Phi}(a)$ is simple,
  with vertices $v_{\sigma}(a)$, $\sigma\in\CB(\Phi,\c)$.  These
  vertices are all distinct, and the   $(N-r)$ edges  of $\Pi_{\Phi}(a)$
  with source at the vertex $v_\sigma(a)$ are the vectors
  $w_k-v_\sigma(\beta_k)$, where $\beta_k\notin \sigma$.
\item If $a\in\overline{\c}$, then the vertices of the convex polytope
  $\Pi_{\Phi}(a)$ are the points $v_{\sigma}(a)$,
  $\sigma\in\CB(\Phi,\c)$. Some of these points may coincide.
\item The faces of dimension $j$ of $\Pi_{\Phi}(a)$ are the sets
  $\rho^{-1}(a)\cap C_{\bf j}^+$, where ${\bf j}$ is a subset of
  $\{1,2,\ldots,N\}$ of cardinality $(n+j)$, and  such that $\c\subset
  C(\Phi_{\bfj})$. Here $\Phi_{\bfj}$ stands for the set of forms
  $\beta_i$, with indices $i\in\bfj$.
\end{itemize}
\end{proposition}

Let $\lambda\in C(\Phi)\cap \Z\Phi$. Consider the function
$k\mapsto \iota_{\Phi}(k\lambda)$, where $k$ is a non-negative
integer. Now we will see that our formula (\ref{eq:ii}) for  the Ehrhart
periodic-polynomial $E[\lambda](k)=\iota_{\Phi}(k\lambda)$ is
actually polynomial in $k$ if all the vertices of
$\Pi_{\Phi}(\lambda)$ have integral coordinates. More generally,
we will show that our formula  is compatible with
some of the results of \cite{eh1,st1,mm2} on the
periodic-polynomial behavior of $E[\lambda](k)$.
\begin{lemma}\label{stanley}
  If $M$ is an integer such that $M\Pi_{\Phi}(\lambda)$ is integral,
  then $E[\lambda](k)=\sum_{\zeta^M=1}\zeta^{k}P_\zeta(k)$, where
  $\zeta$ varies over $M$th roots of unity and $P_\zeta$ is a
  polynomial.  Furthermore, if each $j$-face of $\Pi_{\Phi}(\lambda)$
  is reticular, then, with the exception of $\zeta= 1$, each polynomial $P_\zeta$ is
  of degree strictly less than  $j$.  {\rm(}The degree of $0$ is set to be   $-1$
  {\rm )}.
\end{lemma}
\begin{proof}
After Theorem  \ref{thm:vpartition}, we have $$E[\lambda](k)=
\sum_{p\in RG(\Phi,\c,\Gamma)}e^{-k\langle
    \lambda,p\rangle} P_{(p)}[\lambda](k),$$
  where
  $$P_{(p)}[\lambda](k)=\JK\c{\Tres\left(\frac{e^{k\langle\lambda,z\rangle}}
  {\prod_{i=1}^N
        (1-e^{\langle\beta_i,p\rangle}e^{-\langle\beta_i,z\rangle})}\right)}\volga.$$

  We first show that all exponentials $e^{-\langle \lambda,p\rangle}$
  are $M$-th roots of unity. For each $p\in RG(\Phi,\c,\Gamma)$, there
  exists $\sigma\in \CB(\Phi,\c)$ such that $p$ is a solution of the
  equations $e^{\langle\beta_i,p\rangle}=1$ for all $\beta_i\in
  \sigma$. Since $\c\subset C(\sigma)$, we can write
  $\lambda=\sum_{\beta_i\in {\sigma}} x_i \beta_i$, where each $x_i$
  is a rational nonnegative number.  The point
  $v_\sigma(\lambda)=\sum_{\beta_i\in \sigma} x_i w_i$ is a vertex of
  the polytope $\Pi_{\Phi}(\lambda)$. If $M v_\sigma(\lambda)$ is
  integral, then all numbers $Mx_i$ are integers, so we have
  $e^{M\langle\lambda,p\rangle}=1$ as $p$ is a solution of the
  equations $e^{\langle\beta_i,p\rangle}=1$ for all $\beta_i\in
  \sigma$.

For each $p\in RG(\Phi,\c,\Gamma)$, we denote
$${\bf i}(p) =\{i\in\{1,2,\ldots,N\}|\;e^{\langle\beta_i,p\rangle}=1\}.$$
The set ${\bf i}(p)$ is of cardinality $n+r$ with $r\geq 0$ and
$\c\subset C(\Phi_{{\bf i}(p)})$. The function $z\mapsto
\prod_{i=1}^N
(1-e^{\langle\beta_i,p\rangle}e^{-\langle\beta_i,z\rangle})$ is
divisible by $\prod_{i\in {\bf i}(p)}\beta_i$, so that in the
computation of the total residue, we need to take the expansion of
$e^{\langle\lambda,z\rangle}$ only up to order $r$. The polynomial
$P_{(p)}[\lambda](k)$ is thus of degree less or equal to $r$. To
prove the second assertion of the lemma, we need to show that for
$r\geq j$, the complex number $e^{\langle\lambda,p\rangle}$ in
front of this polynomial is equal to $1$. The $r$-face
$f=\Pi_{\Phi}(\lambda)\cap C_{{\bf i}(p)}^+$ is reticular as
$r\geq j$. Thus we can write $\lambda=\sum_{i\in {\bf i}(p)}l_i
\beta_i$ where $l_i$ are in $\Z$. Again, we see that we have
$e^{\langle\lambda,p\rangle}=1$ as $p$ is a solution of the
equations $e^{\langle\beta_i,p\rangle}=1$, for all $i\in {\bf
i}(p)$.
\end{proof}

The notion of big chambers in $C(\Phi)$ is closely related to the
Minkowski sum of the corresponding partition polytopes as follows.

\begin{lemma}
  Let $a,b\in C(\Phi)$. The Minkowski sum $\Pi_{\Phi}(a)+\Pi_{\Phi}(b)$ of
  the polytopes $\Pi_{\Phi}(a)$ and $\Pi_{\Phi}(b)$ is equal to the
  polytope $\Pi_{\Phi}(a+b)$ if and only if there exists a big chamber
  $\c$ contained in $C(\phi)$ such that $a,b \in \overline{\c}$.
\end{lemma}
\begin{proof}
  Clearly the polytope $\Pi_{\Phi}(a)+\Pi_{\Phi}(b)$ is a subset of
  the polytope $\Pi_{\Phi}(a+b)$.

  Let $\c$ be a chamber contained in $C(\Phi)$ such that $a,b\in \overline{\c}$. Hence
  $a+b$ is in $\overline{\c}$. Let us prove that $\Pi_{\Phi}(a+b)$ is equal
  to $\Pi_{\Phi}(a)+\Pi_{\Phi}(b)$.  By the description of the
  vertices given in Proposition \ref{pro:vertex}, any element $\bx$ of
  the polytope $\Pi_{\Phi}(a+b)$ can be written as $\sum_{\sigma\in
    \CB(\Phi,\c)} t_\sigma v_\sigma (a+b)$, with $\sum t_\sigma=1$.
  Then we may write $\bx=\bx_1+\bx_2$, with $\bx_1=\sum_{\sigma\in
    \CB(\Phi,\c)} t_\sigma v_\sigma (a)$ and $\bx_2=\sum_{\sigma\in
    \CB(\Phi,\c)} t_\sigma v_\sigma (b)$, with $\bx_1\in \Pi_{\Phi}(a)$
  and $\bx_2\in \Pi_{\Phi}(b)$.

  Conversely, let $a,b\in C(\Phi)$ such that
  $\Pi_{\Phi}(a)+\Pi_{\Phi}(b)=\Pi_{\Phi}(a+b)$.  Consider then a
  chamber $\c$ contained in $C(\Phi)$ such that $a+b\in \overline{\c}$.  Let $\sigma$ such
  that $\c\subset C(\sigma)$.  The point $v_{\sigma}(a+b)$, being in
  $\Pi_{\Phi}(a+b)$, can be written as $\bx_1+\bx_2$ with $\bx_1\in
  \Pi_{\Phi}(a)$ and $\bx_2\in \Pi_{\Phi}(b)$. Since those coordinates
  of $v_\sigma(a+b)$ corresponding to  $\beta_k\notin \sigma$ are
  equal to $0$,
  we see that the $k$th coordinate of $\bx_1, \bx_2$ vanish when
  $\beta_k\notin\sigma$. This implies that $\bx_1=v_\sigma(a)$
  and $\bx_2=v_\sigma(b)$. Thus $a,b\in \cap_{\sigma\in
    \CB(\Phi,\c)}C(\sigma)=\overline{\c}$. The Lemma is proved.
   \end{proof}

   When $(\lambda_1,\lambda_2,\ldots,\lambda_s)$ are elements of
   $\overline\c\cap \Gamma^*$, and $k_i$ are nonnegative integers, the
   polytope $\Pi_{\Phi}(k_1 \lambda_1+k_2 \lambda_2+\cdots +k_s
   \lambda_s)$ is a rational convex polytope which is the weighted
   Minkowski sum $k_1\Pi_{\Phi}(\lambda_1)+\cdots +
   k_s\Pi_{\Phi}(\lambda_s)$. As $(k_1\lambda_1+\cdots +k_s \lambda_s)$
   varies in $\overline{\c}$, the function
   $\iota_{\Phi}(k_1\lambda_1+\cdots +k_s \lambda_s)$ is a
   periodic-polynomial function of $k_i$. This extension of Ehrhart's
   result is well-known \cite{mm1}. As in Lemma \ref{stanley}, if the
   polytopes $\Pi_{\Phi}(\lambda_k)$ have integral vertices, then the
   function $(k_1,k_2,\ldots,k_s)\mapsto
   \iota_{\Phi}(k_1\lambda_1+k_2\lambda_2+\cdots +k_s\lambda_s )$ is a
   polynomial function of $k_1,k_2,\ldots,k_s$.

   Now recall briefly (cf.\cite{bv2}) how any convex polytope $\Pi$
   can be embedded in a family $\Pi_{\Phi}(a)$ of partition polytopes.

Let $E$ be a real vector space of dimension $r$. Let $\Pi\subset E$ be
a convex polytope. We can always choose $N$ vectors $u_k\in E^*$ and a
sequence of real numbers $\bh=(h_1,h_2,\ldots,h_N)\in \R^N$ such that
$\Pi=\Pi(\bh)$, where
$$\Pi(\bh)=\{v\in E~\vert~\langle u_k,v\rangle+h_k\geq 0,~1\leq k\leq
N\}~.$$

As $\Pi$ is compact, the vectors $u_k$ generate $E^*$.  We do not
necessarily assume here that this set of inequalities is minimal.
Consider the map $U:\R^N\to E^*$ defined by
 $$(x_1,x_2,\ldots,x_N)\mapsto x_1u_1+x_2u_2+\cdots +x_N u_N,$$
 and let $V$ be the $n=(N-r)$-dimensional vector space
 $V=U^{-1}(0)$. The restrictions $\beta_i$ of the linear coordinates
 $x_i$ to the vector space $V$ form a system $\Phi$ of elements of
 $V^*$.  The elements $\beta_i$ of the system $\Phi$ satisfy
 the equation $\langle u_1,v\rangle\beta_1+\cdots +\langle
 u_N,v\rangle\beta_N=0$ for all $v\in E$.

\begin{lemma}
The polytope $\Pi(\bh)$ is isomorphic to the partition polytope
$\Pi_{\Phi}(h_1\beta_1+\cdots +h_N\beta_N)$.

\end{lemma}
\begin{proof}
  A point of the polytope $\Pi_{\Phi}(h_1\beta_1+\cdots +h_N\beta_N)$
  is a point $(l_1,l_2,\ldots,l_N)\in \R_+^N$ such that
  $l_1\beta_1+l_2\beta_2+ \cdots +l_N\beta_N=h_1\beta_1+\cdots
  +h_N\beta_N$. This implies that there exists a unique $v\in E$ such that
  $l_i-h_i=\langle u_i,v\rangle$, so that $\langle
  u_k,v\rangle+h_k=l_k\geq 0$ and $v$ is in $\Pi(\bh)$.
\end{proof}

Assume now that $E$ is endowed with a lattice $\Theta$ and that the
polytope $\Pi$ is rational. Then there exist vectors $u_k\in \Theta^*$
and integers $h_k$ such that
$$\Pi(\bh)=\{v\in E~\vert~\langle u_k,v\rangle+h_k\geq 0,~1\leq k\leq N\}~.$$

We can always assume, adding superfluous elements $u_k$ to $\Theta^*$
if necessary, that $\langle u_k,v\rangle\in \Z$ if and only if $v\in
\Theta^*$. Then the set of integral points in the polytope
$\Pi_{\Phi}(h_1\beta_1+\cdots +h_N\beta_N)$ is in bijection with the
set of integral points in $\Pi(\bh)$.

More generally, let $\Pi_1, \Pi_2,\ldots, \Pi_k$ be a set of
rational convex polytopes in $E$. The Minkowski sum $t_1
\Pi_1+t_2\Pi_2+\cdots +t_s\Pi_s$, where each $t_k$ is a nonnegative
real number, can be described as a set $\{v\in E|~ \langle
u_k,v\rangle+t_1h_k^1+t_2h_k^2+\cdots+t_s h_k^s\geq 0\}$. As
before, we consider the map $U:\R^N\to E^*$ defined by
$$(x_1,x_2,\ldots,x_N)\mapsto x_1u_1+x_2u_2+\cdots +x_N u_N$$
Let $V=U^{-1}(0)$ and $\Phi$ the system of linear forms obtained
by the restrictions of the linear coordinates. Then the points
$\lambda_i=\sum h_k^i\beta_k$ belong to the closure of a chamber
contained in $C(\Phi)$, and the family $t_1 \Pi_1+t_2\Pi_2+\cdots
+t_s\Pi_s$ is a member of the family of partitions polytopes
$\Pi_\Phi(a)$, where $a=t_1\lambda_1+\cdots +t_s\lambda_s$ varies
in the closure of a chamber contained in $C(\Phi)$ .

Thus the results of this article give, in particular, ``explicit
periodic-polynomial formulae'' for mixed enumerators in functions of
the inequalities defining the family of Minkowski polytopes
$t_1\Pi_1+\cdots +t_s \Pi_s$.

\subsection{Sums of exponentials over partition polytopes}

Consider now a point $\by=(y_1,y_2,\ldots,y_N)$ in $\C^N$ and the
exponential  function $e^{\langle \by,\bx\rangle}=e^{\sum_{i=1}^N
x_i y_i}$ over $\R^N$. We consider the function
$$\Sm[e^\by,\Phi](\lambda)=\sum_{\xi\in \Pi_{\Phi}(\lambda)\cap
\Z^N}e^{\langle \by,\xi\rangle}.$$

Let
$$F_{\Phi,\by}(z)=\frac{1}{\prod_{j=1}^N(1-e^{y_j}e^{\langle\beta_j,z\rangle})}.$$

Almost by definition, the expansion $r^+(F_{\Phi,\by})$ is the
generating function for $\Sm[e^\by,\Phi]$:

$$ r^+(F_{\Phi,\by})=\sum_{\lambda\in \Gamma^*}
\Sm[e^\by,\Phi](\lambda)e^{\lambda}.$$

Let $\sigma\in\CB(\Phi)$. We introduce the set
$$G(\sigma,\by, \Gamma)=\{p\in  V_\C|\, e^{\langle\beta_j,p\rangle}=e^{-y_j}
\text{ for all } \beta_j\in \sigma\}.$$

Clearly, if $\gamma\in 2\pi i\Gamma$ and $p\in
G(\sigma,\by,\Gamma)$, then $p+\gamma\in  G(\sigma,\by, \Gamma)$
and we may choose a finite, reduced set of elements
$RG(\sigma,\by,\Gamma)\subset G(\sigma,\by,\Gamma)$, which is in
one-to-one correspondence with the finite coset
$G(\sigma,\by,\Gamma)/2\pi i\Gamma$.

For a chamber $\c$, we define $RG(\Phi,\by,\c,\Gamma)$ to be the
union of the sets $RG(\sigma,\by,\Gamma)$ over all bases
$\sigma\in\CB(\Phi)$ such that $\c\subset C(\sigma)$.

Applying our Theorem \ref{thm:main}, we obtain:
\begin{theorem}
  Let $\c$ be a big chamber of a sequence $\Phi=[\beta_1,\ldots,\beta_N]$,
  and let $\by\in \C^N$. Denote by $\iota[\c,\by,\Phi]$ the
  exponential-polynomial function on $\Gamma^*$ equal to
$$\sum_{p\in RG(\Phi,\by,\c,\Gamma)}
 e^{-\langle\lambda,p\rangle}\JK\c
{\Tres\left(\frac{e^{\langle\lambda,z\rangle}}{\prod_{i=1}^N
(1-e^{\langle\beta_i,p\rangle}e^{y_i}e^{-\langle\beta_i,z\rangle})}\right)}\volga.$$
 Then, for any $\lambda\in (\c-\Box(\Phi))\cap \Gamma^*$,
 the function $\lambda\mapsto \Sm[e^\by,\Phi](\lambda)$ is given
by the exponential-polynomial formula
$$\Sm[e^\by,\Phi](\lambda)=\iota[\c,\by,\Phi](\lambda).$$
\end{theorem}

Let us compare this expression to the ``explicit'' formula of
\cite{b,ba} for sums of exponentials over the integral
points of a convex polytope for sufficiently generic $\by$.

Let $\sigma\in \CB(\Phi,\c)$ and assume that $\by$ is sufficiently
generic. Then
 for every  $p\in G(\sigma,\by, \Gamma)$,
   we have $e^{y_j}e^{\langle\beta_j,p\rangle}=1 $ for all $\beta_j\in
 \sigma$, while
$e^{y_k}e^{\langle\beta_k,p\rangle}\neq 1 $ for all $\beta_k\notin
\sigma$. Thus, for $p\in G(\sigma,\by, \Gamma)$, the function
 $z\mapsto F_{\Phi,\by}(p-z)$ is equal to
 $$\frac{1}{\prod_{\beta_j\in
 \sigma}(1-e^{-\langle\beta_j,z\rangle})}\prod_{\beta_k\notin
 \sigma}\frac{1}{(1-e^{y_k}e^{\langle\beta_k,p\rangle} e^{-\langle\beta_k,z\rangle})},$$
and we obtain by Lemma \ref{poles}
$$\Tres\left( e^{\langle\lambda,z-p\rangle}F_{\Phi,\by}(p-z)\right)
=e^{-\langle\lambda,p\rangle}\prod_{k\notin
 \sigma}\frac{1}{(1-e^{y_k}e^{\langle\beta_k,p\rangle})}\times
 \frac{1}{\prod_{\beta_j\in \sigma} \beta_j}.$$

For $\by$ generic, all the subsets $G(\sigma,\by, \Gamma)$ are
disjoint as $\sigma$ varies in $\CB(\Phi)$.
 Thus for generic $\by$ we obtain a formula
$\Sm[e^\by,\Phi](\lambda)$ as a linear combination of the pure
exponential functions  $\lambda \mapsto
e^{-\langle\lambda,p\rangle}$ associated to the elements $p\in
RG(\Phi,\by,\c,\Gamma)$.
\begin{theorem}
\label{exp}
 Let $\c$ be a chamber and $\by$ a generic element of
$\C^N$. Let $E[\c,\by,\Phi](\lambda)$ be the function of $\lambda$
defined by
$$\sum_{\sigma\in \CB(\Phi,\c)}\frac{1}{\vol_{\Gamma^*}(\sigma)}
 \sum_{p\in RG(\sigma,\by,
\Gamma)}e^{-\langle\lambda,p\rangle}\prod_{\beta_k\notin
 \sigma}\frac{1}{(1-e^{y_k}e^{\langle\beta_k,p\rangle})}.$$
Then, for $\lambda\in (\c-\Box(\Phi))\cap \Gamma^*$, we have the
''explicit'' formula
$$\Sm[e^\by,\Phi](\lambda)=E[\c,\by,\Phi](\lambda).$$
\end{theorem}

\begin{remark}
On the set $\overline{\c}\cap \Gamma^*$, it is possible to deduce
this formula from the Baum-Fulton-MacPherson  equivariant
Riemann-Roch formula applied to the (possibly singular) toric
variety and its holomorphic line bundle associated with the
polytope $\Pi_{\Phi}(\lambda)$, at least when this polytope is
integral. This dictionary between toric varieties and rational
polytopes is used in several proofs of formulae for sums of
functions over integral points of convex integral polytopes~
\cite{cs1,ms,g}.
\end{remark}

Let us rewrite the formula of Theorem \ref{exp} in  geometric
terms in the case when $\Phi$ is a unimodular system and $\lambda$
is in an open chamber $\c$ of $C(\Phi)$. First we note that for any
$\sigma\in \CB(\Phi,\c)$ each set $RG(\sigma,\by, \Gamma)$
consists of just one element and the number
$\vol_{\Gamma^*}(\sigma)$ is equal to $1$. Thus the formula for
$E[\c,\by,\Phi]$ is simply indexed by the set $\CB(\Phi,\c)$,
which also indexes the vertices of the polytope
$\Pi_{\Phi}(\lambda)$. Let $\sigma\in \CB(\Phi,\c)$ and $p$ be an
element such that $e^{\langle\beta_i,p\rangle}=e^{-y_i}$, for all
$\beta_i\in \sigma$ . If $\lambda=\sum_{\beta_i\in \sigma}
x_i\beta_i$, then
$e^{-\langle\lambda,p\rangle}=e^{\sum_{\beta_i\in \sigma} x_i
  y_i}=e^{\langle \by,v_\sigma(\lambda)\rangle}$ is the value of the
exponential function $e^\by$ at the vertex $v_\sigma(\lambda)$ of
the polytope $\Pi_{\Phi}(\lambda)$. Similarly, the  edges
$a_k^{\sigma}=w_k-s_\sigma(\beta_k)$ passing through
$v_\sigma(\lambda)$ are such that $e^{\langle
\by,a_k^{\sigma}\rangle}= e^{y_k}e^{\langle\beta_k,p\rangle}$.
Thus, for the simple polytope $\Pi_{\Phi}(\lambda)$ associated to
an unimodular system $\Phi$,  we obtain
$$\sum_{\xi\in \Pi_{\Phi}(\lambda)\cap Z^N}e^{\langle \by,\xi\rangle}=
\sum_{v}\frac {e^{\langle
\by,v\rangle}}{\prod_{a_j(v)}(1-e^{\langle
    a_j(v),\by\rangle})},$$
where $v$ varies over the vertices of the polytope
$\Pi_{\Phi}(\lambda)$ and   $a_j(v)$ varies over the primitive
edges of the polytope with source at the vertex $v$. One may
recognize here the localization formula for the equivariant index
applied to the smooth toric variety and its holomorphic line
bundle associated with the polytope $\Pi_{\Phi}(\lambda)$.

In the general  case, Theorem \ref{exp} implies Formula 3.4.1 of
Brion-Vergne (\cite{bv2}). Again, our results here imply that this
formula holds on a larger set of $\lambda$s, on which the elements
$v_\sigma(\lambda)$ are not necessarily vertices of the polytope
$\Pi_{\Phi}(\lambda)$.

\subsection{Summing the values of an exponential-polynomial function
over
partition polytopes}

For  an exponential-polynomial function $f$ on $\R^N$, we denote
$\Sm[f,\Phi](\lambda)=\sum_{\xi\in \Pi_{\Phi}(\lambda)\cap\Z^N}f(\xi)$.

 Recall the
definition of the polynomial functions
$$c(x,h)=\frac{1}{(h-1)!}(x+1)(x+2)\cdots (x+(h-1)), $$
(where $c(x,1)=1$)  which form a basis of polynomial functions on
$\R$ as $h$ runs through the positive integers. Let again
$\Gamma\subset V$ be a lattice in an $n$-dimensional vector space,
and $\Phi$ be a sequence $[\beta_1,\ldots,\beta_N]$ of linear forms
from $\Gamma^*$ lying on the same side of an hyperplane and
generating $V^*$. Also, fix $\by=(y_1,y_2,\ldots,y_N)\in \C^N$ and
let ${\bf h}=(h_1,h_2,\ldots,h_N)$ be a list of positive integers.
Consider the exponential-polynomial function on $\R^N$ given by
$$f_{{\bf h},\by}(\bx)=e^{\langle \by,\bx\rangle}\prod_{i=1}^N c(x_i,h_i).$$

The generating function for the function $\Sm[f_{{\bf
h},\by},\Phi]=\sum_{\xi\in \Pi_{\Phi}(\lambda)}f_{{\bf
h},\by}(\xi)$ is the function
$$F_{\Phi,\by,{\bf h}}(z)= \frac{1}{\prod_{i=1}^N
(1-e^{y_i}e^{\langle\beta_i,z\rangle})^{h_i}}.$$

Note that the set $\Box(F_{\Phi,\by,{\bf h}})$ is the set
$\Box(\Phi,{\bf h})= \sum_{i=1}^N [0,1] h_i\beta_i$. It
always contains $\Box(\Phi)$.

Theorem \ref{thm:main} states:
\begin{theorem}\labelm{thm:final}
Let $\c$ be a chamber, $\by\in \C^N$ and ${\bf h}\in \N^N$. Let
$\Smm[\c,\by,{\bf h},\Phi]$ be the exponential-polynomial function
on $\Gamma^*$ equal to
$$\sum_{p\in RG(\Phi,\by,\c,\Gamma)}
 e^{-\langle\lambda,p\rangle}\JK\c
{\Tres\left(\frac{e^{\langle\lambda,z\rangle}}{\prod_{i=1}^N
(1-e^{\langle\beta_i,p\rangle}e^{y_i}e^{-\langle\beta_i,z\rangle})^{h_i}}\right)}
\volga.$$

Then, for any $\lambda\in (\c-\Box(\Phi,{\bf h}))\cap \Gamma^*$,
the function $\lambda\mapsto \Sm[f_{{\bf h},\by},\Phi](\lambda)$
is given by the exponential-polynomial formula
$$\Sm[f_{{\bf h},\by},\Phi]=\Smm[\c,\by,{\bf h}, \Phi](\lambda).$$

In particular, if $\c$ is a chamber contained in $C(\Phi)$, then
for any exponential-polynomial function $f\in \EP(\R^N)$, the
function $\lambda\mapsto \Sm[f,\Phi](\lambda)$ is given by an
exponential-polynomial function $\Smm[\c,f,\Phi]$ for $\lambda\in
(\c-\Box(\Phi))\cap \Gamma^*$. The set  $(\c-\Box(\Phi))$ contains
$\overline \c$.

\end{theorem}

\section{Appendix: Examples}

Let $V$ be a $2$-dimensional real vector space with basis $(e_1,e_2)$; then
the dual vector space $V^*$ has basis $e^*_1,e^*_2$. Sometimes we denote a vector
 $a_1e^*_1+a_2e^*_2$ in $V^*$  simply by $(a_1,a_2)$; similarly
 $(z_1,z_2)$ stands for $z_1 e_1+z_2 e_2$ in $V$. 
We take the lattice $\Gamma$ to be $\Z e_1\oplus \Z e_2$.

\subsection{The arrangement $A_2$}
Let $$\Phi=\{e^*_1, e^*_2, e^*_1+e^*_2 \}.$$
 The space $R_{\CA(\Phi)}$ consists of rational functions
 $f(z_1,z_2)$  on $V_\C$ with denominator a product of powers of the linear forms
 $z_1$, $z_2$, $z_1+z_2$. The system $\Phi$ is unimodular.

 The closed cone $C(\Phi)$ generated by $\Phi$ is the first quadrant
 $\{a_1\geq 0, a_2\geq 0\}$. There are three big chambers for the
 system $\Phi$: the exterior of the cone $C(\Phi)$ denoted by
 $\cnull$, the chamber $\c^1=\{a_2>0,a_1>a_2\}$ and  the chamber
 $\c^2=\{a_1>0, a_2>a_1\}$.

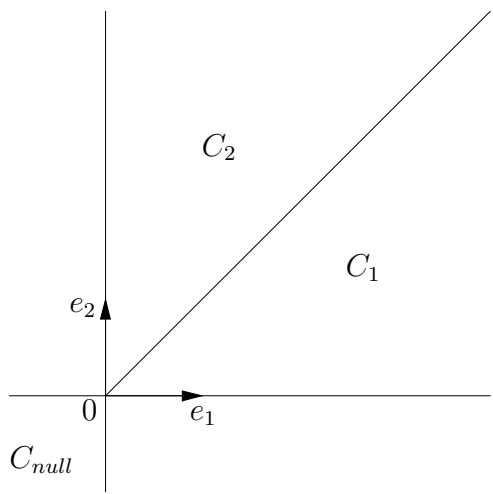
\begin{figure}[ht]
\begin{center}
\setlength{\unitlength}{0.00083333in}
\begingroup\makeatletter\ifx\SetFigFont\undefined%
\gdef\SetFigFont#1#2#3#4#5{%
  \reset@font\fontsize{#1}{#2pt}%
  \fontfamily{#3}\fontseries{#4}\fontshape{#5}%
  \selectfont}%
\fi\endgroup%
{\renewcommand{\dashlinestretch}{30}
\begin{picture}(3024,3039)(0,-10)
\path(12,612)(3012,612)
\path(612,612)(1212,612)
\blacken\path(1092.000,582.000)(1212.000,612.000)(1092.000,642.000)(1092.000,582.000)
\path(612,612)(612,1212)
\blacken\path(642.000,1092.000)(612.000,1212.000)(582.000,1092.000)(642.000,1092.000)
\path(612,612)(3012,3012)
\path(612,12)(612,3012)
\put(2112,1362){\makebox(0,0)[lb]{\smash{{{\SetFigFont{12}{14.4}{\rmdefault}{\mddefault}{\updefault}$C_1$}}}}}
\put(1212,2112){\makebox(0,0)[lb]{\smash{{{\SetFigFont{12}{14.4}{\rmdefault}{\mddefault}{\updefault}$C_2$}}}}}
\put(12,162){\makebox(0,0)[lb]{\smash{{{\SetFigFont{12}{14.4}{\rmdefault}{\mddefault}{\updefault}$C_{null}$}}}}}
\put(387,1137){\makebox(0,0)[lb]{\smash{{{\SetFigFont{12}{14.4}{\rmdefault}{\mddefault}{\updefault}$e_2$}}}}}
\put(1137,462){\makebox(0,0)[lb]{\smash{{{\SetFigFont{12}{14.4}{\rmdefault}{\mddefault}{\updefault}$e_1$}}}}}
\put(462,462){\makebox(0,0)[lb]{\smash{{{\SetFigFont{12}{14.4}{\rmdefault}{\mddefault}{\updefault}$0$}}}}}
\end{picture}
}
\caption{The chambers}
\end{center}
\end{figure}

The linear forms $J(\c_1,da)$ and $J(\c_2,da)$ are easily
computed. For a rational function $f(z_1, z_2)$ in the space
$R_{\CA(\Phi)}$, we have
$$\jk{\c_1}{\Tres f}{\volga}=\res_{z_2=0}\res_{z_1=0}(f(z_1,z_2)\;dz_1\,dz_2),$$ while
$$\jk{\c_2}{\Tres f}{\volga}=\res_{z_1=0}\res_{z_2=0}(f(z_1,z_2)\;dz_1\,dz_2).$$

We denote by $\Phi_n$ the system of $3n$ vectors where each linear
form $e^*_1$, $e^*_2$, $e^*_1+e^*_2$ has  multiplicity $n$.

The polytope $\Box(\Phi)$ is the convex hull of the six points
$0$, $e^*_1$, $e^*_2$, $2e^*_1+e^*_2$, $2e^*_2+e^*_1$, $2e^*_1+2e^*_2$.
The polytope $\Box(\Phi_n)$ is the dilated convex
polytope $n\Box(\Phi)$.

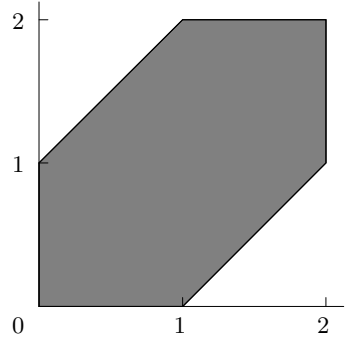
\begin{figure}[ht]
\begin{center}
\setlength{\unitlength}{0.00062500in}
\begingroup\makeatletter\ifx\SetFigFont\undefined%
\gdef\SetFigFont#1#2#3#4#5{%
  \reset@font\fontsize{#1}{#2pt}%
  \fontfamily{#3}\fontseries{#4}\fontshape{#5}%
  \selectfont}%
\fi\endgroup%
{\renewcommand{\dashlinestretch}{30}
\begin{picture}(2787,2802)(0,-10)
\path(225,1425)(1425,2625)(2625,2625)
	(2625,1425)(1425,225)
\texture{44555555 55aaaaaa aa555555 55aaaaaa aa555555 55aaaaaa aa555555 55aaaaaa 
	aa555555 55aaaaaa aa555555 55aaaaaa aa555555 55aaaaaa aa555555 55aaaaaa 
	aa555555 55aaaaaa aa555555 55aaaaaa aa555555 55aaaaaa aa555555 55aaaaaa 
	aa555555 55aaaaaa aa555555 55aaaaaa aa555555 55aaaaaa aa555555 55aaaaaa }
\shade\path(225,225)(225,1425)(1425,2625)
	(2625,2625)(2625,1425)(1425,225)(225,225)
\path(225,225)(225,1425)(1425,2625)
	(2625,2625)(2625,1425)(1425,225)(225,225)
\path(1425,225)(1425,300)
\path(2625,225)(2625,300)
\path(225,1425)(300,1425)
\path(225,2625)(300,2625)
\path(225,225)(225,2775)
\path(225,225)(2775,225)
\put(0,1350){\makebox(0,0)[lb]{\smash{{{\SetFigFont{9}{10.8}{\rmdefault}{\mddefault}{\updefault}$1$}}}}}
\put(0,2550){\makebox(0,0)[lb]{\smash{{{\SetFigFont{9}{10.8}{\rmdefault}{\mddefault}{\updefault}$2$}}}}}
\put(1350,0){\makebox(0,0)[lb]{\smash{{{\SetFigFont{9}{10.8}{\rmdefault}{\mddefault}{\updefault}$1$}}}}}
\put(2550,0){\makebox(0,0)[lb]{\smash{{{\SetFigFont{9}{10.8}{\rmdefault}{\mddefault}{\updefault}$2$}}}}}
\put(0,0){\makebox(0,0)[lb]{\smash{{{\SetFigFont{9}{10.8}{\rmdefault}{\mddefault}{\updefault}$0$}}}}}
\end{picture}
}
\caption{The polytope $\Box(\Phi)$}
\end{center}
\end{figure}

The set $S_{1,n}=\c_1-\Box(\Phi_n)$ is the interior of the
polyhedron determined by the inequalities
$$a_2\geq -2n,\;a_1\geq -2n,\; a_1-a_2\geq -n.$$
The partition function
$\iota_{\Phi_n}(\lambda)$ is given by a polynomial formula
$\iota[\c_1,\Phi_n](\lambda)$ when $\lambda$ varies in the set
$S_{1,n}\cap \Z^2$.

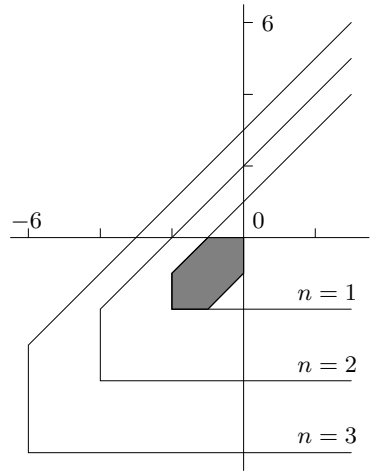
\begin{figure}[ht]
\begin{center}

\setlength{\unitlength}{0.00062500in}
\begingroup\makeatletter\ifx\SetFigFont\undefined%
\gdef\SetFigFont#1#2#3#4#5{%
  \reset@font\fontsize{#1}{#2pt}%
  \fontfamily{#3}\fontseries{#4}\fontshape{#5}%
  \selectfont}%
\fi\endgroup%
{\renewcommand{\dashlinestretch}{30}
\begin{picture}(3024,3939)(0,-10)
\path(12,1962)(3012,1962)
\path(1962,3912)(1962,12)
\path(2862,1362)(1362,1362)(1362,1662)(2862,3162)
\path(2862,762)(762,762)(762,1362)(2862,3462)
\path(2862,162)(162,162)(162,1062)(2862,3762)
\path(162,1962)(162,2037)
\path(762,1962)(762,2037)
\path(1362,1962)(1362,2037)
\path(2562,1962)(2562,2037)
\path(1962,3762)(2037,3762)
\path(1962,3162)(2037,3162)
\path(1962,2562)(2037,2562)
\texture{44555555 55aaaaaa aa555555 55aaaaaa aa555555 55aaaaaa aa555555 55aaaaaa 
	aa555555 55aaaaaa aa555555 55aaaaaa aa555555 55aaaaaa aa555555 55aaaaaa 
	aa555555 55aaaaaa aa555555 55aaaaaa aa555555 55aaaaaa aa555555 55aaaaaa 
	aa555555 55aaaaaa aa555555 55aaaaaa aa555555 55aaaaaa aa555555 55aaaaaa }
\shade\path(1362,1362)(1662,1362)(1962,1662)
	(1962,1962)(1662,1962)(1362,1662)(1362,1362)
\path(1362,1362)(1662,1362)(1962,1662)
	(1962,1962)(1662,1962)(1362,1662)(1362,1362)
\put(2412,237){\makebox(0,0)[lb]{\smash{{{\SetFigFont{9}{10.8}{\rmdefault}{\mddefault}{\updefault}$n=3$}}}}}
\put(2412,837){\makebox(0,0)[lb]{\smash{{{\SetFigFont{9}{10.8}{\rmdefault}{\mddefault}{\updefault}$n=2$}}}}}
\put(2412,1437){\makebox(0,0)[lb]{\smash{{{\SetFigFont{9}{10.8}{\rmdefault}{\mddefault}{\updefault}$n=1$}}}}}
\put(12,2037){\makebox(0,0)[lb]{\smash{{{\SetFigFont{9}{10.8}{\rmdefault}{\mddefault}{\updefault}$-6$}}}}}
\put(2112,3687){\makebox(0,0)[lb]{\smash{{{\SetFigFont{9}{10.8}{\rmdefault}{\mddefault}{\updefault}$6$}}}}}
\put(2037,2037){\makebox(0,0)[lb]{\smash{{{\SetFigFont{9}{10.8}{\rmdefault}{\mddefault}{\updefault}$0$}}}}}
\end{picture}
}
\caption{The polyhedron $S_{1,n}$}
\end{center}
\end{figure}

The set $S_{2,n}=\c_2-\Box(\Phi_n)$ is the interior of the
polyhedron determined by the inequalities
$$a_1\geq -2n,\;a_2\geq-2n,\;a_2-a_1\geq -n.$$ The  partition function
$\iota_{\Phi_n}(\lambda)$ is given by a polynomial function
$\iota[\c_2,\Phi_n](\lambda) $ when $\lambda$ varies in the set
$S_{2,n}\cap \Z^2$.

\begin{figure}[ht]
\begin{center}

\setlength{\unitlength}{0.00062500in}
\begingroup\makeatletter\ifx\SetFigFont\undefined%
\gdef\SetFigFont#1#2#3#4#5{%
  \reset@font\fontsize{#1}{#2pt}%
  \fontfamily{#3}\fontseries{#4}\fontshape{#5}%
  \selectfont}%
\fi\endgroup%
{\renewcommand{\dashlinestretch}{30}
\begin{picture}(3537,3339)(0,-10)
\path(2175,12)(2175,3312)
\path(3525,1962)(225,1962)
\path(3375,2862)(1875,1362)(1575,1362)(1575,2862)
\path(3375,2562)(1575,762)(975,762)(975,2862)
\path(3375,2262)(1275,162)(375,162)(375,2862)
\path(2775,1962)(2775,2037)
\path(3375,1962)(3375,2037)
\path(2175,1362)(2250,1362)
\path(2175,762)(2250,762)
\path(2175,162)(2250,162)
\path(2175,2562)(2250,2562)
\texture{44555555 55aaaaaa aa555555 55aaaaaa aa555555 55aaaaaa aa555555 55aaaaaa 
	aa555555 55aaaaaa aa555555 55aaaaaa aa555555 55aaaaaa aa555555 55aaaaaa 
	aa555555 55aaaaaa aa555555 55aaaaaa aa555555 55aaaaaa aa555555 55aaaaaa 
	aa555555 55aaaaaa aa555555 55aaaaaa aa555555 55aaaaaa aa555555 55aaaaaa }
\shade\path(2175,1962)(1875,1962)(1575,1662)
	(1575,1362)(1875,1362)(2175,1662)(2175,1962)
\path(2175,1962)(1875,1962)(1575,1662)
	(1575,1362)(1875,1362)(2175,1662)(2175,1962)
\put(2250,87){\makebox(0,0)[lb]{\smash{{{\SetFigFont{9}{10.8}{\rmdefault}{\mddefault}{\updefault}$-6$}}}}}
\put(0,2037){\makebox(0,0)[lb]{\smash{{{\SetFigFont{9}{10.8}{\rmdefault}{\mddefault}{\updefault}$-6$}}}}}
\put(2250,2037){\makebox(0,0)[lb]{\smash{{{\SetFigFont{9}{10.8}{\rmdefault}{\mddefault}{\updefault}$0$}}}}}
\put(150,2937){\makebox(0,0)[lb]{\smash{{{\SetFigFont{9}{10.8}{\rmdefault}{\mddefault}{\updefault}$n=3$}}}}}
\put(750,2937){\makebox(0,0)[lb]{\smash{{{\SetFigFont{9}{10.8}{\rmdefault}{\mddefault}{\updefault}$n=2$}}}}}
\put(1350,2937){\makebox(0,0)[lb]{\smash{{{\SetFigFont{9}{10.8}{\rmdefault}{\mddefault}{\updefault}$n=1$}}}}}
\put(2250,2487){\makebox(0,0)[lb]{\smash{{{\SetFigFont{9}{10.8}{\rmdefault}{\mddefault}{\updefault}$2$}}}}}
\end{picture}
}
\caption{The polyhedron $S_{2,n}$}
\end{center}
\end{figure}
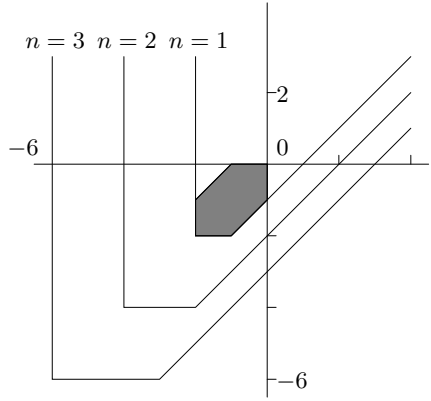

We see that the set $\cnull\cap S_{1,n}$ contains the $(2n-1)$
half-lines $p_{j}+t e^*_1$, where $t\geq 0$ and
\begin{equation}
p_j =
\begin{cases}
(1-n-j)e^*_1-je^*_2,&\text { if }1\leq j\leq n,\\
  (1-2n)e^*_1-je_2^*,&\text { if }n\leq j\leq 2n-1.
\end{cases}
\end{equation}

The function $\iota[\c_1,\Phi_n]$ vanishes on all the integral
points contained in  $\cnull\cap S_{1,n}$, as the partition
function $\iota_{\Phi_n}$ is identically $0$ on $\cnull$ . The set
of integral points in these half-lines is Zariski dense in the
affine line $a_2+j=0$, so that the polynomial function
$\iota[\c_1,\Phi_n]$ is divisible by $(a_2+1)(a_2+2)\cdots
(a_2+(2n-1))$.  Similarly the polynomial function
$\iota[\c_2,\Phi_n]$ is divisible by $(a_1+1)(a_1+2)\cdots
(a_1+(2n-1))$. These divisibility properties are also clear from
the Ehrhart reciprocity formula.

The set $S_{1,n}\cap S_{2,n}$ on which both formulae
$\iota[\c_1,\Phi_n]$ and $\iota[\c_2,\Phi_n]$ agree contains the  half
lines $q_k+t(e^*_1+e^*_2)$ with $t\geq 0$, where
$$
q_k =
\begin{cases}
(1-2n-j)e^*_1+(1-2n)e^*_2,&\text { if }1-n\leq j\leq 0,\\
  (1-2n)e^*_1+(1-2n+j)e_2^*,&\text { if }0\leq j\leq n-1.
\end{cases}
$$
By the same density argument, we see that the polynomial function
$\iota[\c_1,\Phi_n]-\iota[\c_2,\Phi_n]$ is divisible by
$$(a_1-a_2-(n-1))\cdots (a_1-a_2-1)(a_1-a_2)(a_1-a_2+1)\cdots
(a_1-a_2+(n-1))
.$$

\begin{figure}[ht]
\begin{center}

\setlength{\unitlength}{0.00062500in}
\begingroup\makeatletter\ifx\SetFigFont\undefined%
\gdef\SetFigFont#1#2#3#4#5{%
  \reset@font\fontsize{#1}{#2pt}%
  \fontfamily{#3}\fontseries{#4}\fontshape{#5}%
  \selectfont}%
\fi\endgroup%
{\renewcommand{\dashlinestretch}{30}
\begin{picture}(3924,4239)(0,-10)
\path(1962,4212)(1962,12)
\path(12,1962)(3912,1962)
\path(3162,2862)(1662,1362)(1362,1362)
	(1362,1662)(3162,3462)
\path(3162,2562)(1362,762)(762,762)
	(762,1362)(3162,3762)
\path(3162,2262)(1062,162)(162,162)
	(162,1062)(3162,4062)
\path(1962,2562)(2037,2562)
\path(1962,3162)(2037,3162)
\path(1962,3762)(2037,3762)
\path(1962,1362)(2037,1362)
\path(1962,762)(2037,762)
\path(1962,162)(2037,162)
\path(1362,1962)(1362,2037)
\path(762,1962)(762,2037)
\path(162,1962)(162,2037)
\path(2562,1962)(2562,2037)
\path(3162,1962)(3162,2037)
\path(3762,1962)(3762,2037)
\texture{44555555 55aaaaaa aa555555 55aaaaaa aa555555 55aaaaaa aa555555 55aaaaaa 
	aa555555 55aaaaaa aa555555 55aaaaaa aa555555 55aaaaaa aa555555 55aaaaaa 
	aa555555 55aaaaaa aa555555 55aaaaaa aa555555 55aaaaaa aa555555 55aaaaaa 
	aa555555 55aaaaaa aa555555 55aaaaaa aa555555 55aaaaaa aa555555 55aaaaaa }
\shade\path(1662,1962)(1362,1662)(1362,1362)
	(1662,1362)(1962,1662)(1962,1962)(1662,1962)
\path(1662,1962)(1362,1662)(1362,1362)
	(1662,1362)(1962,1662)(1962,1962)(1662,1962)
\put(3237,2187){\makebox(0,0)[lb]{\smash{{{\SetFigFont{9}{10.8}{\rmdefault}{\mddefault}{\updefault}$n=3$}}}}}
\put(3237,2487){\makebox(0,0)[lb]{\smash{{{\SetFigFont{9}{10.8}{\rmdefault}{\mddefault}{\updefault}$n=2$}}}}}
\put(3237,2787){\makebox(0,0)[lb]{\smash{{{\SetFigFont{9}{10.8}{\rmdefault}{\mddefault}{\updefault}$n=1$}}}}}
\put(12,2037){\makebox(0,0)[lb]{\smash{{{\SetFigFont{9}{10.8}{\rmdefault}{\mddefault}{\updefault}$-6$}}}}}
\put(2112,87){\makebox(0,0)[lb]{\smash{{{\SetFigFont{9}{10.8}{\rmdefault}{\mddefault}{\updefault}$-6$}}}}}
\put(2037,2037){\makebox(0,0)[lb]{\smash{{{\SetFigFont{9}{10.8}{\rmdefault}{\mddefault}{\updefault}$0$}}}}}
\put(2112,3687){\makebox(0,0)[lb]{\smash{{{\SetFigFont{9}{10.8}{\rmdefault}{\mddefault}{\updefault}$6$}}}}}
\end{picture}
}
\caption{The polyhedron  $S_{1,n}\cap S_{2,n}$}
\end{center}
\end{figure}
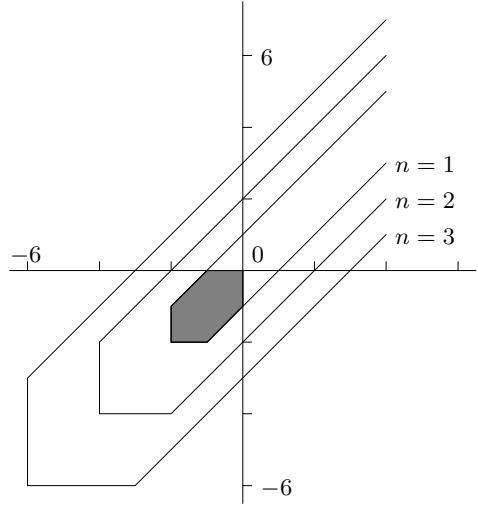

Below we give the formulas for the  cases $n=1$, $2$, $3$; the
functions $\iota[\c_1,\Phi_n]$ and $\iota[\c_2,\Phi_n]$ can easily
be computed from our formula,  with some help from Maple. One can
easily see that the appropriate functions vanish on the lines
indicated above. To simplify our formulas we use binomial
coefficients. Note that $\binom{a+m}k$, where we consider $a$ to
be the variable, vanishes at $a=-m,\ldots,k-m-1$.

\medskip

{\bf Case $n=1$}.
$$\iota[\c_1,\Phi_1]=(a_2+1),\quad \iota[\c_2,\Phi_1]=(a_1+1).$$
We also have
$$\iota[\c_1,\Phi_1]-\iota[\c_2,\Phi_2]=(a_2-a_1),$$
which vanishes on the line  $a_2=a_1$.

\medskip

{\bf Case $n=2$}.
$$\iota[\c_1,\Phi_2]=\frac{1}{2}\binom{a_2+3}3(2a_1-a_2+2),\quad
\iota[\c_2,\Phi_2]=\frac{1}{2}\binom{a_1+3}3(2a_2-a_1+2).$$
Again, we see that  the function
$$\iota[\c_1,\Phi_2]-\iota[\c_2,\Phi_2]=
 \frac {1}{2} \binom{a_1 - a_2 + 1}{3}(a_1 + a_2 + 4)
$$
 vanishes on the lines $a_1-a_2=-1$, $0$, $1$.

The example of $n=3$ is described in the introduction.

\subsection{A non-unimodular example}
Keeping the same vector space and lattice, we now consider a
non-unimodular system
$$\Phi=\{e^*_1, e^*_2, e^*_1+2e^*_2 \}.$$
The closed cone $C(\Phi)$
generated by $\Phi$ is still the first quadrant $\{a_1\geq 0, a_2\geq
0\}$.

Again, there are three open chambers for the system $\Phi$: The
exterior of the cone $C(\Phi)$ denoted by $\cnull$,
  the chamber $\c^1=\{a_2>0,2a_1>a_2\}$ and the chamber
$\c^2=\{a_1>0, a_2>2a_1\}$.

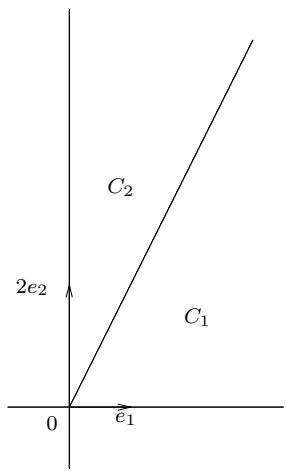
\begin{figure}[ht]
\begin{center}
\setlength{\unitlength}{0.00054167in}
\begingroup\makeatletter\ifx\SetFigFont\undefined%
\gdef\SetFigFont#1#2#3#4#5{%
  \reset@font\fontsize{#1}{#2pt}%
  \fontfamily{#3}\fontseries{#4}\fontshape{#5}%
  \selectfont}%
\fi\endgroup%
{\renewcommand{\dashlinestretch}{30}
\begin{picture}(2724,4539)(0,-10)
\texture{44555555 55aaaaaa aa555555 55aaaaaa aa555555 55aaaaaa aa555555 55aaaaaa 
    aa555555 55aaaaaa aa555555 55aaaaaa aa555555 55aaaaaa aa555555 55aaaaaa 
    aa555555 55aaaaaa aa555555 55aaaaaa aa555555 55aaaaaa aa555555 55aaaaaa 
    aa555555 55aaaaaa aa555555 55aaaaaa aa555555 55aaaaaa aa555555 55aaaaaa }
\path(612,12)(612,4512)
\path(612,12)(612,4512)
\path(12,612)(2712,612)
\path(12,612)(2712,612)
\path(612,612)(2412,4212)
\path(612,612)(2412,4212)
\path(612,612)(1212,612)
\path(1092.000,582.000)(1212.000,612.000)(1092.000,642.000)
\path(612,612)(612,1812)
\path(642.000,1692.000)(612.000,1812.000)(582.000,1692.000)
\put(387,387){\makebox(0,0)[lb]{\smash{{{\SetFigFont{8}{9.6}{\rmdefault}{\mddefault}{\updefault}$0$}}}}}
\put(1737,1437){\makebox(0,0)[lb]{\smash{{{\SetFigFont{8}{9.6}{\rmdefault}{\mddefault}{\updefault}$C_1$}}}}}
\put(987,2712){\makebox(0,0)[lb]{\smash{{{\SetFigFont{8}{9.6}{\rmdefault}{\mddefault}{\updefault}$C_2$}}}}}
\put(1062,462){\makebox(0,0)[lb]{\smash{{{\SetFigFont{8}{9.6}{\rmdefault}{\mddefault}{\updefault}$e_1$}}}}}
\put(87,1737){\makebox(0,0)[lb]{\smash{{{\SetFigFont{8}{9.6}{\rmdefault}{\mddefault}{\updefault}$2e_2$}}}}}
\end{picture}
}
\caption{The chambers}
\end{center}
\end{figure}

 The set $RG(\Phi,\c_1,\Gamma)$ consists of the elements
 $\{(0,0),(0,i\pi)\}$. The set $RG(\Phi,\c_2,\Gamma)$ is reduced to the
 element $\{(0,0)\}$.

 The polytope $\Box(\Phi)$ is the convex hull of the six points
 $0$, $e^*_1$, $e^*_2$, $2e^*_1+2e^*_2$, $3e^*_2+e^*_1$, $2e^*_1+3e^*_2$.

\begin{figure}[ht]
\begin{center}
\setlength{\unitlength}{0.00083333in}
\begingroup\makeatletter\ifx\SetFigFont\undefined%
\gdef\SetFigFont#1#2#3#4#5{%
  \reset@font\fontsize{#1}{#2pt}%
  \fontfamily{#3}\fontseries{#4}\fontshape{#5}%
  \selectfont}%
\fi\endgroup%
{\renewcommand{\dashlinestretch}{30}
\begin{picture}(2826,4002)(0,-10)
\texture{44555555 55aaaaaa aa555555 55aaaaaa aa555555 55aaaaaa aa555555 55aaaaaa 
	aa555555 55aaaaaa aa555555 55aaaaaa aa555555 55aaaaaa aa555555 55aaaaaa 
	aa555555 55aaaaaa aa555555 55aaaaaa aa555555 55aaaaaa aa555555 55aaaaaa 
	aa555555 55aaaaaa aa555555 55aaaaaa aa555555 55aaaaaa aa555555 55aaaaaa }
\shade\path(225,1425)(1425,3825)(2625,3825)
	(2625,2625)(1425,225)(225,225)(225,1425)
\path(225,1425)(1425,3825)(2625,3825)
	(2625,2625)(1425,225)(225,225)(225,1425)
\path(225,225)(2775,225)
\path(225,225)(2775,225)
\path(225,225)(225,3975)
\path(225,225)(225,3975)
\path(1425,225)(1425,300)
\path(2625,225)(2625,300)
\path(225,1425)(300,1425)
\path(225,2625)(300,2625)
\path(225,3825)(300,3825)
\put(2550,0){\makebox(0,0)[lb]{\smash{{{\SetFigFont{12}{14.4}{\rmdefault}{\mddefault}{\updefault}$2$}}}}}
\put(1350,0){\makebox(0,0)[lb]{\smash{{{\SetFigFont{12}{14.4}{\rmdefault}{\mddefault}{\updefault}$1$}}}}}
\put(0,1425){\makebox(0,0)[lb]{\smash{{{\SetFigFont{12}{14.4}{\rmdefault}{\mddefault}{\updefault}$1$}}}}}
\put(0,0){\makebox(0,0)[lb]{\smash{{{\SetFigFont{12}{14.4}{\rmdefault}{\mddefault}{\updefault}$0$}}}}}
\put(0,2550){\makebox(0,0)[lb]{\smash{{{\SetFigFont{12}{14.4}{\rmdefault}{\mddefault}{\updefault}$2$}}}}}
\put(0,3750){\makebox(0,0)[lb]{\smash{{{\SetFigFont{12}{14.4}{\rmdefault}{\mddefault}{\updefault}$3$}}}}}
\end{picture}
}
\caption{The polytope $\Box(\Phi)$}
\end{center}
\end{figure}
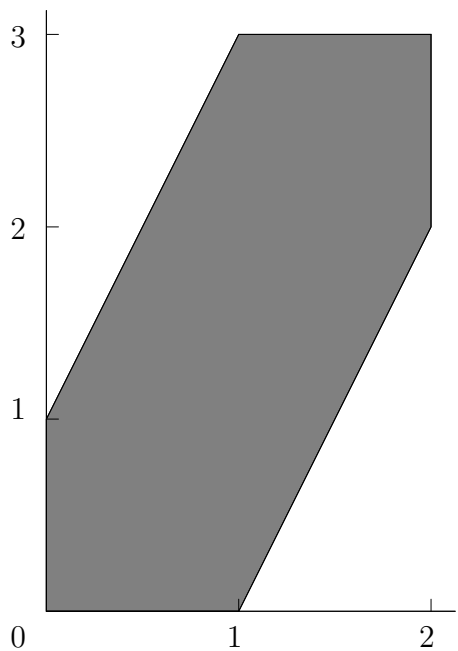

 We consider the system $\Phi_n$, where each of the three vectors
 $e^*_1$, $e^*_2$ and $e^*_1+2e^*_2$ has multiplicity $n$.

The set $S_{1,n}=\c_1-\Box(\Phi_n)$ is the interior of the
polyhedron determined by the inequalities
$$a_2\geq -3n,\; a_1\geq-2n,\; 2a_1-a_2\geq -2n.$$
The  partition function
$\iota_{\Phi_n}(\lambda)$ is given by a  periodic-polynomial formula
$\iota[\c_1,\Phi_n](\lambda) $ when $\lambda$ varies in the set
$S_{1,n}\cap \Z^2$. By our results, the periodic-polynomial
$\iota[\c_1,\Phi_n]$ is of the form $P(a_1,a_2)+\exp(i\pi a_2)
Q(a_1,a_2)$ where $P(a_1,a_2)$ and $Q(a_1,a_2)$ are polynomials.
 We denote by
$\iota[\c_1,\Phi_n,\even](a_1,a_2)$ the polynomial function
$P(a_1,a_2)+Q(a_1,a_2)$ on $\R^2$, which is equal to
$\iota[\c_1,\Phi_n](a_1,a_2)$ when $a_2$ is an even integer, and by
$\iota[\c_1,\Phi_n,\odd](a_1,a_2)$ the polynomial function
$P(a_1,a_2)-Q(a_1,a_2)$ on $\R^2$, which is equal to
$\iota[\c_1,\Phi_n](a_1,a_2)$ when $a_2$ is an odd integer.

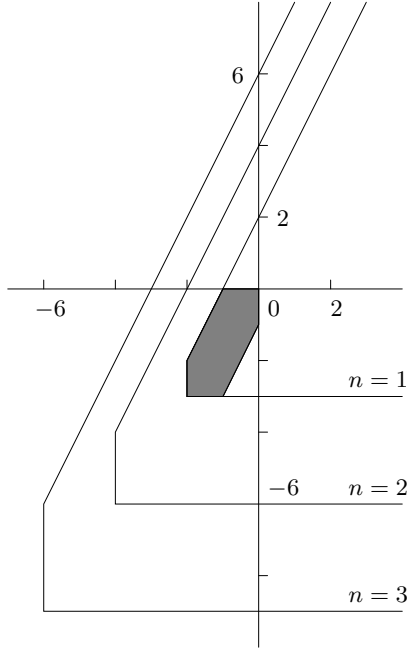
\begin{figure}[ht]
\begin{center}

\setlength{\unitlength}{0.00062500in}
\begingroup\makeatletter\ifx\SetFigFont\undefined%
\gdef\SetFigFont#1#2#3#4#5{%
  \reset@font\fontsize{#1}{#2pt}%
  \fontfamily{#3}\fontseries{#4}\fontshape{#5}%
  \selectfont}%
\fi\endgroup%
{\renewcommand{\dashlinestretch}{30}
\begin{picture}(3338,5439)(0,-10)
\path(12,3012)(3312,3012)
\path(3312,2112)(1512,2112)(1512,2412)(3012,5412)
\path(3312,1212)(912,1212)(912,1812)(2712,5412)
\path(3312,312)(312,312)(312,1212)(2412,5412)
\path(2112,12)(2112,5412)
\path(312,3087)(312,3012)
\path(912,3012)(912,3087)
\path(2712,3012)(2712,3087)
\path(1512,3012)(1512,3087)
\path(2112,2412)(2187,2412)
\path(2112,3612)(2187,3612)
\path(2112,4212)(2187,4212)
\path(2112,4812)(2187,4812)
\path(2112,1812)(2187,1812)
\path(2112,612)(2187,612)
\texture{44555555 55aaaaaa aa555555 55aaaaaa aa555555 55aaaaaa aa555555 55aaaaaa 
	aa555555 55aaaaaa aa555555 55aaaaaa aa555555 55aaaaaa aa555555 55aaaaaa 
	aa555555 55aaaaaa aa555555 55aaaaaa aa555555 55aaaaaa aa555555 55aaaaaa 
	aa555555 55aaaaaa aa555555 55aaaaaa aa555555 55aaaaaa aa555555 55aaaaaa }
\shade\path(2112,3012)(1812,3012)(1512,2412)
	(1512,2112)(1812,2112)(2112,2712)(2112,3012)
\path(2112,3012)(1812,3012)(1512,2412)
	(1512,2112)(1812,2112)(2112,2712)(2112,3012)
\put(2862,387){\makebox(0,0)[lb]{\smash{{{\SetFigFont{9}{10.8}{\rmdefault}{\mddefault}{\updefault}$n=3$}}}}}
\put(2862,1287){\makebox(0,0)[lb]{\smash{{{\SetFigFont{9}{10.8}{\rmdefault}{\mddefault}{\updefault}$n=2$}}}}}
\put(2862,2187){\makebox(0,0)[lb]{\smash{{{\SetFigFont{9}{10.8}{\rmdefault}{\mddefault}{\updefault}$n=1$}}}}}
\put(2187,2787){\makebox(0,0)[lb]{\smash{{{\SetFigFont{9}{10.8}{\rmdefault}{\mddefault}{\updefault}$0$}}}}}
\put(2712,2787){\makebox(0,0)[lb]{\smash{{{\SetFigFont{9}{10.8}{\rmdefault}{\mddefault}{\updefault}$2$}}}}}
\put(2262,3537){\makebox(0,0)[lb]{\smash{{{\SetFigFont{9}{10.8}{\rmdefault}{\mddefault}{\updefault}$2$}}}}}
\put(237,2787){\makebox(0,0)[lb]{\smash{{{\SetFigFont{9}{10.8}{\rmdefault}{\mddefault}{\updefault}$-6$}}}}}
\put(2187,1287){\makebox(0,0)[lb]{\smash{{{\SetFigFont{9}{10.8}{\rmdefault}{\mddefault}{\updefault}$-6$}}}}}
\put(1887,4737){\makebox(0,0)[lb]{\smash{{{\SetFigFont{9}{10.8}{\rmdefault}{\mddefault}{\updefault}$6$}}}}}
\end{picture}
}
\caption{The polyhedron $S_{1,n}$}
\end{center}
\end{figure}

The set $S_{2,n}=\c_2-\Box(\Phi_n)$ is the interior of the
polyhedron determined by the inequalities
$$a_1\geq -2n,\; a_2\geq-3n,\;2a_1-a_2\leq n.$$ The  partition function
$\iota_{\Phi_n}(\lambda)$ is given by a polynomial formula
$\iota[\c_2,\Phi_n](\lambda) $ when $\lambda$ varies in the set
$S_{2,n}\cap \Z^2$.

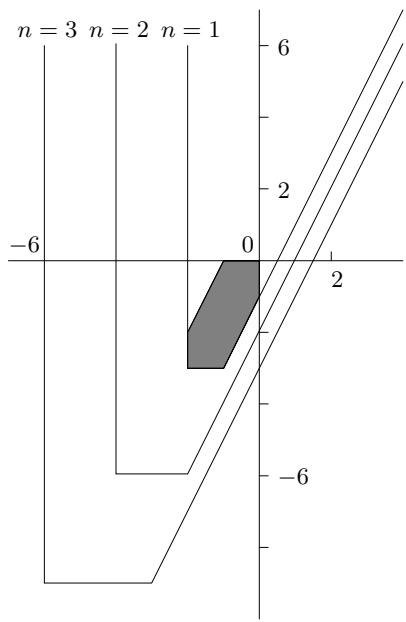
\begin{figure}[ht]
\begin{center}

\setlength{\unitlength}{0.00062500in}
\begingroup\makeatletter\ifx\SetFigFont\undefined%
\gdef\SetFigFont#1#2#3#4#5{%
  \reset@font\fontsize{#1}{#2pt}%
  \fontfamily{#3}\fontseries{#4}\fontshape{#5}%
  \selectfont}%
\fi\endgroup%
{\renewcommand{\dashlinestretch}{30}
\begin{picture}(3324,5139)(0,-10)
\path(312,4812)(312,312)(1212,312)(3312,4512)
\path(2112,2412)(2187,2412)
\path(2112,1812)(2187,1812)
\path(2112,1212)(2187,1212)
\path(2112,612)(2187,612)
\path(2112,3612)(2187,3612)
\path(2112,4212)(2187,4212)
\path(2112,4812)(2187,4812)
\path(2712,3087)(2712,3012)
\path(12,3012)(3312,3012)
\path(2112,5112)(2112,12)
\path(1512,4812)(1512,2112)(1812,2112)(3312,5112)
\path(912,4827)(912,1227)(1512,1227)(3312,4827)
\texture{44555555 55aaaaaa aa555555 55aaaaaa aa555555 55aaaaaa aa555555 55aaaaaa 
	aa555555 55aaaaaa aa555555 55aaaaaa aa555555 55aaaaaa aa555555 55aaaaaa 
	aa555555 55aaaaaa aa555555 55aaaaaa aa555555 55aaaaaa aa555555 55aaaaaa 
	aa555555 55aaaaaa aa555555 55aaaaaa aa555555 55aaaaaa aa555555 55aaaaaa }
\shade\path(2112,3012)(1812,3012)(1512,2412)
	(1512,2112)(1812,2112)(2112,2712)
	(2112,3012)(1812,3012)
\path(2112,3012)(1812,3012)(1512,2412)
	(1512,2112)(1812,2112)(2112,2712)
	(2112,3012)(1812,3012)
\put(687,4887){\makebox(0,0)[lb]{\smash{{{\SetFigFont{9}{10.8}{\rmdefault}{\mddefault}{\updefault}$n=2$}}}}}
\put(1287,4887){\makebox(0,0)[lb]{\smash{{{\SetFigFont{9}{10.8}{\rmdefault}{\mddefault}{\updefault}$n=1$}}}}}
\put(87,4887){\makebox(0,0)[lb]{\smash{{{\SetFigFont{9}{10.8}{\rmdefault}{\mddefault}{\updefault}$n=3$}}}}}
\put(2712,2787){\makebox(0,0)[lb]{\smash{{{\SetFigFont{9}{10.8}{\rmdefault}{\mddefault}{\updefault}$2$}}}}}
\put(2262,3537){\makebox(0,0)[lb]{\smash{{{\SetFigFont{9}{10.8}{\rmdefault}{\mddefault}{\updefault}$2$}}}}}
\put(1962,3087){\makebox(0,0)[lb]{\smash{{{\SetFigFont{9}{10.8}{\rmdefault}{\mddefault}{\updefault}$0$}}}}}
\put(2262,4737){\makebox(0,0)[lb]{\smash{{{\SetFigFont{9}{10.8}{\rmdefault}{\mddefault}{\updefault}$6$}}}}}
\put(2262,1137){\makebox(0,0)[lb]{\smash{{{\SetFigFont{9}{10.8}{\rmdefault}{\mddefault}{\updefault}$-6$}}}}}
\put(12,3087){\makebox(0,0)[lb]{\smash{{{\SetFigFont{9}{10.8}{\rmdefault}{\mddefault}{\updefault}$-6$}}}}}
\end{picture}
}
\caption{The polyhedron  $S_{2,n}$}
\end{center}
\end{figure}

 For the same reasons as before, the periodic-polynomial
 $\iota[\c_1,\Phi_n]$ vanishes on the lines $a_2=-1,-2,\ldots,-(3n-1)$,
 while the polynomial $\iota[\c_2,\Phi_n]$ vanishes on the lines
 $a_1=-1,-2,\ldots,-(2n-1)$; the function
 $\iota[\c_1,\Phi_n]-\iota[\c_2,\Phi_n]$ vanishes on the lines
 $(2a_1-a_2+k)=0$ for $-(n-1)\leq k\leq (2n-1)$. Note that if $k$ is
 even, then the $a_2$ coordinate of an integral point on the line
 $(2a_1-a_2+k)=0$ is even, while if $k$ is odd, then this coordinate
 is odd.

\begin{figure}[ht]
\begin{center}

\setlength{\unitlength}{0.00062500in}
\begingroup\makeatletter\ifx\SetFigFont\undefined%
\gdef\SetFigFont#1#2#3#4#5{%
  \reset@font\fontsize{#1}{#2pt}%
  \fontfamily{#3}\fontseries{#4}\fontshape{#5}%
  \selectfont}%
\fi\endgroup%
{\renewcommand{\dashlinestretch}{30}
\begin{picture}(3924,5439)(0,-10)
\path(12,3012)(3912,3012)
\path(2112,5412)(2112,12)
\path(3012,5412)(1512,2412)(1512,2112)
	(1812,2112)(3312,5112)
\path(2712,5412)(912,1812)(912,1212)
	(1512,1212)(3312,4812)
\path(2412,5412)(312,1212)(312,312)
	(1212,312)(3312,4512)
\path(2712,3012)(2712,3087)
\path(3312,3012)(3312,3087)
\path(1512,3012)(1512,3087)
\path(912,3012)(912,3087)
\path(312,3012)(312,3087)
\path(2112,2412)(2187,2412)
\path(2112,1812)(2187,1812)
\path(2112,1212)(2187,1212)
\path(2112,612)(2187,612)
\path(2112,3612)(2187,3612)
\path(2112,4212)(2187,4212)
\path(2112,4812)(2187,4812)
\texture{44555555 55aaaaaa aa555555 55aaaaaa aa555555 55aaaaaa aa555555 55aaaaaa 
	aa555555 55aaaaaa aa555555 55aaaaaa aa555555 55aaaaaa aa555555 55aaaaaa 
	aa555555 55aaaaaa aa555555 55aaaaaa aa555555 55aaaaaa aa555555 55aaaaaa 
	aa555555 55aaaaaa aa555555 55aaaaaa aa555555 55aaaaaa aa555555 55aaaaaa }
\shade\path(2112,3012)(1812,3012)(1512,2412)
	(1512,2112)(1812,2112)(2112,2712)
	(2112,3012)(1812,3012)
\path(2112,3012)(1812,3012)(1512,2412)
	(1512,2112)(1812,2112)(2112,2712)
	(2112,3012)(1812,3012)
\put(1962,3087){\makebox(0,0)[lb]{\smash{{{\SetFigFont{9}{10.8}{\rmdefault}{\mddefault}{\updefault}$0$}}}}}
\put(3387,4437){\makebox(0,0)[lb]{\smash{{{\SetFigFont{9}{10.8}{\rmdefault}{\mddefault}{\updefault}$n=3$}}}}}
\put(3387,4737){\makebox(0,0)[lb]{\smash{{{\SetFigFont{9}{10.8}{\rmdefault}{\mddefault}{\updefault}$n=2$}}}}}
\put(3387,5037){\makebox(0,0)[lb]{\smash{{{\SetFigFont{9}{10.8}{\rmdefault}{\mddefault}{\updefault}$n=1$}}}}}
\put(3312,2862){\makebox(0,0)[lb]{\smash{{{\SetFigFont{9}{10.8}{\rmdefault}{\mddefault}{\updefault}$4$}}}}}
\put(162,2862){\makebox(0,0)[lb]{\smash{{{\SetFigFont{9}{10.8}{\rmdefault}{\mddefault}{\updefault}$-6$}}}}}
\put(1887,4737){\makebox(0,0)[lb]{\smash{{{\SetFigFont{9}{10.8}{\rmdefault}{\mddefault}{\updefault}$6$}}}}}
\put(2262,537){\makebox(0,0)[lb]{\smash{{{\SetFigFont{9}{10.8}{\rmdefault}{\mddefault}{\updefault}$-8$}}}}}
\end{picture}
}
\caption{The polyhedron $S_{1,n}\cap S_{2,n}$}
\end{center}
\end{figure}
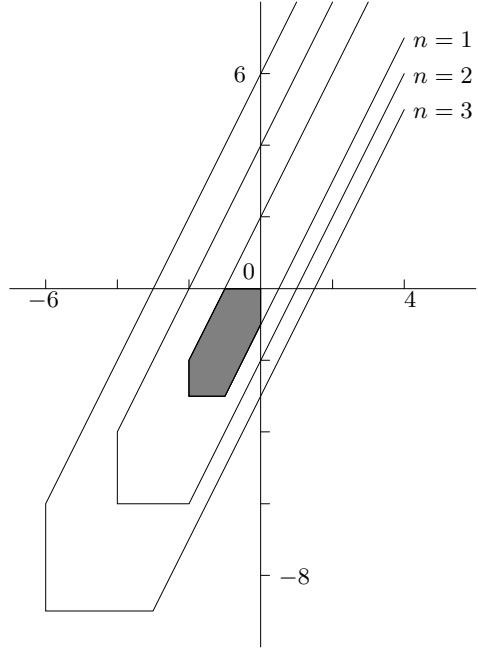

We verify these properties for $n=1$, $2$, $3$.\\

\medskip

{\bf Case $n=1$.}
$$\iota[\c_1,\Phi_1]= \frac {a_2}{2}
 + \frac {3}{4}  +  \frac {e^{i\pi a_2}}{4},
$$
hence
$$\iota[\c_1,\even]=\frac{1}{2}(a_2+2)\quad\text{and}\quad
\iota[\c_1,\odd]=\frac{1}{2}(a_2+1).$$
Thus the function $\iota[\c_1,\Phi_1]$ vanishes on the lines $a_2=-1$,
$a_2=-2$ as stated.

In the other chamber, we have $\iota[\c_2,\Phi_1]=(a_1+1)$. This function vanishes on the
line $a_1=-1$. Then
$$\iota[\c_1,\Phi_1,\even]-\iota[\c_2,\Phi_1]=
\frac{1}{2}(a_2-2a_1),$$ which vanishes when $2a_1-a_2=0$.
Also
$$
\iota[\c_1,\Phi_1,\odd]-\iota[\c_1,\Phi_1]=
-\frac{1}{2}(2a_1-a_2+1),$$ which vanishes when $2a_1-a_2+1=0$.

\medskip

{\bf Case $n=2$}

Here
$$\iota[\c_1,\Phi_2,\even]=
 \frac {1}{96} (a_2 + 2)(a_2 + 4)(4
a_1a_2 - a_2^{2} + 12a_1 + 2a_2 + 12)
$$
and
$$\iota[\c_1,\Phi_2,\odd]= \frac {1}{96} (a_2 + 1)(a_2 + 3)(a_2 + 5)(4a_1 - a_2 + 5).
$$
Thus the periodic-polynomial  function $\iota[\c_1,\Phi_2]$
vanishes on all the lines $a_2=-1$, $-2$, $-3$, $-4$, $-5$.

In the other chamber
$$\iota[\c_2,\Phi_2]=- \frac {1}{6} (a_1 + 1)(a_1 + 2)(a_1 + 3)(a_1 - a_2 - 1).$$
Thus the function $\iota[\c_2,\Phi_2]$ vanishes on all the lines
$a_2=-1$, $-2$, $-3$.

Now we have the polynomial formulas
$$\iota[\c_1,\Phi_2,\even]-\iota[\c_2,\Phi_2]=
 \frac {1}{96} (2a_1 - a_2)(2{a_{1}} - a_2 + 2)(4a_1^{2} - a_2^{2} + 16a_1
 - 6a_2 + 4)
$$
and
$$\iota[\c_1,\Phi_2,\odd]-\iota[\c_2,\Phi_2]
= \frac {1}{96} (2a_1 - a_2 - 1)(2a_1- a_2 + 1)(2a_1 - a_2 + 3)(2a_1 +
a_2 + 7),
$$
thus the function $\iota[\c_2,\Phi_2]-\iota[\c_1,\Phi_2]$ vanishes
on all the lines $2a_1-a_2=-3$, $-2$, $-1$, $0$, $1$.

\medskip

{\bf Case $n=3$.}

Here we have
\begin{multline*}
\iota[\c_1,\Phi_3,\even]=\frac {1}{53760} (a_2 + 2)(a_2
 + 4)(a_2 + 6)(a_2 + 8)  \\ \times
 (28a_1^2a_2 - 14a_1a_2^2 + 2a_2^3 + 70a_1^2 + 70a_1a_2 - 19a_2
^2 + 210a_1 + 44a_2 + 140)
\end{multline*}
and
\begin{multline*}
\iota[\c_1,\Phi_3,\odd]= \frac {1}{53760} (a_2 + 1)(a_2 + 3)(a_2 +
5)(a_2 + 7)  \\ \times
 (28a_1^2a_2 - 14a_1a_2^2 + 2a_2^3 +
  182a_1^2 + 14a_1a_2 - 11a_2^2 + 630a_1 - 52a_2 + 481),
\end{multline*}
thus the function $\iota[\c_1,\Phi_3]$ vanishes for $a_2=-1$,
$-2$, $-3$, $\ldots$, $-8$.

In the other chamber
\begin{multline*}
\iota[\c_2,\Phi_3]=  \frac {1}{1680}(a_1 + 1)(a_1 + 2)(a_1 +
3)(a_1 + 4)(a_1 + 5)  \\ \times (8a_1^2 - 14a_1a_2 +7a_2^2 - 15a_1
+ 21a_2 + 14),
\end{multline*}
thus the function $\iota[\c_2,\Phi_3]$ vanishes for $a_1=-1$,
$-2$, $-3$, $-4$, $-5$.

Now the difference $  \iota[\c_1,\Phi_3,\even]-\iota[\c_2,\Phi_3]$ is
given by
\begin{multline*}
-\frac 1{53760}(2a_1-a_2 - 2)(2a_1 - a_2 )(2a_1 - a_2+2)(2    a_1
- a_2 + 4) \times \\ \times
  (16a_1^3 + 4a_1^2a_2 - 2a_1a_2^2 - 2a_2^3 + 178a_1^2 +
  18a_1a_2 - 29a_2^2 + 598a_1 - 68a_2 + 484),
\end{multline*}
and $\iota[\c_1,\Phi_3,\odd]-\iota[\c_2,\Phi_3]$ by
\begin{multline*}
-\frac1{53760} (2a_1 - a_2 - 1)(2a_1 - a_2 + 1)(2a_1 - a_2 +
3)(2a_1 - a_2 + 5)  \\ \times (16a_1^3 + 4a_1^2a_2 - 2a_1a_2^2 -
2a_2^3 + 146a_1^2 - 6a_1a_2
 - 37a_2^2 + 298a_1 - 212 a_2 - 217).
\end{multline*}
Thus the function $\iota[\c_1,\Phi_3]-\iota[\c_2,\Phi_3]$ vanishes
on the lines $$2a_1-a_2=-5\;, -4\:, -3\;, -2\;, -1\;, 0\;, 1\;, 2.
$$

\end{document}